\documentclass{amsart}
\usepackage{amsmath,amssymb,amsthm,amscd,mathtools}
\usepackage[dvipsnames]{xcolor}
\usepackage{tikz-cd}
\usepackage{enumitem}
\usepackage{array,booktabs,calc}
\usepackage[textwidth=6.1in,textheight=8.9in,centering]{geometry}
\usepackage{microtype}
\usepackage{hyperref}
\hypersetup{
colorlinks=true,
pdftitle={Valuation Descent for Cyclic Groups: Complete RO-Graded Tambara Structure},
pdfauthor={Dongming Zhang},
linkcolor=blue,
anchorcolor=blue,
citecolor=blue
}
\usepackage[nameinlink]{cleveref}

\newtheorem{theorem}{Theorem}[section]
\newtheorem{corollary}[theorem]{Corollary}
\newtheorem{proposition}[theorem]{Proposition}
\newtheorem{lemma}[theorem]{Lemma}

\theoremstyle{definition}
\newtheorem{definition}[theorem]{Definition}
\newtheorem{example}[theorem]{Example}

\theoremstyle{remark}
\newtheorem{remark}[theorem]{Remark}

\makeatletter
\let\c@equation\c@theorem
\makeatother
\numberwithin{equation}{section}

\newcommand{\HZ}{H\underline{\mathbb Z}}
\newcommand{\res}{\operatorname{res}}
\newcommand{\Res}{\operatorname{Res}}
\newcommand{\tr}{\operatorname{tr}}
\newcommand{\Ind}{\operatorname{Ind}}
\newcommand{\sgn}{\operatorname{sgn}}
\newcommand{\coker}{\operatorname{coker}}
\newcommand{\Sp}{\operatorname{Sp}}
\newcommand{\CAlg}{\operatorname{CAlg}}
\newcommand{\holim}{\operatorname*{holim}}
\newcommand{\hocolim}{\operatorname*{hocolim}}

\newcommand{\leadbreak}{\par\ifdim\pagetotal>.75\textheight\newpage\fi}

\begin{document}
\author{Dongming Zhang}
\title[Valuation descent and graded norms]{Complete $RO$-Graded Tambara functor of $H\underline{\mathbb Z}$ for Cyclic Groups}
\begin{abstract}
We compute the $RO$-graded Tambara functor of
$H\underline{\mathbb{Z}}$ for every finite cyclic group.
We introduce a valuation language for homotopy groups to give the explicit formula for products, restrictions,
transfers, Weyl actions, and norms. Besides the computation, We also develop a theory that reconstructs a genuine equivariant spectrum as the homotopy limit
of a finite diagram of homotopy fixed-point spectra, Tate spectra, and only one localization spectrum
at the terminal of the diagram. We call it the Tate compression theory. From it we also give a filtration for $\pi_\star(H\underline{\mathbb{Z}})$, compatible with products, restrictions,
transfers, Weyl actions, and norms.
\end{abstract}
\maketitle
\setcounter{tocdepth}{1}
\tableofcontents
\section{Introduction}\label{sec:introduction}
The Eilenberg--Mac Lane spectrum $H\mathbb Z$ is one of the most basic spectra in homotopy theory. It represents
ordinary homology and cohomology. It connects homology theory with stable homotopy theory. For a space $X$,
$H_r(X;\mathbb Z)=\pi_r(H\mathbb Z\wedge\Sigma^\infty X_+)$, while $H^r(X;\mathbb Z)=[\Sigma^\infty X_+,\Sigma^rH\mathbb Z]$.
Suspensions of Eilenberg--Mac Lane spectra are the layers of
Postnikov towers. Its homotopy groups is simple, we have $\pi_0(H\mathbb Z)=\mathbb Z$ and
$\pi_r(H\mathbb Z)=0$ for $r\ne0$.

For a finite group $G$, let $\HZ$ be the genuine $G$-spectrum
associated to the constant Mackey functor $\underline{\mathbb Z}$.
For $K\leq L\leq G$,
$\underline{\mathbb Z}(G/K)=\mathbb Z$, with $\res_K^L=\operatorname{id}$ and $\tr_K^L=[L:K]$;
all conjugation maps are the identity. In integer degrees,
$\underline\pi_0^G(\HZ)=\underline{\mathbb Z}$, whereas $\underline\pi_r^G(\HZ)=0\;(r\ne0)$.
Write $RO(G)$ for the Grothendieck group of finite-dimensional real
$G$-representations. For $V\in RO(G)$, let $S^V$ be the
representation sphere spectrum.
The coefficients of integral Bredon cohomology are
$\pi_V^G(\HZ)=[S^V,\HZ]^G =H_G^{-V}(\mathrm{pt};\underline{\mathbb Z})$.

These groups depend on the representation, not only on its
dimension. For example, let $\alpha$ be the sign representation
of $C_2$. Then
$\pi_{-1}^{C_2}(\HZ)=0$, while $\pi_{-\alpha}^{C_2}(\HZ)\cong\mathbb Z/2$.
Both degrees have underlying dimension $-1$.
For larger cyclic groups, $RO(G)$ is more complicated. Surprisingly and not surprisingly, the $RO(G)$-degree homotopy groups are extremely difficult to compute.

For $K\leq L\leq G$ and $V\in RO(L)$, restriction and transfer
have the form
\[
\pi_V^L(\HZ) \xrightarrow{\ \res_K^L\ } \pi_{\Res_K^L V}^K(\HZ) \xrightarrow{\ \tr_K^L\ } \pi_V^L(\HZ).
\]
Together with products and Weyl actions, these maps give the graded
Green functor. There is also norm functor for $H\underline{\mathbb{Z}}$.
For $W\in RO(K)$, the norm has the form
\[
 N_K^L\colon\pi_W^K(\HZ)
       \longrightarrow\pi_{\Ind_K^L W}^L(\HZ),
\]
where $\Res$ and $\Ind$ denote restriction and induction of real
representations.
Norms preserve products and satisfy the Tambara
distributivity law. Even when $W\in\mathbb Z\subset RO(K)$,
the induced degree $\Ind_K^LW$ need not lie in
$\mathbb Z\subset RO(L)$.
These groups and operations form the graded Tambara functor.

Let $G$ be a finite cyclic group and put $H=\HZ$.
We compute the complete $RO(G)$-graded Tambara functor of $H$.

\begin{theorem}[Complete Tambara functor]
\label{thm:intro-green-completion}
\label{thm:intro-complete-cyclic}
\label{thm:intro-complete-tambara}
For every finite cyclic group $G$, the groups
$\pi_V^K(H)$, $K\leq G$, $V\in RO(K)$, and the associated graded Tambara
functor are given by the following results.
\begin{enumerate}[label=\textup{(\roman*)}]
\item For $G=C_{p^n}$, with $p$ prime and $n\geq1$, the homotopy groups are given by
Theorems~\ref{thm:prime-power-final} and
\ref{thm:integral-cyclic-reconstruction}.
Products are given by
Proposition~\ref{prop:array-products},
restrictions and transfers by
Proposition~\ref{prop:array-equal-parity} and
Theorem~\ref{thm:cross-maps},
and Weyl actions by Proposition~\ref{prop:adjacent-weyl}.
For arbitrary finite cyclic groups, products, restrictions,
transfers, and Weyl actions are given by
Proposition~\ref{prop:cyclic-green-operations}.
\item The norms $N_K^L$ for all $K\leq L\leq G$ are given by Theorem~\ref{thm:small-data-determination}.
\end{enumerate}
\end{theorem}
Example~\ref{ex:c4-c2-coefficients} shows that our computation is the same as the already computed Green functor of $C_2$ and $C_4$ case.

Our computation of the $RO(G)$-graded Tambara functor of $\HZ$
is relevant to several areas, including chromatic homotopy theory,
motivic homotopy theory, topological Hochschild homology,
and algebraic $K$-theory.

The $RO(G)$-graded homotopy groups of $\HZ$ provide the
$E_2$-page for the slice spectral sequences used by
Hill--Hopkins--Ravenel in their work on the Kervaire invariant
one problem.
For $G=C_{2^n}$, $n\geq1$, let $MU_{\mathbb R}$ be the Real bordism
spectrum. Hill--Hopkins--Ravenel form
$MU^{((G))}=N_{C_2}^G MU_{\mathbb R}$,
where $N_{C_2}^G$ is the norm functor on spectra.
For $G=C_8$, the fixed points of a suitable localization detect
any Kervaire invariant one element in dimension $2^{j+1}-2$
with $j\geq3$. Their vanishing and periodicity results rule out
these elements for $j\geq7$
\cite{HHRNonexistence}.
After localization at $2$, their Slice Theorem expresses the
nonzero slices of $MU^{((G))}$ as wedges of spectra of the form
\[
 G_+\wedge_K\bigl(S^{m\rho_K}\wedge\HZ_{(2)}\bigr).
\]

Our result also has a connection with chromatic
homotopy theory.
Let $E_h$ be the Lubin--Tate spectrum of a height-$h$ formal group
over a perfect field of characteristic $2$.
Let $G$ be a finite subgroup of its Morava stabilizer group.
Assume that $G$ contains the central $C_2$.
Hahn--Shi \cite{HahnShiRealOrientations} construct a $G$-equivariant map
\[
 N_{C_2}^G MU_{\mathbb R}\longrightarrow F(EG_+,E_h).
\]
The spectrum $F(EG_+,E_h)$ is the cofree genuine
$G$-spectrum for the stabilizer action. For $G=C_2$, this map is a Real
orientation. Taking $G$-fixed points gives maps to the homotopy fixed points $E_h^{hG}$.

Hill--Shi--Wang--Xu calculate the slice spectral sequence of
$BP^{((C_4))}\langle2\rangle$ \cite{HSWXHeightFour}.
This spectrum is a quotient of $N_{C_2}^{C_4}BP_{\mathbb R}$,
where $BP_{\mathbb R}$ is the Real Brown--Peterson spectrum.
They use the calculation to study the homotopy fixed points
of a height-$4$ Lubin--Tate spectrum with $C_4$-action.
The slices are built from suspensions
of $\HZ_{(2)}$. Their calculation uses norms to
relate differentials at different subgroups.

The spectrum $\HZ$ also appears in algebraic $K$-theory.
Let $k/\mathbb F_q$ be a finite extension of finite fields, and let
$G=\operatorname{Gal}(k/\mathbb F_q)$. This group is cyclic.
Write $K_G(k)$ for the genuine $G$-spectrum of Galois-equivariant
algebraic $K$-theory. Chan--Vogeli prove that, for every
$V\in RO(G)$,
\[
 \pi_V^G K_G(k)\cong
 \pi_V^G(\HZ)\oplus\pi_V^G\bigl(\tau_{\geq1}K_G(k)\bigr)
\]
\cite[Theorem~5.2.3]{ChanVogeli}.
Here $\tau_{\geq1}K_G(k)$ is the fiber of the truncation
$K_G(k)\to\HZ$.
The second summand is computed from the positive $K$-groups of $k$,
with Galois action twisted by $\det(V)\colon G\to\{\pm1\}$, and vanishes
for $\dim V\leq0$. Chan--Vogeli compute the products; the first summand
has the products of $\HZ$ \cite[Proposition~6.2.1]{ChanVogeli}. Thus, for $G=C_N$, the $RO(C_N)$-graded groups computed here
give the explicit formula of it.

Our calculations also give explicit formulas for the $RO(G)$-graded
homotopy groups of the classical genuine $S^1$-spectrum
$\operatorname{THH}(\mathbb F_p)$ for every finite cyclic group
$G=C_N$.
For a prime $p$, write $N=p^n m$ with $(m,p)=1$.
The restriction-tower description in \cite[Section~3]{BHKKS25} gives
\[
 \Phi^{C_p}(\HZ)
 \simeq
 \Phi^{C_p}(H\underline{\mathbb Z}_p)
 \simeq
 \operatorname{THH}(\mathbb F_p)|_{C_{p^n}},
\]
where the Eilenberg--Mac Lane spectra are taken over $C_{p^{n+1}}$.
For odd $p$, the integral comparison also appears in
\cite[Corollary~10.8]{AMGR}.
Let $V\in RO(C_N)$ and put
$V_d=(V^{C_d})|_{C_{p^n}}$ for $d\mid m$.
Let $q\colon C_{p^{n+1}}\twoheadrightarrow C_{p^n}$ be the quotient
map, and let $\lambda$ be the underlying real representation of
a faithful complex character of $C_{p^{n+1}}$.
Combining this comparison with the $p$-local idempotent splitting
\cite[Theorem~3.4 and Proposition~3.15]{Bohme19},
the cyclotomic equivalences
\cite[Theorems~4.6--4.7]{ABGHLM18},
and Euler localization shows that, for $r\in\mathbb Z_{\geq0}$ with
$2r>\max_{d\mid m}\dim V_d$,
\[
 \pi_V^{C_N}\operatorname{THH}(\mathbb F_p)
 \cong
 \prod_{d\mid m}
 \left(
   \pi_{q^*V_d-r\lambda}^{C_{p^{n+1}}}\HZ
 \right)^{C_{m/d}},
\]
where $g\in C_{m/d}\leq C_{N/d}$ acts on the $d$th factor by
$x\mapsto\det(g|_{V^{C_d}})\,x$.
The bound on $r$ follows from the Euler cofibration
$S(\lambda)_+\to S^0\to S^\lambda$.
In fact, forming the restriction tower in genuine commutative
ring spectra and applying
\cite[Theorems~5.2--5.3]{UllmanTambara} and
\cite[Theorem~E]{BlumbergMandell24}
makes the comparison compatible with all norms.
Together with Euler localization, the norm double-coset formula,
and the Adams operations of \cite[Section~10]{ABGHLM18},
this recovers all products, restrictions, transfers, Weyl actions,
and norms from the corresponding operations on $\HZ$,
and hence determines the complete $RO(G)$-graded Tambara functor
of $\operatorname{THH}(\mathbb F_p)$.

This computation also applies in motivic homotopy theory. For a finite group $G$, write $\mathbf{SH}^G(\mathbb C)$ and
$\mathbf{SH}^G$ for the equivariant motivic and genuine equivariant
stable homotopy categories, respectively. Complex Betti realization gives a functor
$\mathrm{Re}_B^G: \mathbf{SH}^G(\mathbb C)\longrightarrow\mathbf{SH}^G$.
Heller--Voineagu--{\O}stv{\ae}r prove that the equivariant
motivic integral Eilenberg--Mac Lane spectrum $\mathrm M\mathbb Z$
satisfies
\[
\mathrm{Re}_B^G(\mathrm M\mathbb Z)\simeq\HZ
\]
\cite[Theorem~A.31]{HellerVoineaguOstvaer}.
Thus, for $G=C_N$, the $RO(C_N)$-graded groups computed here
give the homotopy group of the topological realization of
equivariant motivic integral cohomology.
For $G=C_2$, a related realization also arises from complex
conjugation on real varieties
\cite[Theorem~4.17]{HellerOrmsby}.

The new language we introduce is the valuation language. This language enable us to write the homotopy tambara functor of $\HZ$ in closed formulas.

We also develop a theory that reconstructs
a genuine equivariant spectrum as the homotopy limit
of a finite diagram of homotopy fixed-point spectra,
Tate spectra, and only one localization spectrum at the terminal of the diagram. We call it the Tate compression theory.

Write $\Sp^G$ for the $\infty$-category of genuine $G$-spectra, $\Sp=\Sp^e$,
$\CAlg(\Sp^G)$ for its commutative algebras, and
$\operatorname{Mod}_R(\Sp^G)$ for $R$-modules.

\begin{definition}
Let a family $\mathcal F$ of subgroups of $G$ be closed under subgroups
and conjugation. Its universal $G$-space satisfies
\[
 (E\mathcal F)^K\simeq
 \begin{cases}*,&K\in\mathcal F,\\ \varnothing,&K\notin\mathcal F.\end{cases}
 \qquad
 \widetilde E\mathcal F=\operatorname{cofib}((E\mathcal F)_+\to S^0).
\]
For $X\in\Sp^G$, its family localization and completion are
$\widetilde E\mathcal F\wedge X$ and $F((E\mathcal F)_+,X)$.
We call $X$ $\widetilde E\mathcal F$-local or $\mathcal F$-complete
when the corresponding unit map is an equivalence.
\end{definition}

The motivation of Tate compression theorey is mainly the family Tate squares, coming from
\cite[Section~17]{GreenleesMayTate}.
Yan uses them for cyclic $2$-groups
\cite[Section~1.2]{Yan}.
Related Tate gluing appears in
\cite[Theorem~7.2]{AMGR} and \cite{QuigleyShah}.

Tate compression diagram is the core of Tate compression theory.
Let $G=C_{p^n}$, with $n\geq1$, and set
$\mathcal F_i=\{C_{p^j}:0\leq j\leq i\}$ for $0\leq i<n$.
Put $\mathcal F_{-1}=\varnothing$.
For $X\in\Sp^G$, put
\[
 \begin{aligned}
 X^{\Phi_i}&=\widetilde E\mathcal F_{i-1}\wedge X
       &&(0\leq i\leq n),\\
 X^{h_i}&=F((E\mathcal F_i)_+,X^{\Phi_i}),\qquad
 X^{t_i}=\widetilde E\mathcal F_i\wedge X^{h_i}
       &&(0\leq i<n).
 \end{aligned}
\]
We write the subgroup geometric-fixed-point functor as
$\Phi^{C_{p^i}}$.
Isotropy separation gives a homotopy pullback square
\[
 \begin{tikzcd}
 X^{\Phi_i}\ar[r]\ar[d]&X^{h_i}\ar[d]\\
 X^{\Phi_{i+1}}\ar[r]&X^{t_i}.
 \end{tikzcd}
\]
Under the condition \eqref{eq:uniform-family-nilpotence}
with $R=X\in\CAlg(\Sp^G)$, the maps
$X^{h_{i+1}}\to X^{t_i}$ for $0\leq i<n-1$, thus we construct a diagram
\[
X^{h_0}\to X^{t_0}\leftarrow X^{h_1}
 \to\cdots\to X^{t_{n-1}}\leftarrow X^{\Phi_n}
.
\]
Our Tate compression glues these squares and
eliminates all intermediate terms
$X^{\Phi_i}$ for $1\le i<n$.
The original spectrum is reconstructed from a finite
diagram containing only the homotopy-fixed-point
terms $X^{h_i}$, the Tate terms $X^{t_i}$,
and the terminal term $X^{\Phi_n}$.
This reconstruction
also applies to modules and commutative algebras.

\begin{theorem}[Tate compression]
\label{thm:intro-sharp-compression}
Let $X\in\CAlg(\Sp^{C_{p^n}})$.
Suppose $X^{t_i}$ is $\mathcal F_{i+1}$-complete for
$0\leq i<n-1$.
Then the maps $X^{\Phi_{i+1}}\to X^{t_i}$ factor through
$X^{\Phi_{i+1}}\to X^{h_{i+1}}$ in commutative algebras.
We have
\[
 X\simeq\operatorname*{holim}
 \bigl(
 X^{h_0}\to X^{t_0}\leftarrow X^{h_1}
 \to\cdots\to X^{t_{n-1}}\leftarrow X^{\Phi_n}
 \bigr)
\]
\textup{(Theorem~\ref{thm:sharp-path-compression})}.
Under the Euler nilpotence hypothesis
\eqref{eq:uniform-family-nilpotence},
the functor in \eqref{eq:module-category-descent} gives an equivalence
between $\operatorname{Mod}_X(\Sp^{C_{p^n}})$ and the diagram category
$\mathcal D_X$ of \eqref{eq:local-section-category}, with $R=X$.

For every nontrivial finite cyclic $G$, the map $H\to B_G$ of
\eqref{eq:cyclic-sylow-vertex-algebra} satisfies
$\exp_{B_G}(H)=2$
\textup{(Definition~\ref{def:b-exponent-amitsur})}.
\end{theorem}

For $G=C_{p^n}$ and $X=H$,
Lemma~\ref{lem:local-euler-height} verifies
\eqref{eq:uniform-family-nilpotence} with $R=H$.
Hence
\[
 H\simeq\holim\bigl(
 H^{h_0}\to H^{t_0}\leftarrow H^{h_1}
 \to\cdots\to H^{t_{n-1}}\leftarrow H^{\Phi_n}
 \bigr).
\]

All norms are determined by the following.

\begin{theorem}[Complete calculation of norms]
\label{thm:small-data-determination}
Let $G$ be a finite cyclic group and let $H=\HZ$.
For every $K\leq L\leq G$ and $V\in RO(K)$, we compute the norm
\[
N_K^L\colon\pi_V^K(H)
\longrightarrow\pi_{\Ind_K^L V}^L(H)
\]
on every class of homotopy group.
The formulas are as follows.
\begin{enumerate}[label=\textup{(\roman*)}]
\item The homotopy groups are given in Sections~\ref{sec:reconstruction} and
\ref{sec:sylow-assembly}.
\item For every prime $p$ and $m\in\mathbb Z$,
$N_e^{C_p}(m)=m^p$.
\item For $C_{p^k}<C_{p^{k+1}}\leq G$, $k\geq1$, the norms in
natural degrees is given by Theorem~\ref{thm:tambara-natural-kernel}.
\item For $C_{p^k}<C_{p^{k+1}}\leq G$, $k\geq1$, the norms in shifted degrees is given by
Theorem~\ref{thm:tambara-shifted-exact}.
\item Norms involve the classes $v_D$ is given by
Proposition~\ref{prop:norm-full-degree-coherence}, so this extend the norms from $RO_\mathrm{sel}(C_{p^k})$ to $RO(C_{p^k})$.
\end{enumerate}
Norm transitivity gives all norms $N_K^{C_{p^n}}$ for $K\leq C_{p^n}\leq G$.
Theorem~\ref{thm:prime-local-cyclic-norms}
gives the norms for arbitrary finite cyclic groups.
Equation~\eqref{eq:norm-explicit-exponential} gives norms of sums.
Thus these formulas compute $N_K^L(x)$ for every
$x\in\pi_V^K(H)$.
\end{theorem}

We also have the following result from the computation of norms.

\begin{theorem}[Vanishing of norms]
\label{thm:intro-first-order-normed-descent}
For every finite cyclic $G$, the descent filtration
$F^1\underline\pi_\star H$ defined in \eqref{eq:cyclic-descent-ideal}
is square-zero and is preserved by restrictions, transfers, Weyl
actions, and norms. For $K<L$ and $m\in F^1\pi_V^K(H)$,
\[
 [L:K]N_K^L(m)=0,\qquad
 [L:K]\notin\{2,3\}\ \Longrightarrow\ N_K^L(m)=0.
\]
Both exceptional indices occur, and every composite of two proper
norms on $F^1$ is zero. For $J\leq L$,
\[
 JK<L\ \Longrightarrow\ \res_J^L N_K^L(m)=0.
\]
\end{theorem}

The $RO(C_p)$-graded calculations began with Stong
(unpublished).
Lewis gave the first published account \cite{Lewis1988}.
Both used Burnside coefficients.
For any Mackey functor $\underline{M}$,
Ferland--Lewis \cite{FerlandLewis} computed the
$RO(C_p)$-graded homotopy Mackey functors of
$H\underline{M}$, including $\HZ$.

For $p$ odd, Zeng \cite{Zeng} computed the full
$RO(C_{p^2})$-graded Green functor of $\HZ$,
including the groups, $\operatorname{res}$,
$\operatorname{tr}$, and products.
For $C_4$, Zeng computed the orientable degrees.
Georgakopoulos \cite{Georgakopoulos} computed the full
$RO(C_4)$-graded Green functor of $\HZ$.
Yan \cite{Yan} gave another full $C_4$ computation
using generalized Tate squares.

For cyclic groups $G$ of odd square-free order,
Basu--Ghosh \cite{BasuGhoshSquarefree} computed the
$RO(G)$-graded homotopy Mackey functors of $\HZ$
and their ring structure.
For $p$ odd, Ayala--Mazel-Gee--Rozenblyum
\cite[Theorem~11.5]{AMGR} gave chain complexes
for the Picard-graded cohomology of $C_{p^n}$
with constant coefficients.
They expressed $\operatorname{res}$ and
$\operatorname{tr}$ as chain maps
\cite[Lemmas~11.7--11.9]{AMGR}.
The additive Mackey functors are the cohomology
of these complexes.
For $p$ odd and $n\geq 3$,
Basu--Ghosh \cite{BasuGhosh}
found nontrivial Mackey functor extensions.

Sikora--Yan \cite{SikoraYan} computed the positive
and negative cones of $\HZ$ for cyclic $2$-groups,
in degrees $k-W$ and $k+W$, respectively,
where $k\in\mathbb{Z}$ and $W$ is an actual representation.
Basu--Dey \cite[Theorem~B]{BasuDey} computed the
positive-cone ring for finite cyclic groups
in divisor grading.
For cyclic $p$-groups, they also computed the ring
in degrees $V$ satisfying
$\dim(V^K)\leq 0$ for all $K\leq G$
\cite[Section~8]{BasuDey}.
Dugger--Hazel \cite{DuggerHazelPoint} computed the
positive-cone ring and negative-cone Mackey functors
for cyclic groups of odd order.
They also studied products and mixed degrees.
They gave complete ring descriptions for
$C_p$ and $C_{p^2}$,
and left further computations in the irregular region
as an open problem.

Zeng \cite{Zeng} used Tate squares and induction
through quotient groups.
Dugger--Hazel \cite{DuggerHazelAlg} developed
homological algebra for $\underline{\mathbb{Z}}$-modules.
We also use graded Tambara functors,
norms under localization, power operations,
and tensor descent
\cite{AngeltveitBohmann,BlumbergHill,Carrick,Balderrama,MNN}.

In this article, Section~\ref{sec:sylow} gives the foundation for later assemble from $C_{p^n}$-equvariant to $C_N$-equivariant.
Section~\ref{sec:valuation-algebra}
introduces
the valuation language.
Section~\ref{sec:compression}
develops a theory that reconstructs a genuine equivariant spectrum as the homotopy limit
of a finite diagram. We call it the Tate compression theory.
Section~\ref{sec:local}
translate the homotopy groups and structure maps
of the Tate compression diagram into valuation language, where the computation of homotopy groups are in Appendix~\ref{app:local-calculations}.
Section~\ref{sec:reconstruction}
computes the Green functor
for $C_{p^n}$.
Section~\ref{sec:sylow-assembly}
assembled these results
for all finite cyclic groups.
Section~\ref{sec:tambara}
computes the graded norms.

All tensor products and limits are derived. The words `source' and `target' are frequently used to denote the domain and codomain of the associated morphism.
For an inclusion $K<L$, 'upper' and 'lower' refer to $L$ and $K$, respectively,
and to their homotopy groups and classes.
The functors
$\res_K^G\colon\Sp^G\to\Sp^K$ and $N_K^G\colon\Sp^K\to\Sp^G$
are restriction and the indexed smash product.
On homotopy groups,
$\res_K^G$, $\tr_K^G$, $c_g$ and $N_K^G$ denote restriction, additive
transfer, conjugation and the indexed multiplication when defined.
The symbols $\Res_K^G$ and $\Ind_K^G$ act on representations
and on $RO$-degrees. For
$K\le G$, $E\in\Sp^G$, $V\in RO(K)$ and $W\in RO(G)$, we write
$\pi_V^K(E):=[S^V,\res_K^GE]^K$, $\pi_W^K(E):=\pi_{\Res_K^GW}^K(E)$, $\underline\pi_W^G(E)\colon G/K\longmapsto\pi_W^K(E)$,
and $\pi_\star^K(E)=\bigoplus_V\pi_V^K(E)$ with $V$ running through the
indicated grading group.

\section{Foundations}\label{sec:sylow}
In this section, we introduce two tools for later calculations. First, we show how to recover the Green functor of \(H\underline{\mathbb Z}\) in cyclic group case from the $C_{p^n}$-equivariant cases. Next, we define a type of elements later to extend the degree we essentially compute to the whole $RO(G)$ degree.

Let $G$ be a finite  group and $H=\HZ$.
For fixed $V\in RO(G)$, $\underline\pi_V^G(H)$ is a $G$-Mackey functor, and
$\underline\pi_\star^G(H)=\bigoplus_{V\in RO(G)}\underline\pi_V^G(H)$ is an
$RO(G)$-graded Green functor.

\subsection{Cohomologicality}
We first prove this algebraic property of the Mackey functor:
\begin{proposition}\label{prop:cohomological}\label{prop:cohomological-hz}
For $L\le K\le G$ and $x\in \pi_V^K(H)$,
\[
 \tr_L^K\res_L^K(x)=[K:L]x.
\]
Hence every $\underline\pi_V^G(H)$ is a cohomological Mackey functor.
\end{proposition}

\begin{proof}
Since $\pi_0^K(H)=\mathbb Z$ and $\tr_L^K(1_L)=[K:L]1_K$, Frobenius reciprocity gives
$\tr_L^K\res_L^K(x) =x\tr_L^K(1_L) =[K:L]x$.
\end{proof}
The next theorem shows that the Mackey functor is determined by the Mackey functor after restriction to every Sylow groups and Weyl actions.
\begin{theorem}\label{thm:sylow-reconstruction}\label{thm:sylow-green-reduction}
Let $G$ be finite.  For $G\ne e$ and every prime $p\mid |G|$, choose a Sylow subgroup $P_p$.  The $RO(G)$-graded Green functor $\underline\pi_\star^G(H)$ is determined by
\begin{enumerate}[label=\textup{(\roman*)}]
\item the restrictions of $\underline\pi_\star^G(H)$ to the $P_p$, with the grading read through $RO(G)\to RO(P_p)$;
\item for every $L\le P_p$, the Weyl action of $W_G(L)=N_G(L)/L$ on $\pi^L_{\Res^G_LV}(H)$.
\end{enumerate}
\end{theorem}

\begin{proof}
By Proposition~\ref{prop:cohomological},
\cite[Theorem~1.1]{Angeltveit} reconstructs the additive Mackey functors. For $K=e$, it's clear.
For $K\ne e$, choose a Sylow subgroup $P_p(K)\le K$ for each $p\mid |K|$.
Since $\gcd_p[K:P_p(K)]=1$, choose $b_p\in\mathbb Z$ with
$\sum_pb_p[K:P_p(K)]=1$. Then, for $z\in\pi_V^K(H)$,
\begin{equation}\label{eq:bezout-reconstruction}
 z=\sum_p b_p\tr_{P_p(K)}^K\res_{P_p(K)}^K(z).
\end{equation}
Thus the restriction maps to the Sylow subgroups of $K$ are jointly injective, and for
$x\in\pi_U^K(H)$ and $y\in\pi_V^K(H)$, $U,V\in RO(G)$,
\[
 xy=\sum_pb_p\tr_{P_p(K)}^K
 \bigl(\res_{P_p(K)}^K(x)\res_{P_p(K)}^K(y)\bigr).
\]
Taking $z=1_K$ in \eqref{eq:bezout-reconstruction} reconstructs the unit.
\end{proof}
Thus we get the following corollary.

\begin{corollary}[Cyclic groups]
\label{cor:cyclic-sylow-reduction}
For a nontrivial finite cyclic group $G$, the
$RO(G)$-graded Green functor of $H$ is determined
by the Green functors for $C_{p^n}\leq G$ and corresponding Weyl actions.
\end{corollary}
\begin{proof}
Apply Theorem~\ref{thm:sylow-green-reduction}.
\end{proof}

\begin{remark}
More generally, the corollary holds for
any finite group with cyclic Sylow subgroups.
\end{remark}

\subsection{Elements for degree translation}
\label{sec:oriented-models}

For a cyclic group $C_N$, fix a generator $\gamma_N$; the generator of
$C_M\leq C_N$ is $\gamma_M=\gamma_N^{N/M}$. Write $\lambda_N(k)$ to be the 2-dimensional representation whose underlying space is $\mathbb C$, with basis $(1,i)$, on which
$\gamma_Nz=\exp(2\pi i k/N)z$. For even $N$, write $\sigma_N$ for the one-dimensional sign representation. So we have $\lambda_N(N/2)\cong2\sigma_N$. For $G=C_{p^n}$, we define
\[
 RO_{\mathrm{sel}}(G)=
 \begin{cases}
 \mathbb Z\{1,\lambda(1),\lambda(p),\ldots,\lambda(p^{n-1})\},&p\text{ odd},\\
 \mathbb Z\{1,\sigma,\lambda(1),\ldots,\lambda(2^{n-2})\},&p=2,
 \end{cases}
\]
with the rotation list empty when $p=2,n=1$. We call $RO_{\mathrm{sel}}(G)\leq RO(G)$ the \emph{selected degree}.

For $K\le L$, order $L/K$ by
$1,\gamma_L,\ldots,\gamma_L^{[L:K]-1}$.
Use this order for induced representations and
concatenate the chosen bases for direct sums.
Induced complex representations carry their
complex orientations. In particular, for
$C_{2^k}<C_{2^{k+1}}$, the ordered basis
$(1,\gamma_{2^{k+1}})$ gives an oriented isomorphism
$\Ind_{C_{2^k}}^{C_{2^{k+1}}}\sigma_{2^k}
\cong\lambda_{2^{k+1}}(2^{k-1})$.

Write $X^\vee=F(X,\mathbb S)$ for the dual spectrum.
For a virtual model $V=V^+-V^-$, set
$S^V=S^{V^+}\wedge (S^{V^-})^\vee$. 
Duality is compatible with norm:
\[
 N_K^L\bigl((S^{V^-})^\vee\bigr)
 \xrightarrow{\cong}\bigl(N_K^L S^{V^-}\bigr)^\vee.
\]

An actual representation is \emph{orientable}
if the action of $\gamma_N$ has determinant $1$.
For an orientable actual representation $W$
of $C_{p^n}$,
\cite[Definition~3.12]{HHRNonexistence}
defines the orientation class
\[
 u_W\in\pi_{\dim W-W}^{C_{p^n}}(H),\qquad
 \res_e^{C_{p^n}}(u_W)=1.
\]

The next lemma defines a type of elements so that allow us to extend computations from certain degrees to the full \(RO(G)\)-grading.
\begin{lemma}
\label{lem:normalized-character-units}
Let $G=C_N$ and let $D\in RO(G)$ is orientable.
Suppose $D$ is a finite sum of terms of the following form, with
$\ker\lambda_N(a)=\ker\lambda_N(b)$:
\[
 \lambda_N(b)-\lambda_N(a).
\]
Then
\begin{equation*}
 \res_e^G:\pi_D^G(H)\xrightarrow{\cong}\pi_D^e(H)\cong\mathbb Z.
\end{equation*}
Let $v_D\in\pi_D^G(H)$ be the unique class
satisfying $\res_e^G(v_D)=1$. Then
\begin{equation}\label{eq:character-unit-coherence}
 v_{D+E}=v_Dv_E,\qquad v_{-D}=v_D^{-1},\qquad
 \res_K^L(v_D)=v_{\Res_K^L D},\qquad
 N_K^L(v_D)=v_{\Ind_K^L D}.
\end{equation}
\end{lemma}
\begin{proof}

For $G=C_{p^n}$, $p\nmid r$, and $0\le j<n$,
\cite[Proposition~4.25]{Zeng} gives
\[
S^{\lambda(rp^j)-\lambda(p^j)}\wedge H\simeq H.
\]
By restricting both sides to $e$, the result follows directly.

For $G=C_N\ne e$, take $P_p,b_p$ as in
\eqref{eq:bezout-reconstruction}. Restriction of degree preserves the form 
$\lambda_N(b)-\lambda_N(a)$, where $\ker\lambda_N(a)=\ker\lambda_N(b)$.
Thus the $C_{p^n}$ case and \eqref{eq:bezout-reconstruction} show that
$\operatorname{res}_e^G:\pi_D^G(H)\to\pi_D^e(H)\cong\mathbb Z$
is injective.
Since $\det D=1$, the double-coset formula gives
\[
 \res_e^G\!\left(\sum_p b_p\tr_{P_p}^G
  (v_{\Res_{P_p}^G D})\right)
 =\sum_p b_p[G:P_p]=1.
\]
Hence $\res_e^G$ is also surjective, and thus
an isomorphism.
The class $v_D$ is therefore uniquely defined.
Moreover, $v_Dv_{-D}=1$ by restriction.

The first three identities in
\eqref{eq:character-unit-coherence}
follow by restriction and uniqueness.
Induction preserves the form 
$\lambda_N(b)-\lambda_N(a)$, where $\ker\lambda_N(a)=\ker\lambda_N(b)$,
because if $\chi_b=\chi_a^R$ for
$R\in(\mathbb Z/|L|)^\times$, the extensions
$\psi$ and $\psi^R$ have the same kernel.
Thus $\Ind_K^L D$ satisfies
the second formula in the lemma.
Since
\[
\res_e^L N_K^L(v_D)
=\prod_{g\in L/K}g\cdot\res_e^K(v_D)=1,
\]
the norm identity follows by uniqueness.

\end{proof}
Compare this lemma to
\cite[Theorem~B]{AngeltveitPicard}.
For cyclic groups, \cite[Section~1.2]{BasuDey}
shows that the elements defined in this lemma can be chosen compatibly.
For isomorphic representation  $U\cong W$,
we also have $S^{U-W}\simeq S^0$.
We also define $v_{U-W}$ as the unique class
restricting to $1\in\mathbb Z$.

Unlike $v_D$, orientation classes
$u_W$ generally have different degrees.
For any cyclic group $G$ and orientable actual
representation $W$,
we have
\begin{equation*}
\res_e^G:\pi_{\dim W-W}^G(H)\xrightarrow{\cong}\mathbb Z,
\qquad \res_e^G(u_W)=1.
\end{equation*}
Indeed, the degrees $D$ of $v_D$ satisfy
$\dim D^J=0$ for every $J\le G$, whereas
$\dim(\dim W-W)^J=\dim W-\dim W^J$,
which is nonzero for $J=G$ if $W$ is nontrivial.

For $G=C_{p^n}$, Lemma~\ref{lem:normalized-character-units}
gives
\begin{equation*}
v_D\in\pi_D^G(H)^\times,\qquad
D=\lambda(rp^j)-\lambda(p^j),\qquad
\res_e^G(v_D)=1
\end{equation*}
for $p\nmid r$ and $0\le j<n$.

For $V\in RO(G)$, we can replace each rotation summand
$\lambda(rp^j)$ by $\lambda(p^j)$,
to obtain $V_0\in RO_{\mathrm{sel}}(G)$.
For an arbitrary cyclic group $C_M$, use the same notation $V_0$,
replacing each rotation summand $\lambda_M(k)$ by
$\lambda_M(\gcd(k,M))$.
For even $M$, use $\lambda_M(M/2)=2\sigma_M$.
Then $D=V-V_0$ satisfies
Lemma~\ref{lem:normalized-character-units}.

The following proposition compares the Euler and orientation classes
of $\lambda(rp^j)$ and $\lambda(p^j)$.
\begin{proposition}
\label{prop:character-euler-factor}
Let $G=C_{p^n}$, $0\le j<n$, $r\in\mathbb Z_{>0}$,
$p\nmid r$, and $D=\lambda(rp^j)-\lambda(p^j)$.
For the Euler class $a_W\in\pi_{-W}^G(H)$,
represented by $S^0\to S^W\wedge H$, one has
\begin{equation}\label{eq:character-euler-factor}
a_{\lambda(rp^j)}=r\,v_D^{-1}a_{\lambda(p^j)},
\qquad
u_{\lambda(rp^j)}=v_D^{-1}u_{\lambda(p^j)}.
\end{equation}
\end{proposition}
\begin{proof}
The equivariant map
\[
f:S^{\lambda(p^j)}\longrightarrow
S^{\lambda(rp^j)},\qquad z\longmapsto z^r,
\]
has underlying degree $r$ and fixes $0,\infty$.
Its $H$-Hurewicz image satisfies
\[
[f]_H=r v_D^{-1}\in\pi_{-D}^G(H),\qquad
a_{\lambda(rp^j)}=[f]_H a_{\lambda(p^j)}
\]
by Lemma~\ref{lem:normalized-character-units}.
The formula for $u_{\lambda(rp^j)}$ follows from the
injectivity of $\res_e^G$ in degree $2-\lambda(rp^j)$ and
$\res_e^G(u_{\lambda(rp^j)}) =1=\res_e^G(v_D^{-1}u_{\lambda(p^j)})$.
\end{proof}

We also need the following formula.
For a $G$--linear isomorphism
$f:U\xrightarrow{\cong}W$ between oriented $G$-representations,
let
\begin{equation}\label{eq:oriented-model-determinant}
[S^f]_H\in\pi_{U-W}^G(H),\qquad
[S^f]_H=\sgn\det(f)\,v_{U-W}.
\end{equation}
The induced map is
\begin{equation}\label{eq:oriented-model-coefficient-transport}
f_*:\pi_U^G(H)\longrightarrow\pi_W^G(H),\qquad
f_*(x)=x\circ S^{f^{-1}}=[S^f]_H^{-1}x.
\end{equation}

\section{Valuation language}\label{sec:valuation-algebra}
In this section, we introduce the valuation language used in our calculations.
We fix a prime $p$. The groups below can been seen as quotients of subgroups
in $\mathbb Z[1/p]$.

\subsection{Algebraic preliminary for valuation language}

For $b,c\in\mathbb Z$ with $b<c$, put
$p^b\mathbb Z/p^c\mathbb Z =\{p^b z\bmod p^c\mathbb Z:z\in\mathbb Z\}$.
We allow $c=\infty$, with $p^\infty\mathbb Z=0$ and
$\infty+r=\infty$ for every integer $r$.
Multiplication by $p^{-b}$ gives an isomorphism from this group to
$\mathbb Z/p^{c-b}$ when $c<\infty$, and to $\mathbb Z$ otherwise.

Let $S$ be an additive submonoid of an abelian group, with $0\in S$.
For $V\in S$ choose $b(V)\in\mathbb Z$ and
$c(V)\in\mathbb Z\cup\{\infty\}$, with $b(V)<c(V)$ and $b(0)=0$.
Suppose that, for $V,W\in S$,
\begin{equation*}
\begin{split}
 b(V+W)&\leq b(V)+b(W),\\
 c(V+W)&\leq\min\{c(V)+b(W),\ b(V)+c(W)\}.
\end{split}
\end{equation*}
Write $t^V$ for the monoid-algebra basis element of degree $V$, with
$t^Vt^W=t^{V+W}$. Inside $\mathbb Z[1/p][S]$ define
\[
 I=\bigoplus_{V\in S}p^{b(V)}\mathbb Z\,t^V,
 \qquad J=\bigoplus_{V\in S}p^{c(V)}\mathbb Z\,t^V.
\]

We start with the following theorem.

\begin{theorem}
\label{thm:valuation-quotient-algebra}\label{thm:simple-forest-data}
\label{prop:roots-realize}\label{cor:roots}
Then $I$ is a unital subring, $J$ is a homogeneous ideal of $I$, and
\begin{equation*}
 (I/J)_V=p^{b(V)}\mathbb Z/p^{c(V)}\mathbb Z
 \quad(V\in S),\qquad (I/J)_V=0\quad(V\notin S).
\end{equation*}
For $q\in p^{b(V)}\mathbb Z$ and $r\in p^{b(W)}\mathbb Z$, the product is
$[q][r]=[qr]\in(I/J)_{V+W}$.
The unit is the class of $1$ in degree zero.
\end{theorem}
\begin{proof}
The hypotheses give
$1\in I_0$, $I_VI_W\subseteq I_{V+W}$, $J_VI_W+I_VJ_W\subseteq J_{V+W}$.
Thus $I/J$ has the stated grading and multiplication.
\end{proof}

For a nonzero class $x\in p^b\mathbb Z/p^c\mathbb Z$, define
$v_p(x)=v_p(\widetilde x)\in[b,c)$, $v_p(0)=\infty$,
where $\widetilde x\in p^b\mathbb Z$ is any representative.
For homogeneous $x\in(I/J)_V$ and $y\in(I/J)_W$,
with $V=W$ in the sum inequality and $x,y\ne0$ in the product formula,
\begin{equation*}
 \begin{aligned}
 v_p(x+y)&\geq\min\{v_p(x),v_p(y)\},\\
 v_p(xy)&=\begin{cases}
 v_p(x)+v_p(y),&v_p(x)+v_p(y)<c(V+W),\\
 \infty,&v_p(x)+v_p(y)\geq c(V+W).
 \end{cases}
 \end{aligned}
\end{equation*}

We also have the following proposition.

\begin{proposition}
\label{prop:valuation-maps}\label{prop:functoriality}
For $A=p^b\mathbb Z/p^c\mathbb Z$ and
$A'=p^{b'}\mathbb Z/p^{c'}\mathbb Z$, multiplication by
$\varepsilon p^s$, where $s\in\mathbb Z$, $\varepsilon\in\mathbb Z$ and
$p\nmid\varepsilon$, defines a group homomorphism $A\to A'$ provided
\begin{equation*}
 b+s\geq b',\qquad c+s\geq c'.
\end{equation*}
For $c=\infty$, the second inequality is automatic. For $0\ne[q]\in A$,
\[
 v_p([\varepsilon p^s q])=
 \begin{cases}
 v_p(q)+s,&v_p(q)+s<c',\\
 \infty,&v_p(q)+s\geq c'.
 \end{cases}
\]
\end{proposition}
\begin{proof}
These conditions give
$\varepsilon p^sp^b\mathbb Z\subseteq p^{b'}\mathbb Z$ and
$\varepsilon p^sp^c\mathbb Z\subseteq p^{c'}\mathbb Z$.
For finite $A'$, every nonzero group homomorphism sends $[p^b]$ to
$[\varepsilon p^{b+s}]$ with $b'\leq b+s<c'$; the source relation gives
$c+s\geq c'$. Reduction modulo $p^{c'}\mathbb Z$ gives the valuation formula.
\end{proof}

We also have the following proposition to shift valuation.

\begin{proposition}[Shift valuation]\label{prop:valuation-shifts}
For the same setting of Theorem~\ref{thm:valuation-quotient-algebra}, let
$\lambda:S\to\mathbb Z$ be additive. Replacing $b,c$ by
$b+\lambda,c+\lambda$ gives an isomorphic graded ring.
In degree $V$, the isomorphism is
\begin{equation*}
\begin{aligned}
 p^{b(V)}\mathbb Z/p^{c(V)}\mathbb Z
 &\xrightarrow{\cong}
 p^{b(V)+\lambda(V)}\mathbb Z/p^{c(V)+\lambda(V)}\mathbb Z,\\
 [x]&\longmapsto[p^{\lambda(V)}x].
\end{aligned}
\end{equation*}
The integer $b(V)+b(W)-b(V+W)$ is unchanged.
\end{proposition}
\begin{proof}
Additivity of $\lambda$ makes these maps multiplicative.
In degree $V$, the inverse is
$[x]\mapsto[p^{-\lambda(V)}x]$, where
$x\in p^{b(V)+\lambda(V)}\mathbb Z$.
\end{proof}

For $p\nmid\varepsilon$, multiplication by $\varepsilon$ is invertible
on a finite cyclic $p$-group; on $p^b\mathbb Z$, multiplication by
$\varepsilon$ is invertible exactly when $\varepsilon=\pm1$.
\section{Tate compression}\label{sec:compression}\label{sec:descent}
In this section, we develop the Tate compression theory, which is used in later calculations. The main idea is to reconstruct an equivariant spectrum from a diagram. Under suitable conditions, this diagram consists of its generalized homotopy fixed-point spectrum, Tate spectra, and a terminal family localization. Theorem~\ref{thm:family-chain-reconstruction} proves the equivalence of two categories. The homotopy groups of the spectrum can then be computed from these terms and the comparison maps between them.

\subsection{Factorizations}

In a stable $\infty$-category, let
$C\xrightarrow{\iota}Y\xrightarrow{\eta}Y'$ be a fiber sequence and let $\theta:Y\to T$.
Write $\operatorname{Map}$ for mapping spaces, $\operatorname{map}$ for
mapping spectra, and $[Y,T]=\pi_0\operatorname{Map}(Y,T)$. Write $\operatorname{Path}_{x,y}Z$ for paths from $x$ to $y$ in $Z$.

We first give a basic result.

\begin{proposition}\label{prop:factorization-space}\label{prop:sharp-obstruction}
There is a natural equivalence
\[
 \operatorname{fib}_{\theta}
 \bigl(\operatorname{Map}(Y',T)\longrightarrow
                    \operatorname{Map}(Y,T)\bigr)
 \simeq
 \operatorname{Path}_{\theta\iota,0}\operatorname{Map}(C,T).
\]
A factorization $(\rho:Y'\to T,\rho\eta\simeq\theta)$ exists if and only if
$[\theta\iota]=0\in[C,T]$. When nonempty, the space of factorizations is
a torsor for $\operatorname{Map}(\Sigma C,T)$, acting by concatenation
with loops at $0$ in $\operatorname{Map}(C,T)$.
\end{proposition}

\begin{proof}
Applying the mapping-spectrum functor to the fiber sequence gives
$\operatorname{map}(Y',T)\longrightarrow\operatorname{map}(Y,T) \longrightarrow\operatorname{map}(C,T)$.
Apply $\Omega^\infty$ and take the fiber over $\theta$; any path from
$\theta\iota$ to $0$ gives an equivalence between the factorization
space and $\Omega\operatorname{Map}(C,T)\simeq\operatorname{Map}(\Sigma C,T)$.
\end{proof}

For a family $\mathcal F$ of subgroups, the completion
$\eta_Y:Y\to F((E\mathcal F)_+,Y)$ has fiber
$F(\widetilde E\mathcal F,Y)$.

We also have the following result.

\begin{corollary}\label{cor:complete-target-factor}
If the completion $T\to F((E\mathcal F)_+,T)$ is an equivalence,
precomposition with $\eta_Y$ gives
\begin{equation}\label{eq:complete-mapping-equivalence}
 \operatorname{Map}(F((E\mathcal F)_+,Y),T)
 \xrightarrow{\simeq}\operatorname{Map}(Y,T).
\end{equation}
Here $\operatorname{Map}$ may be taken in $\Sp^G$, in
$\operatorname{Mod}_R(\Sp^G)$ for a ring $G$-spectrum $R$, or in
$\CAlg(\Sp^G)$ or $\CAlg(\operatorname{Mod}_R(\Sp^G))$ for commutative $R$.
\end{corollary}
\begin{proof}
The projections $E\mathcal F\times E\mathcal F\rightrightarrows E\mathcal F$
are homotopy inverses to the diagonal. Hence
\[
 \begin{aligned}
 F((E\mathcal F)_+,F((E\mathcal F)_+,Y))
 &\simeq F((E\mathcal F)_+\wedge(E\mathcal F)_+,Y)\\
 &\simeq F((E\mathcal F)_+,Y),\\
 F((E\mathcal F)_+,\eta_Y)&\simeq\eta_{F((E\mathcal F)_+,Y)}.
 \end{aligned}
\]
Thus completion is a localization, and an inverse to
\eqref{eq:complete-mapping-equivalence} is
\[
 \theta\longmapsto\eta_T^{-1}F((E\mathcal F)_+,\theta).
\]
The inverse homotopies follow from naturality and idempotence.
The diagonal makes completion lax symmetric monoidal and $\eta$
multiplicative. For $R$-modules, restrict scalars along
$R\to F((E\mathcal F)_+,R)$. The same identities hold in all four categories,
where equivalences are detected on underlying spectra.
\end{proof}

\subsection{Morphisms}\label{sec:incidence-category}

In this section, a commutative algebra means an object of
$\CAlg(\Sp^G)$ for the symmetric
monoidal $\infty$-category of genuine $G$-spectra; for
$R\in\CAlg(\Sp^G)$ we write $\operatorname{Mod}_R=\operatorname{Mod}_R(\Sp^G)$ and
$\CAlg_R=\CAlg(\operatorname{Mod}_R)\simeq\CAlg(\Sp^G)_{R/}$.

Fix such a commutative algebra $R$ and a chain of
families of subgroups of $G$,
$\varnothing=\mathcal F_{-1} \subseteq\mathcal F_0\subseteq\cdots \subseteq\mathcal F_{n-1}$, $n\ge 1$. For an $R$-module $M$ write
\begin{equation}\label{eq:general-family-layers}
 \begin{aligned}
 M^{\Phi_i}&=\widetilde E\mathcal F_{i-1}\wedge M
 &&(0\leq i\leq n),\\
 M^{h_i}&=F((E\mathcal F_i)_+,M^{\Phi_i}),\qquad
 M^{t_i}=\widetilde E\mathcal F_i\wedge M^{h_i}
 &&(0\leq i<n).
 \end{aligned}
\end{equation}

Here $M^{\Phi_0}=M$; $\Phi^KM$, $M^{hK}$, and $M^{tK}$ denote
subgroup geometric fixed points, homotopy fixed points, and Tate
constructions. Isotropy separation gives the family Tate square
\cite[Section~17, Proposition~17.2]{GreenleesMayTate};
compare \cite[Section~1.2]{Yan}:

\begin{equation}\label{eq:general-family-box}
 M^{\Phi_i}\simeq
 M^{h_i}\mathop{\times}^{h}_{M^{t_i}}M^{\Phi_{i+1}}.
\end{equation}
The $R^{\Phi_i}$-action on $M^{\Phi_i}$, lax monoidal completion, and
smashing localization give $M^{h_i}$ and $M^{t_i}$ their
$R^{h_i}$- and $R^{t_i}$-module structures. Both restrict to the
original $R$-actions.

For $X\in\Sp^G$, apply \eqref{eq:general-family-layers}. For
$0\leq i<n$, let $\ell_i:X^{h_i}\to X^{t_i}$ be localization and let
$\theta_i:X^{\Phi_{i+1}}\to X^{t_i}$ be the other map into
$X^{t_i}$ in \eqref{eq:general-family-box}.

Let $\mathcal I$ be the incidence category of the path, with objects
$v_0,\ldots,v_n,e_0,\ldots,e_{n-1}$ and nonidentity arrows
\[
 v_i\longrightarrow e_i\longleftarrow v_{i+1}\qquad(0\leq i<n).
\]
For an $\mathcal I$-diagram, call the values at $v_i$ and $e_i$ its
vertex and edge objects, respectively.

The next proposition gives a necessary and sufficient condition
for $\theta_i$ to factor through
$X^{\Phi_{i+1}}\to X^{h_{i+1}}$.

\begin{proposition}\label{prop:path-factorization-criterion}
For $0\leq i<n-1$, the map $\theta_i$ factors in $\Sp^G$ through
$\eta_{i+1}:X^{\Phi_{i+1}}\to X^{h_{i+1}}$ if and only if
\begin{equation*}
 [F(\widetilde E\mathcal F_{i+1},X^{\Phi_{i+1}})
       \longrightarrow X^{\Phi_{i+1}}\xrightarrow{\theta_i}X^{t_i}]=0.
\end{equation*}
\end{proposition}
\begin{proof}
Apply Proposition~\ref{prop:factorization-space} to the fiber of $\eta_{i+1}$.
\end{proof}

For $0\leq i<n-1$, choose maps $\rho_i:X^{h_{i+1}}\to X^{t_i}$ and homotopies
$\rho_i\eta_{i+1}\simeq\theta_i$. Together with the maps $\ell_i$ and
$\rho_{n-1}=\theta_{n-1}:X^{\Phi_n}\to X^{t_{n-1}}$, these maps define
an $\mathcal I$-diagram in $\Sp^G$, with value $X^{\Phi_n}$ at $v_n$.

The chosen factorizations give the following reconstruction.
\begin{theorem}[Reconstruction from factorizations]\label{thm:sharp-path-compression}
For the diagram defined above, there is an equivalence
\begin{equation*}
 X\simeq\holim
 (X^{h_0}\to X^{t_0}\leftarrow X^{h_1}\to\cdots
       \to X^{t_{n-1}}\leftarrow X^{\Phi_n}).
\end{equation*}
\end{theorem}
\begin{proof}
For $i<n-1$, the specified homotopy and associativity of pullbacks give
\[
 \begin{aligned}
 X^{\Phi_i}
 &\simeq X^{h_i}\mathop{\times}^{h}_{X^{t_i}}
       \left(X^{h_{i+1}}\mathop{\times}^{h}_{X^{t_{i+1}}}
                        X^{\Phi_{i+2}}\right)\\
 &\simeq\holim
   \left(X^{h_i}\to X^{t_i}\leftarrow X^{h_{i+1}}
                     \to X^{t_{i+1}}\leftarrow X^{\Phi_{i+2}}\right).
 \end{aligned}
\]
Descending induction from
$X^{\Phi_{n-1}}\simeq X^{h_{n-1}}\times^h_{X^{t_{n-1}}}X^{\Phi_n}$
proves the formula in the theorem.
\end{proof}
\subsection{The Tate compression}
Motivated by the above theorem, we now develop the Tate compression theory.
We need the following lemma.

\begin{lemma}\label{lem:euler-nilpotent-completion}
Let $W$ be a $G$-representation with $\widetilde E\mathcal F\simeq S^{\infty W}
=\hocolim_k S^{kW}$, let $A\in\CAlg(\Sp^G)$ (or, more generally, an
associative ring $G$-spectrum), and suppose that, for some $N\geq1$,
the following equality holds in $[A,S^{NW}\wedge A]_A$:
\[
 [A\xrightarrow{\ a_W^N\ }S^{NW}\wedge A]=0,
\]
where $a_W^N$ is induced by the Euler map $a_{NW}\colon S^0\to S^{NW}$.
Then every $A$-module $T$ satisfies
\[
 F(\widetilde E\mathcal F,T)\simeq0,
 \qquad \widetilde E\mathcal F\wedge T\simeq0,
 \qquad T\xrightarrow{\simeq}F((E\mathcal F)_+,T).
\]
\end{lemma}

\begin{proof}
The free-module adjunctions give
\[
 \begin{gathered}
 [A,S^{NW}\wedge A]_A
 \cong\pi_{-NW}^G A
 \cong[A,S^{NW}\wedge A]_{A^{\mathrm{op}}},\\
 [a_W^N]_A=0\quad\Longleftrightarrow\quad a_{NW}\cdot1_A=0
 \quad\Longleftrightarrow\quad[a_W^N]_{A^{\mathrm{op}}}=0,\\
 [T\xrightarrow{a_W^N}S^{NW}\wedge T]^G
 =[(A\xrightarrow{a_W^N}S^{NW}\wedge A)\otimes_A T]^G=0.
 \end{gathered}
\]
In
\[
 F(S^{\infty W},T)\simeq\holim_k\Sigma^{-kW}T,
 \qquad S^{\infty W}\wedge T\simeq\hocolim_k\Sigma^{kW}T,
\]
the $N$-fold transition composites are null. For $K\leq G$ and
$q\in\mathbb Z$, therefore,
\[
 \begin{gathered}
 \lim_k\pi_q^K(\Sigma^{-kW}T)
 =\lim\nolimits_k^1\pi_{q+1}^K(\Sigma^{-kW}T)=0,
 \qquad\operatorname*{colim}_k\pi_q^K(\Sigma^{kW}T)=0,\\
 \pi_q^K F(S^{\infty W},T)
 =\pi_q^K(S^{\infty W}\wedge T)=0,\\
 F(\widetilde E\mathcal F,T)\simeq0,
 \qquad\widetilde E\mathcal F\wedge T\simeq0.
 \end{gathered}
\]
The completion fiber sequence gives
$T\xrightarrow{\simeq}F((E\mathcal F)_+,T)$.
\end{proof}

From now on, assume that, for $0\leq i<n-1$, there are representations
$W_{i+1}$ and integers $N_i\geq1$ such that
\begin{equation}\label{eq:uniform-family-nilpotence}
 \widetilde E\mathcal F_{i+1}\simeq S^{\infty W_{i+1}},
 \qquad
 [R^{t_i}\xrightarrow{a_{W_{i+1}}^{N_i}}
 S^{N_iW_{i+1}}\wedge R^{t_i}]_{R^{t_i}}=0.
\end{equation}
For $0\leq i<n$, let $\mathcal C_i\subseteq\operatorname{Mod}_R$ be the
full $\infty$-subcategory of modules $A$ satisfying
\begin{equation*}
 A\simeq\widetilde E\mathcal F_{i-1}\wedge A,
 \qquad A\xrightarrow{\simeq}F((E\mathcal F_i)_+,A),
\end{equation*}
and let $\mathcal C_n\subseteq\operatorname{Mod}_R$ consist of modules
$A\simeq\widetilde E\mathcal F_{n-1}\wedge A$. Define
$\mathcal D_R$ to be the full subcategory of
$\operatorname{Fun}(\mathcal I,\operatorname{Mod}_R)$ on diagrams
\begin{equation}\label{eq:local-section-category}
 A_0\longrightarrow T_0\longleftarrow A_1
 \longrightarrow\cdots\longrightarrow T_{n-1}\longleftarrow A_n
\end{equation}
such that the value $A_i$ at $v_i$ belongs to $\mathcal C_i$ and, for
$0\leq i<n$, the arrow $\ell_i:A_i\to T_i$ is the localization
$A_i\to\widetilde E\mathcal F_i\wedge A_i$ under an equivalence
$T_i\simeq\widetilde E\mathcal F_i\wedge A_i$.
The arrow $\rho_i:A_{i+1}\to T_i$ is part of the diagram.

The next theorem gives the equivalence of two categories, which is the core of so-called Tate compression theory.

\begin{theorem}[Tate compression]
\label{thm:family-chain-reconstruction}
Under \eqref{eq:uniform-family-nilpotence}, there is an equivalence of $\infty$-categories
\begin{equation}\label{eq:module-category-descent}
 \operatorname{Mod}_R(\Sp^G)
 \xrightarrow{\ \simeq\ }\mathcal D_R,
 \qquad
 M\longmapsto
 (M^{h_0}\to M^{t_0}\leftarrow\cdots
            \to M^{t_{n-1}}\leftarrow M^{\Phi_n}).
\end{equation}
The inverse is the finite homotopy limit.
The same holds for commutative $R$-algebras,
with underlying module diagrams in $\mathcal D_R$.
\end{theorem}

\begin{proof}
Each $M^{t_i}$ is an $R^{t_i}$-module, so, for $i<n-1$,
Lemma~\ref{lem:euler-nilpotent-completion} gives
\[
 M^{t_i}\xrightarrow{\simeq}F((E\mathcal F_{i+1})_+,M^{t_i}).
\]
For $i<n-1$, Corollary~\ref{cor:complete-target-factor} gives natural
factorizations
\[
 M^{\Phi_{i+1}}\longrightarrow M^{h_{i+1}}
 \xrightarrow{\rho_i}M^{t_i}.
\]
Set $\rho_{n-1}=\theta_{n-1}:M^{\Phi_n}\to M^{t_{n-1}}$.
For $0\leq i<n$ and $L\in\mathcal F_{i-1}$,
$\res_L^G M^{\Phi_i}\simeq0$, and hence
$\res_L^G M^{h_i}\simeq0$.
Since $M^{h_i}$ is $\mathcal F_i$-complete,
it lies in $\mathcal C_i$.
Also $M^{\Phi_n}\in\mathcal C_n$ and
$M^{t_i}\simeq\widetilde E\mathcal F_i\wedge M^{h_i}$.
These define the functor in \eqref{eq:module-category-descent}.

Conversely, take an object of $\mathcal D_R$ as in
\eqref{eq:local-section-category}.
Since
$A_i\simeq \widetilde{E\mathcal F}_{i-1}\wedge A_i$, $R^{\Phi_i}=\widetilde{E\mathcal F}_{i-1}\wedge R$,
the $R$-action $\mu:R\wedge A_i\to A_i$ extends to
\[
\begin{aligned}
R^{\Phi_i}\wedge A_i
&\simeq \widetilde{E\mathcal F}_{i-1}
   \wedge R\wedge A_i\\
&\xrightarrow{1\wedge\mu}
\widetilde{E\mathcal F}_{i-1}\wedge A_i
\simeq A_i.
\end{aligned}
\] Lax monoidal completion and
$A_i\simeq F((E\mathcal F_i)_+,A_i)$ give an $R^{h_i}$-action;
localization gives an $R^{t_i}$-action on $T_i$.
These actions restrict to the $R$-actions. Thus $T_i$ is
$\mathcal F_{i+1}$-complete for $i<n-1$.

Define
\[
 M_n=A_n,\qquad
 M_i=A_i\mathop{\times}^{h}_{T_i}M_{i+1}\quad(i=n-1,\ldots,0),
\]
using the composite
$M_{i+1}\to A_{i+1}\xrightarrow{\rho_i}T_i$. The first map is in the
pullback diagram, or the identity when $i+1=n$.
A $\widetilde E\mathcal F_i$-local module $Z$ satisfies
$F((E\mathcal F_i)_+,Z)\simeq0$, since
$\res_L^G Z\simeq0$ for $L\in\mathcal F_i$.
By descending induction,
\[
 \begin{aligned}
 \widetilde E\mathcal F_{i-1}\wedge M_i&\simeq M_i,\\
 \widetilde E\mathcal F_i\wedge M_i
 &\simeq T_i\mathop{\times}^{h}_{T_i}M_{i+1}\simeq M_{i+1},\\
 F((E\mathcal F_i)_+,M_i)
 &\simeq A_i\mathop{\times}^{h}_{0}0\simeq A_i.
 \end{aligned}
\]
Under these equivalences, the projection $M_i\to A_i$ is the
completion $M_i\to F((E\mathcal F_i)_+,M_i)$.
Hence $M_i\simeq M_0^{\Phi_i}$ for $0\leq i\leq n$, and
$A_i\simeq M_0^{h_i}$, $T_i\simeq M_0^{t_i}$ for $i<n$.
For $i<n-1$, completeness gives
\[
 \operatorname{Map}_R(A_{i+1},T_i)
 \xrightarrow{\simeq}\operatorname{Map}_R(M_{i+1},T_i).
\]
Thus the composite $M_{i+1}\to T_i$ recovers $\rho_i$ and its
factorization homotopy with contractible choice. For $i=n-1$,
$M_n=A_n$. These natural equivalences recover the diagram.
Theorem~\ref{thm:sharp-path-compression} recovers $M$ from the terms of its Tate compression diagram,
so the two functors are inverse.

For commutative $R$-algebras, the units and localizations are
multiplicative, Corollary~\ref{cor:complete-target-factor} gives
contractible spaces of factorizations in $\CAlg_R$, and limits are created in
underlying modules. The same construction gives the inverse in $\CAlg_R$.
\end{proof}

Theorem~7.2 of \cite{AMGR} gives symmetric monoidal Tate
reconstruction for $\mathbb Z$-linear derived Mackey functors.
The theorem above concerns genuine $R$-modules and commutative
$R$-algebras under the Euler nilpotence hypothesis
\eqref{eq:uniform-family-nilpotence}.
\begin{definition}[Tate compression diagram]\label{def:tate-compression-diagram}
For an $R$-module $M$ under \eqref{eq:uniform-family-nilpotence},
the diagram in \eqref{eq:module-category-descent} is called the
\emph{Tate compression diagram} of $M$.
Its \emph{vertex objects} are $M^{h_i}$ for $0\leq i<n$ and $M^{\Phi_n}$,
and its \emph{edge objects} are $M^{t_i}$ for $0\leq i<n$.
The \emph{incidence map} is
\[
 \begin{gathered}
 \Delta:\left(\prod_{i<n}M^{h_i}\right)\times M^{\Phi_n}
 \longrightarrow\prod_{i<n}M^{t_i},\\
 \Delta(x)_i=\ell_i(x_i)-\rho_i(x_{i+1}).
 \end{gathered}
\]
For $V\in RO(G)$, write $\Delta_V$ for the induced map of Mackey functors
\begin{equation}\label{eq:ch4-finite-difference}
 \Delta_V:\left(\prod_{i<n}\underline\pi_VM^{h_i}\right)
              \times\underline\pi_VM^{\Phi_n}
 \longrightarrow\prod_{i<n}\underline\pi_VM^{t_i}.
\end{equation}
For each $K\leq G$, elements of the source and target of
$\Delta_V(G/K)$ are called \emph{vertex arrays} and
\emph{edge arrays}, respectively.
\end{definition}

Thus the intermediate terms
$M^{\Phi_i}$ for $1\leq i<n$ have been eliminated
from the reconstruction.

The homotopy limit gives the following fiber sequence.

\begin{corollary}[The fiber sequence]\label{cor:module-incidence-fiber}
For every $R$-module $M$, there is a natural fiber sequence
\begin{equation}\label{eq:universal-module-fibre}
 M\longrightarrow
 \left(\prod_{i<n}M^{h_i}\right)\times M^{\Phi_n}
 \xrightarrow{\Delta}\prod_{i<n}M^{t_i}.
\end{equation}
\end{corollary}
\begin{proof}
Apply Theorem~\ref{thm:family-chain-reconstruction} to the finite path.
Its homotopy limit is the fiber of $\Delta$.
\end{proof}
The following definition is from \cite[Definition~4.8]{MNN}.
\begin{definition}[$B$-exponent and Amitsur descent]
\label{def:b-exponent-amitsur}
Let $R\to B$ be a map of commutative algebras and let $M$ be an
$R$-module. Write $I=\operatorname{fib}(R\to B)$, with $j:I\to R$.
With $\min\varnothing=\infty$, the $B$-exponent is
\[
 \exp_B(M)=\min\{d\geq0:
 j^{\otimes_R d}\otimes_R\operatorname{id}_M\simeq0\}
\]
\cite[Definition~4.8]{MNN}. The module $M$ satisfies Amitsur descent
along $R\to B$ if the canonical map
\[
 M\longrightarrow
 \operatorname{Tot}\bigl(M\otimes_R B^{\otimes_R(\bullet+1)}\bigr)
\]
is an equivalence. Here \(\operatorname{Tot}\) denotes the limit of the Amitsur cosimplicial diagram associated to \(R\to B\). Its coface maps are induced by the unit \(R\to B\), and its codegeneracy maps are induced by the multiplication \(B\otimes_R B\to B\).
\end{definition}

We now give the $B$-exponent bound.

\begin{proposition}[B-exponent]\label{prop:uniform-tensor-descent}
Let
\begin{equation}\label{eq:vertex-tensor-descent}
 B=\left(\prod_{i<n}R^{h_i}\right)\times R^{\Phi_n},\qquad
 I=\operatorname{fib}(R\xrightarrow{u}B),\qquad j:I\longrightarrow R.
\end{equation}
For every $R$-module $M$,
\begin{equation}\label{eq:tensor-descent-exponent}
 [j^{\otimes_R2}\otimes_R\operatorname{id}_M]
 =0\in[I\otimes_R I\otimes_R M,M]_R.
\end{equation}
Thus $\exp_B(R)\leq2$.
\end{proposition}
\begin{proof}
For $M=R$, \eqref{eq:universal-module-fibre} gives
$I\simeq\Omega\prod_{i<n}R^{t_i}$ as $R$-modules.
The algebra map
\[
 B\xrightarrow{(\ell_i\operatorname{pr}_i)_{i<n}}\prod_{i<n}R^{t_i}
\]
equips $I$ with a $B$-module structure extending this $R$-action. Let
$\mu:B\otimes_R I\to I$ to be this action, then
\[
 j\otimes_R\mathrm{id}_I
 =\mu\bigl((uj)\otimes_R\mathrm{id}_I\bigr)\simeq0,
 \qquad
 j^{\otimes_R2}=j\,(j\otimes_R\mathrm{id}_I)\simeq0.
\]
Tensoring with $M$ proves
\eqref{eq:tensor-descent-exponent}.

\end{proof}

For $G=C_{p^n}$ and $R=H$, this implies that
the products of classes in shifted degrees (defined in
Section~\ref{sec:reconstruction}) vanish, since
$F^1\pi_\star^G(H)=\pi_\star^G(H)^{\mathrm{edge}}$
by Theorem~\ref{thm:padic-homotopy-ring}.

\begin{remark}\label{rem:amitsur-base-change}
Write $\operatorname{Tot}_q$ for the limit over cosimplicial degrees
$0,\ldots,q$.
In particular,
\begin{equation}\label{eq:amitsur-module-descent}
 M\xrightarrow{\simeq}
 \operatorname{Tot}\bigl(M\otimes_R B^{\otimes_R(\bullet+1)}\bigr).
\end{equation}
For $S\in\CAlg_R$ and every $S$-module $N$,
$\exp_{S\otimes_R B}(N)\leq2$, and
\[
 N\xrightarrow{\simeq}
 \operatorname{Tot}\bigl(N\otimes_S(S\otimes_R B)^{\otimes_S(\bullet+1)}\bigr).
\]

By \cite[Proposition~2.14]{MNN},
\[
 \operatorname{fib}\!\left(M\to
 \operatorname{Tot}_q(M\otimes_R B^{\otimes_R(\bullet+1)})\right)
 \simeq I^{\otimes_R(q+1)}\otimes_R M.
\]
The transition maps are
$(j\otimes_R\operatorname{id}_I)\otimes_R
\operatorname{id}_{I^{\otimes_R q}\otimes_R M}\simeq0$
by Proposition~\ref{prop:uniform-tensor-descent}.
Hence the fiber tower has zero limit, proving
\eqref{eq:amitsur-module-descent}.
After base change,
\[
 \operatorname{fib}(S\to S\otimes_R B)\simeq S\otimes_R I,
 \qquad
 (S\otimes_R j)^{\otimes_S2}
 \simeq S\otimes_R j^{\otimes_R2}\simeq0.
\]
Same argument applies over $S$.
\end{remark}

We also have the following result for the localization functors.

\begin{corollary}\label{cor:zero-composite-gluing}
For $0\leq i<n$ and $A\in\mathcal C_i$, naturally in $A$,
\[
 \widetilde E\mathcal F_i\wedge A\in\mathcal C_{i+1},
 \qquad
 \widetilde E\mathcal F_{i+1}\wedge\widetilde E\mathcal F_i\wedge A\simeq0
 \quad(i<n-1).
\]
\end{corollary}

\begin{proof}
For $A\in\mathcal C_i$, the module $\widetilde E\mathcal F_i\wedge A$ is
$\widetilde E\mathcal F_i$-local. For $i<n-1$, this localization is an
$R^{t_i}$-module, so Lemma~\ref{lem:euler-nilpotent-completion} gives
the result.
\end{proof}

This gives the following corollary.

\begin{corollary}[Mapping spectra]
\label{cor:mapping-incidence}
For $M,N\in\operatorname{Mod}_R$, there is a natural fiber sequence
\begin{equation*}
 \begin{gathered}
 \operatorname{map}_R(M,N)\longrightarrow
 \left(\prod_{i<n}\operatorname{map}_R(M^{h_i},N^{h_i})\right)
 \times\operatorname{map}_R(M^{\Phi_n},N^{\Phi_n})\\
 \xrightarrow{d}
 \left(\prod_{i<n-1}\operatorname{map}_R(M^{h_{i+1}},N^{t_i})\right)
 \times\operatorname{map}_R(M^{\Phi_n},N^{t_{n-1}}),
 \end{gathered}
\end{equation*}
where $f=(f_0,\ldots,f_n)$,
$f_i:M^{h_i}\to N^{h_i}$ for $i<n$, $f_n:M^{\Phi_n}\to N^{\Phi_n}$, and
\[
 d(f)_i=(\widetilde E\mathcal F_i\wedge f_i)\rho_i^M
                         -\rho_i^N f_{i+1}.
\]
\end{corollary}
\begin{proof}
By (4.4),
$M^{t_i}\simeq \widetilde{E\mathcal F_i}\wedge M^{h_i}$, $\widetilde{E\mathcal F_i}\wedge N^{t_i} \simeq N^{t_i}$.
Thus $N^{t_i}$ is
$\widetilde{E\mathcal F_i}$-local.
The localization adjunction gives
\[
\operatorname{map}_R(M^{t_i},N^{t_i})
\xrightarrow[\simeq]{\,(\ell_i^M)^*\,}
\operatorname{map}_R(M^{h_i},N^{t_i}).
\]
Thus a map between the Tate compression diagrams of $M$ and $N$
consists of maps $f_i$ and nullhomotopies of $d(f)_i$.
The formula in the corollary then follows from
\eqref{eq:module-category-descent}.
\end{proof}

The same fiber sequence detects nilpotence.

\begin{proposition}[Nilpotence detection]
\label{prop:nilpotence-detection}
Let $f:M\to N$ and $g:N\to P$ be $R$-module maps.
If $f^{h_i}=g^{h_i}=0$ for $i<n$ and
$f^{\Phi_n}=g^{\Phi_n}=0$ in the homotopy category,
then $gf=0$. Thus the subgroups
\[
\{[f]\in[M,N]_R:f^{h_i}=0\ (i<n),\
f^{\Phi_n}=0\}
\]
form a square-zero ideal, independently of $n$.

For $B$ in \eqref{eq:vertex-tensor-descent},
if $B\otimes_R f\simeq0$ and $B\otimes_R g\simeq0$,
then $f\otimes_R g=0$ and $gf=0$.
For $r\geq1$,
\begin{equation}\label{eq:tensor-nilpotence-detection}
(B\otimes_R f)^{\otimes_B r}=0
\quad\Longrightarrow\quad
f^{\otimes_R 2r}=0.
\end{equation}
All assertions hold after commutative base change.
\end{proposition}
\begin{proof}
By \eqref{eq:universal-module-fibre},
$f=\partial_N\widetilde f$ for some
$\widetilde f:M\to\Omega\prod_{i<n}N^{t_i}$,
where
$\partial_N:\Omega\prod_{i<n}N^{t_i}\to N$
is the boundary map.
Since $g^{h_i}=0$, also $g^{t_i}=0$. Hence
\[
gf=g\partial_N\widetilde f
=\partial_P\,\Omega\!\left(
\prod_{i<n}g^{t_i}\right)\widetilde f=0.
\]
If $B\otimes_R f\simeq0$, then
$(u\otimes_R N)f\simeq0$.
Hence
$f=(j\otimes_R N)\widetilde f$
for some $\widetilde f:M\to I\otimes_R N$.
Similarly,
$g=(j\otimes_R P)\widetilde g$.
By \eqref{eq:tensor-descent-exponent},
\[
f\otimes_R g=0,\qquad gf=0.
\]
Since
$B\otimes_R f^{\otimes_R r} \simeq (B\otimes_R f)^{\otimes_B r}$,
the implication
$B\otimes_R h\simeq0\Rightarrow
h^{\otimes_R2}\simeq0$
applied to $h=f^{\otimes_R r}$ proves
\eqref{eq:tensor-nilpotence-detection}.

For $S\in\CAlg_R$, the maps
$R^{t_i}\to S^{t_i}$ preserve Euler classes,
so \eqref{eq:uniform-family-nilpotence}
holds for $S$.
The preceding proof in $\operatorname{Mod}_S$
gives $gf=0$ whenever $f$ and $g$
vanish on all $h_i$- and $\Phi_n$-terms.

By Proposition~\ref{prop:uniform-tensor-descent}
and Remark~\ref{rem:amitsur-base-change},
$(S\otimes_R j)^{\otimes_S2}\simeq0$.
Thus, if
$(S\otimes_R B)\otimes_S f\simeq0$ and
$(S\otimes_R B)\otimes_S g\simeq0$,
then $f\otimes_S g\simeq0$ and $gf=0$.
Applying this to $f^{\otimes_S r}$
proves the tensor-nilpotence formula
after base change.
\end{proof}

\subsection{The descent filtration}

For an $R$-module $M$ under \eqref{eq:uniform-family-nilpotence} and
$V\in RO(G)$, regard $\Delta_V$ as a complex of Mackey functors in cohomological degrees $0,1$.

We now define the descent filtration.

\begin{definition}\label{def:descent-filtration}
For an $R$-module $M$, the descent filtration is $F^0=\underline\pi_VM$,\newline
$F^1=\ker\!\left(\underline\pi_VM\to
 (\prod_{i<n}\underline\pi_VM^{h_i})\times\underline\pi_VM^{\Phi_n}\right)$,
and $F^s=0$ for $s\geq2$.
\end{definition}

The fiber sequence gives its two pieces.

\begin{proposition}\label{prop:descent-filtration-ses}
The fiber sequence \eqref{eq:universal-module-fibre} induces short exact sequences
\begin{equation}\label{eq:Mackey-cellular-ses}
 0\longrightarrow\coker\Delta_{V+1}
 \longrightarrow\underline\pi_VM
 \longrightarrow\ker\Delta_V\longrightarrow0.
\end{equation}
\end{proposition}
\begin{proof}
Apply the Mackey functors of homotopy group to \eqref{eq:universal-module-fibre}.
\end{proof}

\leadbreak
For algebras, the descent filtration is compatible with products.

\begin{theorem}[Products in the descent filtration]
\label{thm:Mackey-cellular}\label{thm:filtered-green-descent}
\label{thm:mackey-sheaf-descent}
Let $A\in\CAlg_R$. Then $F^1\underline\pi_\star A$ is a
homogeneous Green ideal of the graded Green functor $\underline\pi_\star A$ and, for $a,b\geq0$,
\begin{equation}\label{eq:square-zero-descent}
 F^aF^b\subseteq F^{a+b},\qquad (F^1)^2=0;
\end{equation}
in particular $\underline\pi_\star A/F^1\cong\bigoplus_V\ker\Delta_V$ is a
graded Green functor with product $(x_i)_i(y_i)_i=(x_iy_i)_i$. Moreover, the product of
$x\in\pi_U^KA$ with image $(x_i)\in(\ker\Delta_U)(G/K)$ and of
$[z]\in(\coker\Delta_{V+1})(G/K)=F^1\pi_V^KA$,
$z=(z_i)_i\in\prod_{i<n}\pi_{V+1}^KA^{t_i}$, is
\begin{equation}\label{eq:mixed-AW}
 x[z]=[(\ell_i(x_i)z_i)_i]
      =[(\rho_i(x_{i+1})z_i)_i]\in F^1\pi_{U+V}^KA .
\end{equation}
\end{theorem}
\begin{proof}
For $A\in\CAlg_R$, the terms in
\eqref{eq:module-category-descent} are
commutative $A$-algebras. Thus \eqref{eq:universal-module-fibre} for $M=A$ is a
fiber sequence of $A$-modules. Write
\[
d:A\longrightarrow
\left(\prod_{i<n}A^{h_i}\right)\times A^{\Phi_n},
\qquad
\partial:\Omega\prod_{i<n}A^{t_i}\longrightarrow A
\]
for the vertex map and boundary map, respectively.
Both are maps of $A$-modules, and $d$ is also multiplicative.
Moreover, $F^1=\ker d_*=\operatorname{im}\partial_*$.
For $a\in\pi_U^K A$ and
$z\in\pi_{V+1}^K\prod_{i<n}A^{t_i}$, we have
\[
a\,\partial_*(z)
=\partial_*(a\cdot z)
=\partial_*\bigl((\ell_i(d_*(a)_i)z_i)_i\bigr).
\]
The $A$-action on $F(S^1,\prod_{i<n}A^{t_i})$
is induced by the $A$-action on $\prod_{i<n}A^{t_i}$.
The convention
$S^{V+1}=S^V\wedge S^1$ leaves $S^1$ on the right, so no sphere
symmetry occurs in the above formula.

Because $d$ is a ring map, $F^1=\ker d_*$ is a homogeneous Green
ideal, with quotient $\operatorname{im}d_*=\bigoplus_V\ker\Delta_V$.
For $a=\partial_*(\tilde a)$ and $b\in F^1$,
\[
 ba=\partial_*\bigl((\ell_i(d_*(b)_i)\tilde a_i)_i\bigr)=0,
\]
proving \eqref{eq:square-zero-descent}. Applying
the second formula in this proof to $x$ and $z$ gives
\eqref{eq:mixed-AW}, since $\ell_i(x_i)=\rho_i(x_{i+1})$.
\end{proof}

\begin{remark}\label{rem:module-descent-filtration}
For every $R$-module $M$, restriction, transfer, and Weyl actions
preserve $F^\bullet\underline\pi_\star M$. The maps in
\eqref{eq:universal-module-fibre} are maps of genuine $G$-spectra,
so they induce maps of homotopy Mackey functors.
\end{remark}

\begin{remark}
If instead one reads $\pi_{V+1}$ through $S^{V+1}=S^1\wedge S^V$, the
second formula in the preceding proof and \eqref{eq:mixed-AW} have an extra factor.
This factor is the image in $\pi_0^K(A)$ of the Burnside-ring unit
obtained by swapping $S^{\Res_K^G U}$ and $S^1$.
Its mark at $J\leq K$ is $(-1)^{\dim(U^J)}$.
For $A=\HZ$, the map from the Burnside ring to
$\pi_0^K(A)=\mathbb Z$ is augmentation, so the factor is
$(-1)^{\dim U}$. The convention $S^{V+1}=S^V\wedge S^1$ gives factor $1$.
\end{remark}

The extension \eqref{eq:Mackey-cellular-ses} need not split in Mackey
functors; \eqref{eq:universal-module-fibre} determines the extension class.

\subsection{Application for $\HZ$}

We can now apply the result to $\HZ$.

\begin{corollary}[Prime-power Tate reconstruction]
\label{cor:hz-module-descent}\label{cor:HZ-compression}\label{prop:adjacent-completeness}
For $G=C_{p^n}$, $n\geq1$, and $R=\HZ$,
Theorem~\ref{thm:family-chain-reconstruction} holds for the chain
$\mathcal F_i=\{L\leq C_{p^i}\}$ and every genuine $\HZ$-module.
In particular
\begin{equation*}
 \HZ\simeq\holim
 (H^{h_0}\to H^{t_0}\leftarrow\cdots\to H^{t_{n-1}}\leftarrow H^{\Phi_n})
\end{equation*}
as commutative algebras, and the corresponding fiber sequence is
\eqref{eq:universal-module-fibre} with $M=\HZ$.
\end{corollary}
\begin{proof}
For $G=C_{p^n}$, write $\lambda_j=\lambda_{p^n}(p^j)$ for $0\leq j<n$,
as in Section~\ref{sec:oriented-models}.
Lemma~\ref{lem:local-euler-height}, proved independently,
gives Euler nilpotence exponent $n-i$ for $i<n-1$. Hence
\begin{equation*}
 [H^{t_i}\xrightarrow{a_{\lambda_{i+1}}^{n-i}}
 S^{(n-i)\lambda_{i+1}}\wedge H^{t_i}]_{H^{t_i}}=0,
 \qquad \widetilde E\mathcal F_{i+1}=S^{\infty\lambda_{i+1}}.
\end{equation*}
At $p=2$, this includes $\lambda_{n-1}=2\sigma_{2^n}$.
Thus \eqref{eq:uniform-family-nilpotence} holds, and
Theorem~\ref{thm:family-chain-reconstruction} gives reconstruction;
Proposition~\ref{prop:uniform-tensor-descent}
gives $\exp_B(M)\leq2$ for every $H$-module $M$.
\end{proof}

For $C_N$ with distinct primes $p,q\mid N$,
the subgroups do not form a chain.
Thus the Tate compression diagram above does not
directly apply in the same form.
We instead paste all its Sylow subgroups to get similar results of the Tate compression theorem.
The next proposition defines $B_G$ and proves $\exp_{B_G}(H)=2$.

\begin{proposition}[$B_G$-exponent]\label{prop:cyclic-descent-exponent}
Let $G$ be any nontrivial finite cyclic group and
$P_p\leq G$ its Sylow $p$-subgroup. Let $B_{P_p}$ be the algebra $B$ in
\eqref{eq:vertex-tensor-descent}, with $G$ replaced by $P_p$
and $R$ replaced by $\res_{P_p}^G H$.Write
\begin{equation}\label{eq:cyclic-sylow-vertex-algebra}
 B_G=\prod_{p\mid |G|}\operatorname{Coind}_{P_p}^G B_{P_p},\qquad
 j_G:\operatorname{fib}(H\longrightarrow B_G)\longrightarrow H.
\end{equation}
Here $\operatorname{Coind}_{P_p}^G=F_{P_p}(G_+,-)$ is coinduction. Then
\begin{equation*}
 j_G\ne0,\qquad j_G^{\otimes_H2}\simeq0.
\end{equation*}
Thus $\exp_{B_G}(H)=2$.
Every $H$-module $M$ satisfies $\exp_{B_G}(M)\leq2$ and
\eqref{eq:amitsur-module-descent} with $R=H$, $B=B_G$.
The exponent bound and formula \eqref{eq:tensor-nilpotence-detection}
persist after arbitrary commutative base change.
\end{proposition}
\begin{proof}
For arbitrary cyclic $G$, coinduction is lax symmetric monoidal,
so $B_G$ is a commutative $H$-algebra. Projection and the counit give
$\res_{P_p}^G B_G\to B_{P_p}$ compatibly with the units. 
The map $\res_{P_p}^G j_G$ factors through
$\operatorname{fib}(\res_{P_p}^G H\to B_{P_p}) \longrightarrow \res_{P_p}^G H$.

By \eqref{eq:tensor-descent-exponent},
$\res_{P_p}^G(j_G^{\otimes_H2})=0 \;\text{for every }p\mid |G|$.
For $H$-modules $M,N$, $F_H(M,N)$
is an $H$-module. The Mackey functor
$\underline\pi_0F_H(M,N)$ is cohomological by
the Frobenius-reciprocity proof of
Proposition~\ref{prop:cohomological-hz}.
Consequently, if $\sum_p b_p[G:P_p]=1$, every $H$-module map $f$ satisfies
\[
 f=\sum_p b_p\tr_{P_p}^G\res_{P_p}^G(f).
\]
Apply this to $f=j_G^{\otimes_H2}$ to obtain $j_G^{\otimes_H2}\simeq0$.
The Amitsur fiber tower now has null two-step transitions;
the totalization, base-change, and tensor-detection proofs of
Remark~\ref{rem:amitsur-base-change} and
Proposition~\ref{prop:nilpotence-detection} apply.

For $G=C_{p^n}$, take $V=2\lambda_0-3$.
In the coordinates of Theorem~\ref{thm:local-padic-identification},
$V+1$ is $(1,-1,\ldots,-1)$. Hence \eqref{eq:local-height-table}, which is proved independently, gives
\[
 \begin{gathered}
 \pi_{V+1}^G(H^{h_0})=\pi_{V+1}^G(H^{\Phi_n})=0,
 \qquad \pi_{V+1}^G(H^{t_0})=p^{-n}\mathbb Z/\mathbb Z,\\
 \pi_{V+1}^G(H^{h_i})=\pi_{V+1}^G(H^{t_i})
 =p^{-(n-i)}\mathbb Z/\mathbb Z\quad(1\leq i<n).
 \end{gathered}
\]
For $1\leq i<n$, the restriction of $\Delta_{V+1}(G/G)$ to
$\pi_{V+1}^G(H^{h_i})$ sends $a$ to $(-a,a)$ in positions $i-1,i$.
Using the inclusions $p^{-(n-i)}\mathbb Z/\mathbb Z\subseteq
p^{-n}\mathbb Z/\mathbb Z$, summation gives
\[
 \coker\bigl(\Delta_{V+1}(G/G)\bigr)
 \xrightarrow{\ \cong\ }p^{-n}\mathbb Z/\mathbb Z,
 \qquad [z_0,\ldots,z_{n-1}]\longmapsto\sum_i z_i,
\]
with inverse $z\mapsto[z,0,\ldots,0]$.
Since $\dim(V^G)=-3$, \eqref{eq:local-height-table} gives $\pi_V^G(B_G)=0$;
thus \eqref{eq:Mackey-cellular-ses} gives
$\pi_V^G(H)\cong\coker(\Delta_{V+1}(G/G))$.

Now let $G$ have arbitrary cyclic order and take
$V=2\lambda_{|G|}(1)-3$. Its restriction to each $P_p$ is the preceding
degree. Rotations are equivariantly homotopic to the identity,
so the Weyl action is trivial. Since $\dim V=1$, one has $\pi_V^e(H)=0$.
On $\pi_V^{P_p}(H)\to\pi_V^{P_q}(H)$, the double-coset formula gives
\[
 \res_{P_q}^G\tr_{P_p}^G=
 \begin{cases}
 [G:P_p],&q=p,\\
 0,&q\ne p.
 \end{cases}
\]
The B\'ezout identity gives $b_q[G:P_q]\equiv1\pmod{|P_q|}$.
For $x_p\in\pi_V^{P_p}(H)\cong\mathbb Z/|P_p|$,
\[
 \res_{P_q}^G\!\left(\sum_p b_p\tr_{P_p}^G(x_p)\right)
 =b_q[G:P_q]x_q=x_q,
 \qquad\sum_p b_p\tr_{P_p}^G\res_{P_p}^G=\mathrm{id}.
\]
Together with coinduction adjunction, we have
\[
 \pi_V^G(H)\cong\prod_{p\mid |G|}\mathbb Z/|P_p|
 \cong\mathbb Z/|G|\ne0,\qquad
 \pi_V^G(B_G)\cong\prod_{p\mid |G|}\pi_{\Res_{P_p}^GV}^{P_p}(B_{P_p})=0.
\]
Hence $j_G$ has nonzero image on homotopy groups and cannot be null.
\end{proof}

In later sections, we assemble both the Green functor
and the descent filtration from the $C_{p^n}$-equivariant cases to the $C_N$-equivariant case.
The equality $j_G^{\otimes_H 2}\simeq0$ suggests
the square-zero property $(F^1)^2=0$, which is proved
in Theorem~\ref{thm:first-order-normed-descent}.
\section{Translation to valuation language}\label{sec:local}\label{sec:local-layers}
Fix $G=C_{p^n}$, where $p$ is prime and $n\geq1$, and still put $H=\HZ$.
We apply Tate compression to $H$ as in Section~\ref{sec:compression}.
We recall the terms of its Tate compression diagram,
then compute their homotopy groups and structure maps in valuation language (Theorem~\ref{thm:local-forest-identification}).

\subsection{Translation to valuation language}

Write $\mathcal F_i=\{C_{p^j}:0\leq j\leq i\}$, and
$\mathcal F_{-1}=\varnothing$. For $0\leq j\leq n$ and $0\leq i<n$,
\eqref{eq:general-family-layers} gives
\begin{equation}\label{eq:layer-definitions}
\begin{aligned}
 H^{\Phi_j}&=\widetilde E\mathcal F_{j-1}\wedge H,
 &H_{h_i}&=(E\mathcal F_i)_+\wedge H^{\Phi_i},\\
 H^{h_i}&=F((E\mathcal F_i)_+,H^{\Phi_i}),
 &H^{t_i}&=\widetilde E\mathcal F_i\wedge H^{h_i}.
\end{aligned}
\end{equation}

By \eqref{eq:general-family-box},
\[
 H^{\Phi_i}\simeq H^{h_i}\times^h_{H^{t_i}}H^{\Phi_{i+1}}.
\]

\begin{definition}\label{def:degree-coordinates}
We describe degrees in $RO_{\mathrm{sel}}(G)$ by fixed-point
coordinates, compatible with the valuation language.
For $V\in RO_{\mathrm{sel}}(G)$ with $\dim(V^G)$ even,
write $\mathbf a=(a_0,\ldots,a_n)$ for the coordinates of $V$, where
$2a_j=\dim(V^{C_{p^j}})$.
These coordinates range over $D_n$:
\[
 D_n^{(0)}=\mathbb Z^{n+1},\qquad
 D_n^{(1)}=(\mathbb Z+\tfrac12)^n\times\mathbb Z,
 \qquad
 D_n=\begin{cases}D_n^{(0)},&p\text{ odd},\\
 D_n^{(0)}\sqcup D_n^{(1)},&p=2.
 \end{cases}
\]
Addition in $D_n$ is coordinatewise, and we write
$\epsilon(\mathbf a)=\epsilon$ for $\mathbf a\in D_n^{(\epsilon)}$.

Write $\lambda_j=\lambda_{p^n}(p^j)$ for $0\leq j<n$ 
as in Section~\ref{sec:oriented-models}.
For odd $p$, write
\begin{equation}\label{eq:odd-forest-coordinate}
 \rho_p(\mathbf a)
 =2a_n+\sum_{j=0}^{n-1}(a_j-a_{j+1})\lambda_j.
\end{equation}
For $p=2$, we have $\lambda_{n-1}=2\sigma_{2^n}$, and define
\begin{equation}\label{eq:two-forest-coordinate}
 \rho_2(\mathbf a)
 =2a_n+2(a_{n-1}-a_n)\sigma_{2^n}
       +\sum_{j=0}^{n-2}(a_j-a_{j+1})\lambda_j.
\end{equation}
For the selected degree $RO_{\mathrm{sel}}(G)$ of
Section~\ref{sec:oriented-models} and $0\leq j\leq n$, we have
\[
 \rho_p:D_n\hookrightarrow RO_{\mathrm{sel}}(G),\qquad
 \dim\bigl(\rho_p(\mathbf a)^{C_{p^j}}\bigr)=2a_j.
\]
\end{definition}

One will see that for $V\in RO_{\mathrm{sel}}(G)\setminus\rho_p(D_n)$
and $0\leq i<n$,
\[
 \pi_V^G(H^{h_i})=\pi_V^G(H^{t_i})=\pi_V^G(H^{\Phi_n})=0.
\]

To reduce the length of the article, we have to define the following variables.
For $\mathbf a=(a_0,\ldots,a_n)\in D_n^{(\epsilon)}$
and $0\leq i<n$, write
\begin{equation*}
 \begin{aligned}
 S_i(\mathbf a)&=\sum_{k=i+1}^n a_k,\qquad
 t_i(\mathbf a)=a_i+\frac{\epsilon(\mathbf a)}2\in\mathbb Z,\\
 S_i(\mathbf a)
                 +\epsilon(\mathbf a)\frac{n-i-1}{2}
 &=a_n+\sum_{k=i+1}^{n-1}t_k(\mathbf a)\in\mathbb Z.
 \end{aligned}
\end{equation*}
Define
\[
 A_i(\mathbf a)=
 \frac{p^{b_{A_i}(\mathbf a)}\mathbb Z}
      {p^{c_{A_i}(\mathbf a)}\mathbb Z},
 \qquad
 T_i(\mathbf a)=
 \frac{p^{b_{T_i}(\mathbf a)}\mathbb Z}
      {p^{c_{T_i}(\mathbf a)}\mathbb Z},
\]
Here we call $b_E$ is the lower valuation bound and $e_E=c_E-b_E$
is the length; equivalently $c_E=b_E+e_E$. Put $a^+=\max(a,0)$,
and let $\mathbf1_P$ be $1$ if $P$ holds and $0$ otherwise.
Define the explicit variables above to be 
\begin{equation}\label{eq:local-height-table}
\begin{aligned}
 (b_{A_0},e_{A_0})&=
 \begin{cases}
  (S_0,n),&\epsilon=0,\ a_0<0,\\
  (S_0,\infty),&\epsilon=0,\ a_0=0,\\
  (S_0+\tfrac{n-1}{2},1),&\epsilon=1,\ t_0\leq0,
 \end{cases}\\[1mm]
 (b_{A_i},e_{A_i})&=
 \begin{cases}
  (S_i+a_i^+,n-i+\mathbf1_{a_i\geq0}),&\epsilon=0,\\
  (S_i+\tfrac{n-i-1}{2}+t_i^+,1),&\epsilon=1,
 \end{cases}
 &&1\leq i<n,\\[1mm]
 (b_{T_i},e_{T_i})&=
 \begin{cases}
  (S_i,n-i),&\epsilon=0,\\
  (S_i+\tfrac{n-i-1}{2},1),&\epsilon=1,
 \end{cases}
 &&0\leq i<n,\\[1mm]
 (b_{A_n},e_{A_n})&=(a_n,1),&&a_n\geq0.
\end{aligned}
\end{equation}
The cases $\epsilon=1$ occur only at $p=2$.
In cases not listed in the table, set $E(\mathbf a)=0$.
The support of $E$ is $\{\mathbf a:E(\mathbf a)\ne0\}$.
Write $A_i=\bigoplus_{\mathbf a\in D_n}A_i(\mathbf a)$ and
$T_i=\bigoplus_{\mathbf a\in D_n}T_i(\mathbf a)$.

We apply Theorem~\ref{thm:valuation-quotient-algebra}
to obtain the following graded rings. In fact, for $0\leq i<n$, $A_i$ and $T_i$ will be equal to  $\pi_\star^G(H^{h_i})$ and $\pi_\star^G(H^{t_i})$, respectively, $A_n$ will be equal to $\pi_\star^G(H^{\Phi_n})$.

\begin{theorem}
\label{thm:local-forest-data}\label{thm:local-padic-data}
With the same construction in Theorem~\ref{thm:valuation-quotient-algebra}, we define unital commutative graded rings for
$E=A_0,\ldots,A_n,T_0,\ldots,T_{n-1}$.
\end{theorem}

\begin{proof}
By Theorem~\ref{thm:valuation-quotient-algebra},
we need to verify
\[
\begin{aligned}
b_E(\mathbf a+\mathbf d)
&\leq b_E(\mathbf a)+b_E(\mathbf d),\\
c_E(\mathbf a+\mathbf d)
&\leq c_E(\mathbf a)+b_E(\mathbf d),\\
c_E(\mathbf a+\mathbf d)
&\leq b_E(\mathbf a)+c_E(\mathbf d)
\end{aligned}
\]
whenever $E(\mathbf a),E(\mathbf d)\ne0$.

Put
$\delta_E(\mathbf a,\mathbf d) =b_E(\mathbf a)+b_E(\mathbf d) -b_E(\mathbf a+\mathbf d)$.
Since $c_E=b_E+e_E$, the three inequalities
are equivalent to
\[
\begin{aligned}
\delta_E(\mathbf a,\mathbf d)&\geq0,\\
e_E(\mathbf a+\mathbf d)
&\leq\min\{e_E(\mathbf a),e_E(\mathbf d)\}
+\delta_E(\mathbf a,\mathbf d),
\end{aligned}
\]
where $\infty+r=\infty$.

Put $\kappa=\epsilon(\mathbf a)\epsilon(\mathbf d)$.
Then
$t_i(\mathbf a+\mathbf d) =t_i(\mathbf a)+t_i(\mathbf d)-\kappa$.
For $1\leq i<n$, we have
\[
\begin{aligned}
\delta_{A_i}(\mathbf a,\mathbf d)
&=t_i(\mathbf a)^++t_i(\mathbf d)^+\\
&\quad-(t_i(\mathbf a)+t_i(\mathbf d)-\kappa)^+
+\kappa(n-i-1)\geq0.
\end{aligned}
\]
If both inputs have $\epsilon=0$, the only
nontrivial length case, up to symmetry, is
$t_i(\mathbf a)<0\leq t_i(\mathbf a+\mathbf d)$.
Then
\[
\begin{aligned}
e_{A_i}(\mathbf a+\mathbf d)-e_{A_i}(\mathbf a)
&=1\\
&\leq-t_i(\mathbf a)
=\delta_{A_i}(\mathbf a,\mathbf d).
\end{aligned}
\]
If exactly one of $\epsilon(\mathbf a)$ and $\epsilon(\mathbf d)$
is $1$, then $e_{A_i}(\mathbf a+\mathbf d)=1$,
so the required inequality holds.
If $\kappa=1$, then
$e_{A_i}(\mathbf a)=e_{A_i}(\mathbf d)=1$, and
\[
\begin{aligned}
\delta_{A_i}(\mathbf a,\mathbf d)
&\geq n-i-1
+\mathbf1_{t_i(\mathbf a+\mathbf d)\geq0}\\
&=e_{A_i}(\mathbf a+\mathbf d)-1.
\end{aligned}
\]

For $T_i$ and $A_0$, we have
$\delta_{T_i}=\kappa(n-i-1)$, $\delta_{A_0}=\kappa(n-1)$.
If $\kappa=0$, the target length is at most
the smaller input length. If $\kappa=1$,
\[
\begin{aligned}
(e_{T_i}(\mathbf a),e_{T_i}(\mathbf d),
e_{T_i}(\mathbf a+\mathbf d))
&=(1,1,n-i),\\
(e_{A_0}(\mathbf a),e_{A_0}(\mathbf d),
e_{A_0}(\mathbf a+\mathbf d))
&=(1,1,n),
\end{aligned}
\]
where the second equality uses
$a_0,d_0\leq-\tfrac12$ on the support of $A_0$.
Thus the required inequalities hold.
\end{proof}

\begin{remark}
For $E(\mathbf a)\ne0$,
\[
 E(\mathbf a)\cong
 \begin{cases}
 \mathbb Z/p^{e_E(\mathbf a)},&e_E(\mathbf a)<\infty,\\
 \mathbb Z,&e_E(\mathbf a)=\infty.
 \end{cases}
\]
\end{remark}

For nonzero $x\in E(\mathbf a)$ and $y\in E(\mathbf d)$,
\begin{equation}\label{eq:local-valuation-product}
 v_p(xy)=
 \begin{cases}
  v_p(x)+v_p(y),&v_p(x)+v_p(y)<c_E(\mathbf a+\mathbf d),\\
  \infty,&v_p(x)+v_p(y)\geq c_E(\mathbf a+\mathbf d).
 \end{cases}
\end{equation}
By Proposition~\ref{prop:valuation-shifts}, adding the same additive
function to $b_E$ and $c_E$ preserves $e_E$, $\delta_E$, and the graded
isomorphism class. The lengths, and $\delta_E$ determine
the graded ring up to isomorphism.

The following theorem computes the homotopy groups in valuation language.

\begin{theorem}[Translation]
\label{thm:local-forest-identification}\label{thm:local-padic-identification}
\label{thm:odd-forest-translation}\label{thm:two-forest-translation}
For $\mathbf a\in D_n$ and $0\leq i<n$,
the following isomorphisms preserve products:
\begin{align*}
 \pi_{\rho_p(\mathbf a)}^G(H^{h_i})
 &\cong A_i(\mathbf a),
 \\
 \pi_{\rho_p(\mathbf a)}^G(H^{t_i})
 &\cong T_i(\mathbf a),
 \\
 \pi_{\rho_p(\mathbf a)}^G(H^{\Phi_n})
 &\cong A_n(\mathbf a).
\end{align*}
\end{theorem}
\begin{proof}
Propositions~\ref{prop:odd-fractional-lattices} and
\ref{prop:two-fractional-lattices} give monomial generators
$g_E(\mathbf a)$ of additive order $p^{e_E(\mathbf a)}$
when $e_E(\mathbf a)<\infty$, and of infinite order
when $e_E(\mathbf a)=\infty$, where $e_E$ is given by
\eqref{eq:local-height-table}. These generators satisfy
\[
 g_E(\mathbf a)g_E(\mathbf d)
 =p^{\delta_E(\mathbf a,\mathbf d)}g_E(\mathbf a+\mathbf d),
 \qquad g_E(0)=1.
\]
Hence the inverse isomorphisms are
\[
 [p^{b_E(\mathbf a)}m]\longmapsto m g_E(\mathbf a).
\]
At $p=2$, $\dim\rho_2(\mathbf a)\equiv\epsilon(\mathbf a)\pmod2$ and
$2g_E(\mathbf a)=0$ when $\epsilon(\mathbf a)=1$; hence the product
formula satisfies graded commutativity.
\end{proof}

For $E=A_i$ $(0\leq i\leq n)$ or $E=T_i$ $(0\leq i<n)$ with
$E(\mathbf a)\ne0$, write $g_E(\mathbf a)$ for the class corresponding to
$[p^{b_E(\mathbf a)}]\in E(\mathbf a)$.
We prove the following Euler nilpotence as the condition used in
Corollary~\ref{cor:hz-module-descent}.

\begin{lemma}
\label{lem:local-euler-height}
For $0\leq i<n-1$ and $N\geq1$,
\[
 \begin{aligned}
 [H^{t_i},S^{N\lambda_{i+1}}\wedge H^{t_i}]_{H^{t_i}}
 &\cong\pi_{-N\lambda_{i+1}}^G(H^{t_i})
 \cong p^{-N}\mathbb Z/p^{n-i-N}\mathbb Z,\\
 [H^{t_i}\xrightarrow{a_{\lambda_{i+1}}^N}
 S^{N\lambda_{i+1}}\wedge H^{t_i}]_{H^{t_i}}&\longleftrightarrow[1].
 \end{aligned}
\]
Thus $a_{\lambda_{i+1}}^{n-i}=0$ in
$\pi_\star^G(H^{t_i})$.
\end{lemma}
\begin{proof}
The first isomorphism is clear. The coordinates corresponding to the degree of the power of Euler class are
\[
 \rho_p(\underbrace{-N,\ldots,-N}_{i+2},
        \underbrace{0,\ldots,0}_{n-i-1})=-N\lambda_{i+1}.
\]
Equation~\eqref{eq:local-height-table} gives the quotient.
Lemma~\ref{lem:local-euler-verification} shows that the Euler class
corresponds to $[1]$ for $N=1$; multiplication of representatives gives $[1]^N=[1]$.
Hence
\[
 [1]=0\text{ in }p^{-N}\mathbb Z/p^{n-i-N}\mathbb Z
 \quad\Longleftrightarrow\quad N\geq n-i.
\]
\end{proof}

\subsection{Comparison maps in terms of valuation language}

For $0\leq i<n$, define graded maps
\[
 \ell_i:A_i\longrightarrow T_i,\qquad
 \rho_i:A_{i+1}\longrightarrow T_i
\]
by $[q]\mapsto[q]$ on nonzero source groups and zero otherwise.
The following identities hold in every degree for which
the corresponding source group is nonzero:
\begin{equation}\label{eq:local-comparison-displacements}
 \begin{aligned}
 b_{A_0}-b_{T_0}&=0,\qquad b_{A_n}-b_{T_{n-1}}=0,\\
 b_{A_i}-b_{T_i}&=t_i^+ &&(1\leq i<n),\\
 b_{A_{i+1}}-b_{T_i}&=(-t_{i+1})^+ &&(0\leq i<n-1).
 \end{aligned}
\end{equation}
For $E=A_i,A_{i+1}$ and $E(\mathbf a)\ne0$, the map
$E(\mathbf a)\to T_i(\mathbf a)$ is zero if and only if
$b_E(\mathbf a)\geq c_{T_i}(\mathbf a)$.

We next give the comparison maps.

\begin{theorem}
\label{thm:local-forest-comparisons}\label{thm:local-padic-comparisons}
$\ell_i$ and $\rho_i$ are well-defined unital graded ring homomorphisms.
Under Theorem~\ref{thm:local-padic-identification}, $\ell_i$ is exactly induced by
the localization $H^{h_i}\to H^{t_i}$.
For $0\leq i<n-1$, $\rho_i$ is also exactly induced by the factorization
\[
 H^{\Phi_{i+1}}\longrightarrow H^{h_{i+1}}
 \longrightarrow H^{t_i};
\]
$\rho_{n-1}$ is induced by $H^{\Phi_n}\to H^{t_{n-1}}$.
\end{theorem}
\begin{proof}
Equation~\eqref{eq:local-height-table} gives
\eqref{eq:local-comparison-displacements}.
In every degree for which the corresponding source group
is nonzero, the upper valuation bounds satisfy
\[
 \begin{aligned}
 e_{A_i}\geq e_{T_i},\quad b_{A_i}\geq b_{T_i}
 &\ \Longrightarrow\ c_{A_i}\geq c_{T_i} &&(0\leq i<n),\\
 c_{A_{i+1}}-c_{T_i}
 &=(-t_{i+1}-1+\epsilon)^+\geq0 &&(0\leq i<n-1),\\
 c_{A_n}&=c_{T_{n-1}}=a_n+1.
 \end{aligned}
\]
Proposition~\ref{prop:valuation-maps} gives well-defined homomorphisms;
$[q][r]=[qr]$ and $1\mapsto1$ give unital graded ring maps.

Lemma~\ref{lem:local-euler-height} and
Corollary~\ref{cor:hz-module-descent} give the factorizations.
By Lemma~\ref{lem:local-height-normalization}, the induced maps on
homotopy send $[q]$ to $[q]$, including $\ell_0$ and $\rho_{n-1}$.
\end{proof}

\subsection{Restriction, transfer and Weyl action in valuation language}

Let $n\geq2$, $G=C_{p^n}$, and $K=C_{p^{n-1}}$.
For $p=2$, the degree of a nonzero $K$-equivariant homotopy group requires
\[
 (a_0,\ldots,a_n)\in\mathbb Z^n\times\tfrac12\mathbb Z.
\]
Although the restriction of $V=\rho_p(\mathbf a)$ to $K$ is determined
by $(a_0,\ldots,a_{n-1})$, we retain the full coordinates
$\mathbf a=(a_0,\ldots,a_n)$ to describe restriction, transfer,
and the Weyl action.
These $K$-equivariant homotopy groups are independent of $a_n$.
More precisely, in the following formulas, $A_i,T_i$ are given by
\eqref{eq:local-height-table} with $n-1$ in place of $n$,
evaluated at $(a_0,\ldots,a_{n-1})$:
\[
\begin{aligned}
 \pi_{\rho_p(\mathbf a)}^K(H^{h_i})&\cong
 A_i(a_0,\ldots,a_{n-1}),&&0\leq i\leq n-1,\\
 \pi_{\rho_p(\mathbf a)}^K(H^{t_i})&\cong
 T_i(a_0,\ldots,a_{n-1}),&&0\leq i<n-1,\\
 \pi_{\rho_p(\mathbf a)}^K(H^{t_{n-1}})
 &=\pi_{\rho_p(\mathbf a)}^K(H^{\Phi_n})=0.
\end{aligned}
\]
Here \eqref{eq:two-forest-coordinate} defines $\rho_2(\mathbf a)$
also for $a_n\in\tfrac12\mathbb Z$.
For $E=H^{h_i},H^{t_i}$ and $V=\rho_p(\mathbf a)$, we have
$\pi_V^G(E),\pi_V^K(E)\ne\nobreak0\ \Longrightarrow\ a_n\in\mathbb Z$.
When both homotopy groups are nonzero, the lower valuation bound for
$\pi_V^G(E)$ equals that for $\pi_V^K(E)$ plus $a_n$.

We have the following result for restriction and transfer.

\begin{theorem}[Restriction and transfer on $H^{h_i},H^{t_i},H^{\Phi_n}$]
\label{thm:local-forest-mackey}
Let $E=H^{h_i},H^{t_i}$, or $H^{\Phi_n}$ and $V=\rho_p(\mathbf a)$ as above.
If $\pi_V^G(E),\pi_V^K(E)\ne0$, then for classes
$[q]\in\pi_V^G(E)$ and $[r]\in\pi_V^K(E)$,
\begin{equation}\label{eq:local-intrinsic-res-tr}
 \res_K^G([q])=[p^{-a_n}q],\qquad
 \tr_K^G([r])=[p^{a_n+1}r].
\end{equation}
If $\pi_V^G(E)=0$ or $\pi_V^K(E)=0$, then
$\res_K^G=\tr_K^G=0$ in degree $V$.
\end{theorem}
\begin{proof}
Propositions~\ref{prop:odd-local-subgroup-maps},
\ref{prop:two-local-subgroup-maps}, and
\ref{prop:terminal-local-subgroup-values} give restriction and transfer. When $\pi_V^G(E),\pi_V^K(E)\ne0$,
for any representative $r$ of $[r]\in\pi_V^K(E)$,
the isomorphism of Theorem~\ref{thm:local-forest-identification} gives
\[
 \begin{gathered}
 [p^{a_n}r]\in\pi_V^G(E),\qquad
 \res_K^G[p^{a_n}r]=[r],\\
 \tr_K^G[r]=\tr_K^G\res_K^G[p^{a_n}r]
 =p[p^{a_n}r]=[p^{a_n+1}r].
 \end{gathered}
\]
Here $\tr_K^G\res_K^G=p$ by the projection formula
for the $H$-module $E$ and $\tr_K^G(1)=p$ in $\pi_0^G(H)$.
\end{proof}

The Weyl action are as follows.

\begin{proposition}[Weyl action]\label{prop:local-weyl-support}
For $E,V,G,K$ as above and $\pi_V^K(E)\ne0$, the generator $\gamma_{p^n}$ acts by
\begin{equation*}
 c_{\gamma_{p^n}}\big|_{\pi_V^K(E)}=
 \begin{cases}
  \operatorname{id},&p\text{ odd},\\
  (-1)^{2a_n}\operatorname{id},&p=2.
 \end{cases}
\end{equation*}
\end{proposition}
\begin{proof}
Apply Propositions~\ref{prop:odd-local-subgroup-maps},
\ref{prop:two-local-subgroup-maps}, and
\ref{prop:terminal-local-subgroup-values}.
\end{proof}
We only have the following trivial exception.

\begin{proposition}[$e<C_p$]
\label{prop:local-prime-base}
Let $G=C_p$. Then
\[
\res_e^G H^{h_0}\simeq H\mathbb Z,\qquad
\res_e^G H^{t_0}\simeq
\res_e^G H^{\Phi_1}\simeq0.
\]
For $V\in RO(G)$ with $\dim V=0$ and $\det(V)=1$,
write $V_0=\rho_p(0,a_1)$ as in
Section~\ref{sec:oriented-models}.
By Lemma~\ref{lem:normalized-character-units}
and Theorem~\ref{thm:local-padic-identification},
we have
\[
\pi_V^G(H^{h_0})\cong p^{a_1}\mathbb Z,\qquad
\pi_V^e(H^{h_0})\cong\mathbb Z.
\]
For $q\in p^{a_1}\mathbb Z$ and $m\in\mathbb Z$, we have
\[
\res_e^G(q)=p^{-a_1}q,\qquad
\tr_e^G(m)=p^{a_1+1}m.
\]
If $p=2$, $\dim V=0$, and $\det(V)=-1$,
then $\pi_V^G(H^{h_0})=0$, so both maps vanish.
If $\dim V\ne0$, then $\pi_V^e(H^{h_0})=0$.
The Weyl action on nonzero $\pi_V^e(H^{h_0})$
is multiplication by $\det(V)$.
\end{proposition}

\begin{proof}
Apply Proposition~\ref{prop:appendix-prime-base}
and Theorem~\ref{thm:local-padic-identification}.
\end{proof}
\section{Prime-power Green functor}\label{sec:reconstruction}\label{sec:cyclic}
In this section, we compute the $RO(C_{p^n})$-graded Green functor
$\underline{\pi}_\star(H)$, including products,
restriction, transfer, and Weyl action.

As usual, put $G=C_{p^n}$, $n\geq1$, and $H=\HZ$. For
$V=\rho_p(\mathbf a)$ with degree coordinates
$\mathbf a=(a_0,\ldots,a_n)\in D_n$, we have
$2a_j=\dim(V^{C_{p^j}})$, as in Definition~\ref{def:degree-coordinates}.

\subsection{Complete prime-power Green functor}

\begin{theorem}[Complete prime-power Green functor]
\label{thm:complete-prime-power-green}
The following results give the $RO_\mathrm{sel}(G)$-graded Green functor
$\underline{\pi}_\star(H)$:
\begin{enumerate}[label=\textup{(\roman*)},leftmargin=2.8em]
\item homotopy groups: Theorem~\ref{thm:prime-power-final} and
Corollary~\ref{cor:array-subgroup-values};
\item products and unit: Proposition~\ref{prop:array-products};
\item restriction and transfer: Proposition~\ref{prop:array-equal-parity}
for most cases; Propositions~\ref{prop:two-natural-shifted-res},
\ref{prop:two-natural-shifted-transfer}, and
\ref{prop:two-shifted-natural-res} for vertex-edge maps (Definition~\ref{def:degree-element-parity}) at $p=2$;
all others are induced by transitivity;
\item Weyl actions: Proposition~\ref{prop:adjacent-weyl}.
\end{enumerate}
\end{theorem}
Then Corollary~\ref{cor:array-character-transport} extends the Green functor from $RO_\mathrm{sel}(G)$ to $RO(G)$.

\subsection{Reconstruction of equivariant homotopy groups}
\label{sec:ch4-tate-pullback}
\begin{definition}
\label{def:degree-element-parity}
For $L\leq G$ and $V\in RO(G)$, $V$ is a \emph{natural degree} at $G/L$
if $\dim(V^L)$ is even and a \emph{shifted degree} otherwise. The term \emph{shifted} comes from the use of
$\Sigma^{-1}$ to denote such elements in earlier computations;
see \cite[Section~5.4]{Yan}.
Elements of $\pi_V^L(H)$ in natural degrees are called \emph{vertex elements},
and those in shifted degrees are called \emph{edge elements}. Write
\[
 \begin{aligned}
 \pi_V^L(H)^{\mathrm{vertex}}
 &=\begin{cases}\pi_V^L(H),&\dim(V^L)\equiv0\pmod2,\\0,&\text{otherwise},\end{cases}\\
 \pi_V^L(H)^{\mathrm{edge}}
 &=\begin{cases}\pi_V^L(H),&\dim(V^L)\equiv1\pmod2,\\0,&\text{otherwise}.\end{cases}
 \end{aligned}
\]
A \emph{vertex-edge map} is a restriction or transfer map between
natural and shifted degrees.
For an edge array $z$ in degree $U$, write
$[z]\in\coker\Delta_U\cong\pi_{U-1}^G(H)$.
\end{definition}

For a natural degree $U\in RO_{\mathrm{sel}}(G)$, take
$A_i(U)$ and $T_i(U)$, including their bounds
$b_E(U),c_E(U)$ and zero groups, from
\eqref{eq:local-height-table}. Thus the degree coordinates of $U$ lie in $\mathbb Z^{n+1}$, or, at $p=2$,
in $(\mathbb Z+\tfrac12)^n\times\mathbb Z\bigcup\mathbb Z^{n+1}$.

Vertex-edge restriction and transfer maps can occur at $p=2$.

We start with the comparison maps.

\begin{proposition}
\label{prop:adjacent-local-comparison}
For $0\leq i<n-1$, the map $H^{\Phi_{i+1}}\to H^{t_i}$
factors in $\CAlg_H$ as
$H^{\Phi_{i+1}}\to H^{h_{i+1}}\xrightarrow{\rho_i}H^{t_i}$.
The maps $\rho_i$, $\rho_{n-1}:H^{\Phi_n}\to H^{t_{n-1}}$, and
$\ell_i$ induce the homomorphisms $[q]\mapsto[q]$ of
Theorem~\ref{thm:local-padic-comparisons}.
\end{proposition}
\begin{proof}
Apply Corollary~\ref{cor:hz-module-descent} and
Theorem~\ref{thm:local-padic-comparisons}.
\end{proof}

For the rest of this section, write
\[
 C^0=\left(\prod_{i=0}^{n-1}H^{h_i}\right)\times H^{\Phi_n},
 \qquad C^1=\prod_{i=0}^{n-1}H^{t_i},
\]
and let $d:H\to C^0$ be the canonical map and
$\Delta=(\ell_i-\rho_i)_i:C^0\to C^1$.

Use $\Delta_V$ from \eqref{eq:ch4-finite-difference} when $V$ is a natural degree.
For a shifted degree $V$ and $0\leq i<n$,
\[
 \pi_V^G(H^{h_i})=\pi_V^G(H^{t_i})=
 \pi_V^G(H^{\Phi_n})=0,
\]
so $\Delta_V:0\to0$. 
We can now get the homotopy groups.

\begin{theorem}[Prime-power homotopy groups]
\label{thm:prime-power-final}\label{thm:prime-power-array-answer}
\label{thm:ch4-literal-odd}\label{thm:integral-coordinate}
\label{thm:two-natural-ring}\label{thm:half-integral}
\label{core_thm}\label{thm:odd-natural-ring}
\label{thm:padic-homotopy-ring}
We have
\begin{equation}\label{eq:ch4-kernel-cokernel-computation}
 \pi_V^G(H)=\begin{cases}
 \ker\Delta_V,&\dim(V^G)\text{ even},\\
 \coker\Delta_{V+1},&\dim(V^G)\text{ odd}.
 \end{cases}
\end{equation}
For $V\in RO_{\mathrm{sel}}(G)$, we have:
\[
\begin{aligned}
 \ker\Delta_V
 &=\left\{(x_i)\in\prod_{i=0}^nA_i(V):
       \ell_i(x_i)=\rho_i(x_{i+1})\ (0\leq i<n)\right\},\\
 \coker\Delta_{V+1}
 &=\left(\prod_{i=0}^{n-1}T_i(V+1)\right)\big/\operatorname{im}\Delta_{V+1}.
\end{aligned}
\]
\end{theorem}
\begin{proof}
The homotopy groups of $H^{h_i}$, $H^{t_i}$, and $H^{\Phi_n}$, together with the comparison maps, are given by
Theorems~\ref{thm:local-padic-identification} and
\ref{thm:local-padic-comparisons}.
By Corollary~\ref{cor:hz-module-descent}, we have the fiber sequence
\begin{equation}\label{eq:ch4-finite-fibre}
 H\xrightarrow{d}C^0\xrightarrow{\Delta}C^1.
\end{equation}
By Theorem~\ref{thm:local-padic-identification}, we have
\[
 0\longrightarrow\coker\Delta_{V+1}
 \xrightarrow{\partial_*}\pi_V^G(H)
 \xrightarrow{d_*}\ker\Delta_V\longrightarrow0.
\]
By \eqref{eq:local-height-table},
\[
 \begin{aligned}
 \dim(V^G)\equiv0\pmod2
 &\ \Longrightarrow\ \coker\Delta_{V+1}=0,\quad
 d_*:\pi_V^G(H)\xrightarrow{\cong}\ker\Delta_V,\\
 \dim(V^G)\equiv1\pmod2
 &\ \Longrightarrow\ \ker\Delta_V=0,\quad
 \partial_*:\coker\Delta_{V+1}\xrightarrow{\cong}\pi_V^G(H).
 \end{aligned}
\]
Thus \eqref{eq:ch4-kernel-cokernel-computation} holds and
$F^1\pi_\star^G(H)=\ker d_*=
\pi_\star^G(H)^{\mathrm{edge}}$.
\end{proof}

\begin{remark}\label{prop:array-congruences}
For a natural degree $U=\rho_p(\mathbf a)$,
$\mathbf a\in D_n$, Theorem~\ref{thm:prime-power-final}
gives
\[
\pi_U^G(H)\cong
\left\{
(q_0,\ldots,q_n):
\begin{array}{l}
q_i\in p^{b_{A_i}(U)}\mathbb Z/
             p^{c_{A_i}(U)}\mathbb Z,\\
q_i-q_{i+1}\in p^{c_{T_i}(U)}\mathbb Z
\quad(0\leq i<n)
\end{array}
\right\},
\]
and
\[
\pi_{U-1}^G(H)\cong
\frac{
\left\{
(z_0,\ldots,z_{n-1}):
z_i\in p^{b_{T_i}(U)}\mathbb Z/
         p^{c_{T_i}(U)}\mathbb Z
\right\}
}{
\left\{
(q_0-q_1,\ldots,q_{n-1}-q_n):
q_i\in A_i(U)\ (0\leq i\leq n)
\right\}
}.
\]
All differences are taken in $T_i(U)$.
If $A_i(U)=0$, its component is zero.
The cokernel classes in the second quotient are written
$[z_0,\ldots,z_{n-1}]$.
\end{remark}

This also gives the finiteness of each group.

\begin{corollary}[Finiteness]\label{cor:coefficient-finiteness}
Each $\pi_V^G(H)$ is finitely generated, has free rank at most one
for a natural degree $V$, and is finite for a shifted degree $V$.
\end{corollary}
\begin{proof}
The groups $A_i(V)$ and $T_i(V)$ are finitely generated; only $A_0(V)$ can be infinite
cyclic, and every $T_i(V)$ is finite. Hence
$\operatorname{rank}_{\mathbb Z}\ker\Delta_V\leq1$, $|\coker\Delta_V|<\infty$.
\end{proof}
$C_4$-equivariant case have already computed in \cite{Yan}, using Euler class and orientation class as generators of polynomial algebra. The language translation is as follows.
\begin{example}[$C_4$ and $C_2$]
\label{ex:c4-c2-coefficients}
\emph{(1) $C_4$: vertex and edge arrays in every degree.}
Write $RO_\mathrm{sel}(C_4)=\mathbb Z\{1,\alpha,\lambda\}$, with
$\alpha=\sigma_4$ and $\lambda=\lambda_4(1)$.
For a natural degree, fix $(a_0,a_1,a_2)\in\mathbb Z^3\sqcup
((\mathbb Z+\tfrac12)^2\times\mathbb Z)$ and put
\[
 U=2a_2+2(a_1-a_2)\alpha+(a_0-a_1)\lambda.
\]
Thus $(a_0,a_1,a_2)=\tfrac12(\dim U,\dim U^{C_2},\dim U^{C_4})$.
Every degree is uniquely $U$ or $U-1$. Put
\[
 \epsilon=\begin{cases}0,&a_0,a_1\in\mathbb Z,\\1,&a_0,a_1\in\mathbb Z+\tfrac12.\end{cases}
\]
Here the groups $A_i(U)$ and $T_i(U)$ are
\[
\begin{array}{c|c|c}
 &\epsilon=0&\epsilon=1\\ \hline
 A_0(U)&
 \begin{cases}2^{a_1+a_2}\mathbb Z/2^{a_1+a_2+2}\mathbb Z,&a_0<0,\\
 2^{a_1+a_2}\mathbb Z,&a_0=0,\\0,&a_0>0\end{cases}&
 \begin{cases}2^{a_1+a_2+1/2}\mathbb Z/2^{a_1+a_2+3/2}\mathbb Z,&a_0<0,\\0,&a_0>0\end{cases}\\[4mm]
 A_1(U)&
 \begin{cases}2^{a_2}\mathbb Z/2^{a_2+1}\mathbb Z,&a_1<0,\\
 2^{a_1+a_2}\mathbb Z/2^{a_1+a_2+2}\mathbb Z,&a_1\geq0\end{cases}&
 2^{a_2+(a_1+1/2)^+}\mathbb Z/2^{a_2+(a_1+1/2)^++1}\mathbb Z\\[3mm]
 A_2(U)&\multicolumn{2}{c}{
 \begin{cases}2^{a_2}\mathbb Z/2^{a_2+1}\mathbb Z,&a_2\geq0,\\0,&a_2<0\end{cases}}\\[3mm]
 T_0(U)&2^{a_1+a_2}\mathbb Z/2^{a_1+a_2+2}\mathbb Z&2^{a_1+a_2+1/2}\mathbb Z/2^{a_1+a_2+3/2}\mathbb Z\\[1mm]
 T_1(U)&2^{a_2}\mathbb Z/2^{a_2+1}\mathbb Z&2^{a_2}\mathbb Z/2^{a_2+1}\mathbb Z
\end{array}
\]
The incidence map is
\[
 \Delta_U(q_0,q_1,q_2)=(q_0-q_1,q_1-q_2),\qquad
 \pi_U^{C_4}(H)=\ker\Delta_U,\quad
 \pi_{U-1}^{C_4}(H)=\coker\Delta_U.
\]
Unlisted cases are zero.

For $\epsilon=0$, the homotopy groups in natural degrees are
\[
\begin{array}{c|c|c|c|l}
 a_0&a_1&a_2&\pi_U^{C_4}(H)&\text{generators}\\ \hline
 0&a_1\geq1&\text{all}&\mathbb Z&(2^{a_1+a_2},2^{a_1+a_2},0)\\
 0&0&a_2\geq0&\mathbb Z&(2^{a_2},2^{a_2},2^{a_2})\\
 0&0&a_2<0&\mathbb Z&(2^{a_2+1},2^{a_2+1},0)\\
 0&-1&a_2\geq0&\mathbb Z&(2^{a_2},2^{a_2},2^{a_2})\\
 0&a_1\leq-2&a_2\geq0&\mathbb Z\oplus\mathbb Z/2&(2^{a_1+a_2+2},0,0),\ (0,2^{a_2},2^{a_2})\\
 0&a_1<0&a_2<0&\mathbb Z&(2^{a_1+a_2+2},0,0)\\ \hline
 a_0<0&a_1\geq1&\text{all}&\mathbb Z/4&(2^{a_1+a_2},2^{a_1+a_2},0)\\
 a_0<0&0&a_2\geq0&\mathbb Z/4&(2^{a_2},2^{a_2},2^{a_2})\\
 a_0<0&0&a_2<0&\mathbb Z/2&(2^{a_2+1},2^{a_2+1},0)\\
 a_0<0&a_1<0&a_2\geq0&\mathbb Z/2&(2^{a_2},2^{a_2},2^{a_2})\\
 a_0>0&a_1\leq-2&a_2\geq0&\mathbb Z/2&(0,2^{a_2},2^{a_2})
\end{array}
\]
For $\epsilon=0$, the homotopy groups in shifted degrees are
\[
\begin{array}{c|c|c|c|c}
 a_0&a_1&a_2&\pi_{U-1}^{C_4}(H)&\text{generator}\\ \hline
 a_0>0&a_1\leq-2&\text{all}&\mathbb Z/4&[2^{a_1+a_2},0]\\
 a_0>0&-1&a_2<0&\mathbb Z/4&[2^{a_1+a_2},0]\\
 a_0>0&-1&a_2\geq0&\mathbb Z/2&[2^{a_1+a_2},0]\\
 a_0>0&0&a_2<0&\mathbb Z/2&[0,2^{a_2}]\\
 \text{all}&a_1\geq1&a_2<0&\mathbb Z/2&[0,2^{a_2}]
\end{array}
\]
For $\epsilon=1$, all nonzero groups have order $2$:
\[
\begin{array}{c|c|c|c|c}
 &a_0&a_1&a_2&\text{generator}\\ \hline
 \pi_U^{C_4}(H)&\text{all}&a_1\leq-\tfrac32&a_2\geq0&(0,2^{a_2},2^{a_2})\\
 &a_0<0&-\tfrac12&a_2\geq0&(2^{a_2},2^{a_2},2^{a_2})\\
 &a_0<0&a_1\geq\tfrac12&\text{all}&(2^{a_1+a_2+1/2},2^{a_1+a_2+1/2},0)\\ \hline
 \pi_{U-1}^{C_4}(H)&a_0>0&a_1\leq-\tfrac32&\text{all}&[2^{a_1+a_2+1/2},0]\\
 &a_0>0&-\tfrac12&a_2<0&[0,2^{a_2}]\\
 &\text{all}&a_1\geq\tfrac12&a_2<0&[0,2^{a_2}]
\end{array}
\]
\emph{(2) $C_4$: comparison with Yan.}
The degrees of Yan's generators are
\[
\begin{gathered}
 |u_{2\alpha}|=2-2\alpha,\qquad |u_\lambda|=2-\lambda,\\
 |a_\alpha|=-\alpha,\qquad |a_\lambda|=-\lambda,\qquad
 |\tr_{C_2}^{C_4}(e_{3\alpha})|=3\alpha-3.
\end{gathered}
\]
Here $e_{3\alpha}=u_\alpha^{-3}\in\pi_{3\alpha-3}^{C_2}(H)$,
where $u_\alpha\in\pi_{1-\alpha}^{C_2}(H)$ is the orientation class
of $\alpha|_{C_2}\cong\mathbb R$.
Yan's additive decomposition \cite[Theorem~4.10]{Yan} is
\begingroup
\allowdisplaybreaks[1]
\begin{align*}
\pi_\star^{C_4}(H)={}&
 \mathbb Z[u_{2\alpha},u_\lambda,a_\alpha,a_\lambda]/
 (2a_\alpha,4a_\lambda,a_\alpha^2u_\lambda-2a_\lambda u_{2\alpha})\\
 &\oplus\mathbb Z\langle2u_{2\alpha}^{-i}:i\geq1\rangle\\
 &\oplus\mathbb Z\langle u_{2\alpha}^{-k}u_\lambda^i:i,k\geq1\rangle\\
 &\oplus\mathbb Z\langle2u_{2\alpha}^ku_\lambda^{-1}:k\geq1\rangle\\
 &\oplus\mathbb Z[u_{2\alpha}^{-1}]\langle4u_\lambda^{-j}:j\geq1\rangle\\
 &\oplus\mathbb Z\langle4u_{2\alpha}^ku_\lambda^{-j}:k\geq1,\ j\geq2\rangle\\
 &\oplus\mathbb Z/2[u_{2\alpha}^{\pm1}]
       \left\langle\Sigma^{-1}\frac{a_\alpha}{a_\lambda^iu_\lambda^j}:i,j\geq1\right\rangle\\
 &\oplus\mathbb Z/4[u_{2\alpha}^{\pm1}]
       \left\langle\Sigma^{-1}\frac{1}{a_\lambda^iu_\lambda^j}:i\geq1,\ j\geq2\right\rangle\\
 &\oplus\mathbb Z/4[u_{2\alpha}^{-1}]
       \left\langle\Sigma^{-1}\frac{1}{a_\lambda^iu_\lambda}:i\geq1\right\rangle\\
 &\oplus\mathbb Z/2
       \left\langle\Sigma^{-1}\frac{u_{2\alpha}^j}{a_\lambda^iu_\lambda}:i,j\geq1\right\rangle\\
 &\oplus\mathbb Z/2
       \left\langle\Sigma^{-1}\frac{a_\alpha}{a_\lambda^iu_{2\alpha}^j}:i,j\geq1\right\rangle\\
 &\oplus\mathbb Z/2
       \left\langle\Sigma^{-1}\frac{1}{a_\lambda^iu_{2\alpha}^j}:i,j\geq1\right\rangle\\
 &\oplus\mathbb Z/2[u_{2\alpha}]
       \langle a_\alpha^i a_\lambda^{-j}:i\geq3,\ j\geq1\rangle\\
 &\oplus\mathbb Z/2[u_{2\alpha}^{-1},a_\alpha^{-1},a_\lambda^{\pm1}]
       \langle\tr_{C_2}^{C_4}(e_{3\alpha})\rangle\\
 &\oplus\mathbb Z/2[a_\lambda]
       \left\langle\frac{a_\alpha u_\lambda^i}{u_{2\alpha}^j}:i,j\geq1\right\rangle\\
 &\oplus\mathbb Z/2\langle2u_{2\alpha}^{-i}:i\geq1\rangle
       \langle a_\lambda^j:j\geq1\rangle\\
 &\oplus\mathbb Z/4\langle u_{2\alpha}^{-k}u_\lambda^i:i,k\geq1\rangle
       \langle a_\lambda^j:j\geq1\rangle.
\end{align*}
\endgroup
The symbols $2u_{2\alpha}^{-i}$, $2u_{2\alpha}^ku_\lambda^{-1}$,
$4u_{2\alpha}^ku_\lambda^{-j}$ have Yan's meaning; they do not assert the existence
of the corresponding unmultiplied fractions in $\pi_\star^{C_4}(H)$.
Also, $|\Sigma^{-1}f|=|f|-1$.

For a vertex generator $f$, the following tables give the degree
coordinates $(a_0,a_1,a_2)$ of $U=|f|$ and its vertex array.
The positive cone is
\cite[Section~6.2]{Yan}
\[
\begin{aligned}
 &\mathbb Z[u_{2\alpha},u_\lambda]
 \oplus\mathbb Z/4[u_{2\alpha},u_\lambda]\langle a_\lambda^r:r\geq1\rangle\\
 &\quad\oplus\mathbb Z/2[u_{2\alpha},a_\lambda]\langle a_\alpha^r:r\geq1\rangle
 \oplus\mathbb Z/2[u_{2\alpha},a_\lambda]\langle u_\lambda^j:j\geq1\rangle\langle a_\alpha\rangle.
\end{aligned}
\]
Their correspondence is
\[
\begin{array}{c|c|c|c}
 f\quad\text{(range)}&(a_0,a_1,a_2)&\text{vertex array}&\text{order}\\ \hline
 u_{2\alpha}^iu_\lambda^j\ (i,j\geq0)&(0,j,i+j)&
 (2^{i+2j},2^{i+2j},\mathbf1_{j=0}2^i)&\infty\\
 u_{2\alpha}^iu_\lambda^ja_\lambda^r\ (i,j\geq0,\ r\geq1)&(-r,j,i+j)&
 (2^{i+2j},2^{i+2j},\mathbf1_{j=0}2^i)&4\\
 u_{2\alpha}^ia_\alpha^ra_\lambda^h\ (i,h\geq0,\ r\geq1)&
 (-h-\tfrac r2,-\tfrac r2,i)&(2^i,2^i,2^i)&2\\
 u_{2\alpha}^iu_\lambda^ja_\alpha a_\lambda^h\ (i,h\geq0,\ j\geq1)&
 (-h-\tfrac12,j-\tfrac12,i+j)&(2^{i+2j},2^{i+2j},0)&2
\end{array}
\]
The remaining free generators are
\[
\begin{array}{c|c|c}
 f\quad\text{(range)}&(a_0,a_1,a_2)&\text{vertex array}\\ \hline
 2u_{2\alpha}^{-i}\ (i\geq1)&(0,0,-i)&(2^{1-i},2^{1-i},0)\\
 u_{2\alpha}^{-k}u_\lambda^i\ (i,k\geq1)&(0,i,i-k)&(2^{2i-k},2^{2i-k},0)\\
 2u_{2\alpha}^ku_\lambda^{-1}\ (k\geq1)&(0,-1,k-1)&(2^{k-1},2^{k-1},2^{k-1})\\
 4u_{2\alpha}^{-k}u_\lambda^{-j}\ (k\geq0,\ j\geq1)&(0,-j,-k-j)&(2^{2-k-2j},0,0)\\
 4u_{2\alpha}^ku_\lambda^{-j}\ (k\geq1,\ j\geq2)&(0,-j,k-j)&(2^{k-2j+2},0,0)
\end{array}
\]
For the last two rows, use the generators in the image of
$\pi_\star^{C_4}(H_{h_0})\to\pi_\star^{C_4}(H)$, as in
\cite[Proposition~4.8]{Yan}; their $A_1$- and $A_2$-entries vanish.
The remaining vertex torsion generators are
\[
\begin{array}{c|c|c|c}
 f\quad\text{(range)}&(a_0,a_1,a_2)&\text{vertex array}&\text{order}\\ \hline
 u_{2\alpha}^ka_\alpha^ra_\lambda^{-j}\ (k\geq0,\ r\geq3,\ j\geq1)&
 (j-\tfrac r2,-\tfrac r2,k)&(0,2^k,2^k)&2\\
 a_\alpha u_\lambda^iu_{2\alpha}^{-j}a_\lambda^h\ (i,j\geq1,\ h\geq0)&
 (-h-\tfrac12,i-\tfrac12,i-j)&(2^{2i-j},2^{2i-j},0)&2\\
 (2u_{2\alpha}^{-i})a_\lambda^j\ (i,j\geq1)&
 (-j,0,-i)&(2^{1-i},2^{1-i},0)&2\\
 u_{2\alpha}^{-k}u_\lambda^ia_\lambda^j\ (i,j,k\geq1)&
 (-j,i,i-k)&(2^{2i-k},2^{2i-k},0)&4
\end{array}
\]
For an edge generator $f$, put $U=|f|+1$.
The following table gives the degree coordinates of $U$
and the corresponding edge array:
\[
\begin{array}{c|c|c|c}
 f\quad\text{(range)}&(a_0,a_1,a_2)&\text{edge array}&\text{order}\\ \hline
 \Sigma^{-1}\dfrac{u_{2\alpha}^ka_\alpha}{a_\lambda^iu_\lambda^j}
 \ (i,j\geq1,\ k\in\mathbb Z)&
 (i-\tfrac12,-j-\tfrac12,k-j)&[2^{k-2j},0]&2\\[2mm]
 \Sigma^{-1}\dfrac{u_{2\alpha}^k}{a_\lambda^iu_\lambda^j}
 \ (i\geq1,\ j\geq2,\ k\in\mathbb Z)&
 (i,-j,k-j)&[2^{k-2j},0]&4\\[2mm]
 \Sigma^{-1}\dfrac{u_{2\alpha}^{-k}}{a_\lambda^iu_\lambda}
 \ (i\geq1,\ k\geq0)&
 (i,-1,-k-1)&[2^{-k-2},0]&4\\[2mm]
 \Sigma^{-1}\dfrac{u_{2\alpha}^j}{a_\lambda^iu_\lambda}
 \ (i,j\geq1)&
 (i,-1,j-1)&[2^{j-2},0]&2\\[2mm]
 \Sigma^{-1}\dfrac{a_\alpha}{a_\lambda^iu_{2\alpha}^j}
 \ (i,j\geq1)&
 (i-\tfrac12,-\tfrac12,-j)&[2^{-j},0]&2\\[2mm]
 \Sigma^{-1}\dfrac{1}{a_\lambda^iu_{2\alpha}^j}
 \ (i,j\geq1)&
 (i,0,-j)&[2^{-j},0]&2\\[2mm]
 u_{2\alpha}^{-k}a_\alpha^{-r}a_\lambda^\ell\tr_{C_2}^{C_4}(e_{3\alpha})
 \ (k,r\geq0,\ \ell\in\mathbb Z)&
 (\tfrac{r+1}{2}-\ell,\tfrac{r+1}{2},-k-1)&[0,2^{-k-1}]&2
\end{array}
\]
The first six rows use the boundary from $H^{t_0}$, with compatible
choices of its sign. The last row uses the boundary from $H^{t_1}$;
in particular $\tr_{C_2}^{C_4}(e_{3\alpha})\leftrightarrow[0,\tfrac12]$.
For $a_0>0$, $a_2<0$, the rows with
$a_1=0$ or $a_1=-\tfrac12$ satisfy $[2^{a_2},0]=[0,2^{a_2}]$.

\emph{(3) $C_2$:}
Write $RO(C_2)=\mathbb Z\{1,\sigma\}$.
For $a_0\in\tfrac12\mathbb Z$ and $a_1\in\mathbb Z$, put
\[
 U=2a_1+2(a_0-a_1)\sigma.
\]
The groups $A_0(U)$, $A_1(U)$, $T_0(U)$, and the incidence map are
\[
\begin{aligned}
 A_0(U)&=\begin{cases}
 2^{a_1}\mathbb Z/2^{a_1+1}\mathbb Z,&a_0<0,\\
 2^{a_1}\mathbb Z,&a_0=0,\\0,&a_0>0,
 \end{cases}
 &A_1(U)&=\begin{cases}2^{a_1}\mathbb Z/2^{a_1+1}\mathbb Z,&a_1\geq0,\\0,&a_1<0,\end{cases}\\
 T_0(U)&=2^{a_1}\mathbb Z/2^{a_1+1}\mathbb Z,
 &\Delta_U(q_0,q_1)&=q_0-q_1.
\end{aligned}
\]
Consequently, all nonzero groups and their generators are
\[
\begin{array}{c|c|c|c|c}
 &a_0&a_1&\text{group}&\text{generator}\\ \hline
 \pi_U^{C_2}(H)&0&a_1\geq0&\mathbb Z&(2^{a_1},2^{a_1})\\
 &0&a_1<0&\mathbb Z&(2^{a_1+1},0)\\
 &a_0<0&a_1\geq0&\mathbb Z/2&(2^{a_1},2^{a_1})\\ \hline
 \pi_{U-1}^{C_2}(H)&a_0>0&a_1<0&\mathbb Z/2&[2^{a_1}]
\end{array}
\]
The classical computation is
\[
\begin{aligned}
 \pi_\star^{C_2}(H)={}&
 \mathbb Z[u_{2\sigma},a_\sigma]/(2a_\sigma)\\
 &\oplus\mathbb Z\langle2u_{2\sigma}^{-i}:i\geq1\rangle\\
 &\oplus\mathbb Z/2
 \langle\Sigma^{-1}u_{2\sigma}^{-i}a_\sigma^{-j}:i,j\geq1\rangle.
\end{aligned}
\]
With $|u_{2\sigma}|=2-2\sigma$ and $|a_\sigma|=-\sigma$,
the complete correspondence between degree coordinates and
vertex arrays or edge arrays is
\[
\begin{array}{c|c|c|c}
 f\quad\text{(range)}&(a_0,a_1)&\text{vertex array or edge array}&\text{order}\\ \hline
 u_{2\sigma}^i\ (i\geq0)&(0,i)&(2^i,2^i)&\infty\\
 u_{2\sigma}^ia_\sigma^j\ (i\geq0,\ j\geq1)&
 (-\tfrac j2,i)&(2^i,2^i)&2\\
 2u_{2\sigma}^{-i}\ (i\geq1)&(0,-i)&(2^{1-i},0)&\infty\\
 \Sigma^{-1}u_{2\sigma}^{-i}a_\sigma^{-j}\ (i,j\geq1)&
 (\tfrac j2,-i)&[2^{-i}]&2
\end{array}
\]
\end{example}

\subsection{Products and homotopy groups at subgroups}

The product is as follows.

\begin{proposition}[Products]\label{prop:array-products}
\label{prop:ch4-mixed-cokernel-product}\label{prop:complete-algebraic-forest}
For natural degrees $U,W$ and a shifted degree $V$, let $x\in\ker\Delta_U$,
$y\in\ker\Delta_W$, $[z]\in\coker\Delta_{V+1}$, and let $[w]$
be any edge element. Then
\begin{equation*}
\begin{aligned}
 (xy)_i&=x_i y_i\in A_i(U+W),\\
 x[z]&=\left[\bigl(\ell_i(x_i)z_i\bigr)_{i<n}\right]
             \in\coker\Delta_{U+V+1},\\
 [z][w]&=0.
\end{aligned}
\end{equation*}
Products are induced by multiplication of representatives; a zero
factor or target gives zero. The unit is $(1,\ldots,1)$.
\end{proposition}
\begin{proof}
Apply Theorem~\ref{thm:filtered-green-descent} and
\eqref{eq:ch4-kernel-cokernel-computation}.
Apply \eqref{eq:mixed-AW} and \eqref{eq:square-zero-descent}, using
$F^1\pi_\star^G(H)=\pi_\star^G(H)^{\mathrm{edge}}$.
\end{proof}

For a graded algebra $A$ and an $A$-bimodule $M$, write $A\ltimes M$
for $A\oplus M$ with
$(a,m)(b,n)=(ab,an+mb)$.
Consequently we have the following decomposition
\begin{equation*}
 \pi_\star^G(H)\cong
 \left(\bigoplus_{\dim(V^G)\text{ even}}\ker\Delta_V\right)
 \ltimes
 \left(\bigoplus_{\dim(V^G)\text{ odd}}\coker\Delta_{V+1}\right).
\end{equation*}
For every $L\leq G$ and $V\in RO(G)$,
\[
 F^1\pi_V^L(H)=\ker\bigl(d_*:\pi_V^L(H)\to\pi_V^L(C^0)\bigr)
 =\pi_V^L(H)^{\mathrm{edge}}.
\]
Naturality of $d$ makes $F^1\underline\pi_\star(H)$ a square-zero Green
ideal. At $p=2$, the vertex summands need not form a Mackey subfunctor;
the vertex-to-edge maps below imply this failure.

The formulas also give the values at smaller subgroups.

\begin{corollary}[Homotopy groups at subgroups]\label{cor:array-subgroup-values}
At $L=C_{p^m}\leq G$, $m\geq1$, apply Theorem~\ref{thm:prime-power-final} and
Proposition~\ref{prop:array-products} to
$V|_L$, replacing $n$ by $m$. At $e$, $\pi_{V|_e}^e(H)$
is $\mathbb Z$ when $\dim V=0$ and zero otherwise.
\end{corollary}
\begin{proof}
Apply the preceding formulas to $L$ and use $\res_e^G H\simeq H\mathbb Z$.
\end{proof}

\begin{remark}
Let $x\in\ker\Delta_V$ and $y\in\ker\Delta_W$.
For $A_i(V+W)\ne0$,
\[
 v_p((xy)_i)=
 \begin{cases}
 v_p(x_i)+v_p(y_i),&v_p(x_i)+v_p(y_i)<c_{A_i}(V+W),\\
 \infty,&v_p(x_i)+v_p(y_i)\ge c_{A_i}(V+W).
 \end{cases}
\]
A zero factor or target gives zero.
This formula follows from \eqref{eq:local-valuation-product}.
\end{remark}

\subsection{Restriction and transfer}

For $L=C_{p^m}<G$, $1\leq m<n$, the Tate compression diagram of $\res_L^G H$ has terms
$(\res_L^G H)^{h_i}$, $(\res_L^G H)^{t_i}$ $(0\leq i<m)$, and
$(\res_L^G H)^{\Phi_m}$ by \eqref{eq:layer-definitions}, using the
families $\{C_{p^j}:0\leq j\leq i\}$ of subgroups of $L$.
The restricted fiber sequence has
\[
 \begin{aligned}
 \res_L^G C^0&\simeq
 \left(\prod_{i<m}(\res_L^G H)^{h_i}\right)\times (\res_L^G H)^{\Phi_m},\\
 \res_L^G C^1&\simeq\prod_{i<m}(\res_L^G H)^{t_i},\qquad
 \res_L^G\Delta=(\ell_i-\rho_i)_{i<m}.
 \end{aligned}
\]
Indeed, $\res_L^G(H^{h_i})\simeq0$ for $i>m$,
$\res_L^G(H^{t_i})\simeq0$ for $i\geq m$, and
$\res_L^G(H^{h_m})\simeq\res_L^G(H^{\Phi_m})
\simeq(\res_L^G H)^{\Phi_m}$.

\leadbreak
We now get restriction and transfer.

\begin{theorem}\label{thm:adjacent-mackey}\label{thm:equal-parity}
Let $K=C_{p^{n-1}}<G$ and suppose
$\dim(V^G)\equiv\dim(V^K)\pmod2$.
For $J=G,K$, write
$\Delta_V^J:\pi_V^J(C^0)\longrightarrow\pi_V^J(C^1)$
for the map induced by $\Delta$.
The restriction and transfer maps on the homotopy groups of the vertex and edge terms in
Theorem~\ref{thm:local-forest-mackey} and
Proposition~\ref{prop:local-prime-base} satisfy
\[
 \Delta_V^K\res_K^G=\res_K^G\Delta_V^G,
 \qquad
 \Delta_V^G\tr_K^G=\tr_K^G\Delta_V^K.
\]
Taking kernels gives the vertex-to-vertex maps; taking cokernels
in degree $V+1$ gives the edge-to-edge maps.
\end{theorem}

\begin{proof}
Restriction and transfer commute with the homotopy long exact sequence
of \eqref{eq:ch4-finite-fibre}, including the connecting maps.
Then \eqref{eq:ch4-kernel-cokernel-computation} give the maps on kernels in degree $V$ and cokernels in degree $V+1$.
\end{proof}

We next record the maps between adjacent subgroups.

\begin{proposition}[Restriction and transfer]\label{prop:array-equal-parity}
Let $L=C_{p^m}$, $m\geq1$, and $V\in RO_{\mathrm{sel}}(G)$.
For $K=C_{p^{m-1}}<L$, the maps are
$\res_K^L:\pi_{V|_L}^L(H)\longrightarrow\pi_{V|_K}^K(H)$, $\tr_K^L:\pi_{V|_K}^K(H)\longrightarrow\pi_{V|_L}^L(H)$.
If $\dim(V^L)\equiv\dim(V^K)\pmod2$, put
$U=V|_L$ in natural degrees and $U=V|_L+1$ in shifted degrees.
Write $(u_0,\ldots,u_m)$ for the degree coordinates of $U$.
Between the corresponding nonzero homotopy groups in natural or shifted degrees, the maps are
\begin{equation*}
 \res_K^L([q])=[p^{-u_m}q],\qquad
 \tr_K^L([r])=[p^{u_m+1}r].
\end{equation*}
\end{proposition}
\begin{proof}
Apply Theorem~\ref{thm:equal-parity},
\eqref{eq:local-intrinsic-res-tr}, and
Proposition~\ref{prop:local-prime-base}.
\end{proof}

Let $K=C_{p^{n-1}}<G$.
For a vertex element use $U=V$; for an edge element use $U=V+1$.
At $K$, omit
$\pi_U^K(H^{\Phi_n})=\pi_U^K(H^{t_{n-1}})=0$.
On the remaining homotopy groups of the diagram terms, restriction and transfer are given by
the formulas in Proposition~\ref{prop:array-equal-parity} with
$m=n$, using the prefix $(u_0,\ldots,u_{n-1})$ of
the degree coordinates $(u_0,\ldots,u_n)$ of $U$; at $K=e$, use
Proposition~\ref{prop:local-prime-base}.

\subsection{Weyl action}

The next result gives the Weyl action.

\begin{proposition}[Weyl action]
\label{prop:adjacent-weyl}\label{prop:Weyl-action}
\label{rem:array-weyl}
Let $L<G$ be any proper subgroup and let $\gamma$ be the chosen generator
of $G$. For odd $p$, the Weyl action on $\pi_V^L(H)$ is
trivial. For $p=2$, write
$V=m+s\alpha+\sum_\lambda s_\lambda\lambda$,
where $\lambda$ ranges over the irreducible real rotation representations
and $\alpha$ is the sign representation.
Then the generator of $G/L$ acts by
\[
 c_\gamma\big|_{\pi_V^L(H)}=(-1)^s\operatorname{id}.
\]
Equivalently, $c_\gamma$ is $(-1)^m$ on vertex elements and
$(-1)^{m+1}$ on edge elements. The formulas apply also to $L=e$.
\end{proposition}

\begin{proof}
For $f\in[S^V,H]^L$,
$c_\gamma(f)=\gamma_H\circ f\circ\gamma_{S^V}^{-1}$.
The map $\gamma_H$ induces the identity on the constant homotopy
Mackey functor $\underline\pi_0(H)=\underline{\mathbb Z}$; hence
$\gamma_H\simeq_L\mathrm{id}$.
For $p=2$,
\[
 L\leq\ker\alpha,\qquad
 \gamma_{S^\alpha}\simeq_L-\mathrm{id},\qquad
 \gamma_{S^\lambda}\simeq_L\mathrm{id},
\]
where the last homotopy is through rotations. The sphere identities
extend to virtual summands, so
\[
 \gamma_{S^V}\simeq_L(-1)^s\mathrm{id},\qquad
 c_\gamma(f)=(-1)^sf,\qquad
 \dim(V^L)\equiv m+s\pmod2.
\]
For odd $p$, every nontrivial irreducible real representation is a
rotation, giving $c_\gamma=\mathrm{id}$.
\end{proof}

\subsection{Vertex-edge maps for $p=2$}
\label{sec:two-cross-parity}

Let $G=C_{2^n}$, $K=C_{2^{n-1}}$, and let $\alpha$ be the sign
representation of $G/K$. Fix the sign cofiber sequence
\[
 G/K_+\xrightarrow{p}S^0\xrightarrow{\iota}S^\alpha
 \xrightarrow{\delta}\Sigma G/K_+.
\]
Write $\iota_*$ for the map induced by $\iota$.
For $E=A_i$ or $T_i$ with
$E(U),E(U-\alpha)\ne0$,
\begin{equation*}
 \iota_*:E(U)\longrightarrow E(U-\alpha),\qquad [q]\longmapsto[q].
\end{equation*}
If either group is zero, the map is zero. By
\eqref{eq:local-height-table} and \eqref{eq:two-forest-generator}, including
the $H^{h_0}$ and $H^{\Phi_n}$ formulas in
Appendix~\ref{app:local-calculations},
\[
 -\alpha=(-\tfrac12,\ldots,-\tfrac12,0),\qquad
 b_E(-\alpha)=0,\qquad a_\alpha\longleftrightarrow[1].
\]
Multiplication by $[1]$ gives $\iota_*([q])=[q]$.
The remaining vertex-edge maps at $p=2$ are computed
below.
\begin{proposition}\label{prop:two-half-coset-kernel}
For an $H$-module $X\in\Sp^G$ and $V\in RO(G)$, we have
\[
 \ker\!\left(\res_K^G:\pi_V^G(X)\to\pi_V^K(X)\right)
 =\iota_*\pi_{V+\alpha}^G(X),
\]
\[
\begin{aligned}
 &\operatorname{im}\!\left(\tr_K^G:
       \pi_{V-1+\alpha}^K(X)\to\pi_{V-1+\alpha}^G(X)\right)\\
 &\hspace{2em}=\ker\!\left(\iota_*:\pi_{V-1+\alpha}^G(X)
                         \to\pi_{V-1}^G(X)\right),
\end{aligned}
\]
and
\[
\begin{aligned}
 \ker\!\left(\tr_K^G:\pi_V^K(X)\to\pi_V^G(X)\right)
 &=\operatorname{im}\!\left(\res_K^G:
       \pi_{V+1-\alpha}^G(X)\to\pi_V^K(X)\right).
\end{aligned}
\]
Since $\Res_K^G(1-\alpha)=0$, we have
$(V-1+\alpha)|_K=V|_K=(V+1-\alpha)|_K$.
\end{proposition}

\begin{proof}
The sign cofiber sequence, its rotation, and self-duality of $G/K_+$ give
\[
\begin{aligned}
 \pi_{V+\alpha}^G(X)&\xrightarrow{\iota_*}\pi_V^G(X)
 \xrightarrow{\res_K^G}\pi_V^K(X)
 =\pi_{V-1+\alpha}^K(X)\\
 &\xrightarrow{\tr_K^G}\pi_{V-1+\alpha}^G(X)
 \xrightarrow{\iota_*}\pi_{V-1}^G(X).
\end{aligned}
\]
Use the standard orientation $\alpha|_K\cong\mathbb R$.
Exactness gives the first two identities.
Replacing $V$ by $V+1-\alpha$ gives the third.
\end{proof}

For $J=G,K$, use
$\Delta_V^J:\pi_V^J(C^0)\to\pi_V^J(C^1)$ induced by $\Delta$,
with vertex indices $0,\ldots,n$ and edge indices $0,\ldots,n-1$.
Attention
$\pi_V^K(H^{\Phi_n})=\pi_V^K(H^{t_{n-1}})=0$
and the zero maps to or from these zero groups. We can apply $\iota_*,\res_K^G,\tr_K^G$ separately
to each vertex entry and each edge entry.

We have the following formula for restriction.

\begin{proposition}[Vertex-to-edge restriction]
\label{prop:two-natural-shifted-res}
Let $V\in(\mathbb Z+\tfrac12)^n\times\mathbb Z$ and
$x\in\ker\Delta_V^G=\pi_V^G(H)^{\mathrm{vertex}}$.
Choose a vertex array $b\in\pi_{V+\alpha}^G(C^0)$ with
$\iota_* b=x$.
There is a unique edge array
$c\in\prod_{i<n}\pi_{V+\alpha}^K(H^{t_i})
=\prod_{i<n}\pi_{V+1}^K(H^{t_i})$ such that
\[
 \tr_K^G(c)=\Delta_{V+\alpha}^G(b).
\]
Then
\begin{equation*}
 \res_K^G(x)
 =[c]\in\coker\Delta_{V+1}^K
       =\pi_V^K(H)^{\mathrm{edge}}.
\end{equation*}
The class $[c]$ is independent of the chosen vertex lift $b$.
\end{proposition}

\begin{proof}
For $E=C^0,C^1$, \eqref{eq:local-height-table} gives
$\pi_{V+1}^G(E)=\pi_V^K(E)=0$. The sign exact sequence becomes
\[
 0\longrightarrow\pi_{V+\alpha}^K(E)
 \xrightarrow{\tr_K^G}\pi_{V+\alpha}^G(E)
 \xrightarrow{\iota_*}\pi_V^G(E)\longrightarrow0.
\]
Choose $b$ by surjectivity for $C^0$. Then
\[
 \iota_*\Delta_{V+\alpha}^G(b)=\Delta_V^G(x)=0
 \ \Longrightarrow\
 \exists!\,c:\tr_K^G(c)=\Delta_{V+\alpha}^G(b).
\]
For another lift $b'$, exactness gives
\[
 b'-b=\tr_K^G(e),\quad e\in\pi_{V+\alpha}^K(C^0)
 \ \Longrightarrow\
 c'-c=\Delta_{V+\alpha}^K(e)=\Delta_{V+1}^K(e).
\]
Thus $[c]$ is independent of $b$.

Let $\partial:\Sigma^{-1}C^1\to H$ be the boundary of
\eqref{eq:ch4-finite-fibre}. Regard
$x:S^V\to H$ and $b:S^{V+\alpha}\to C^0$ as representatives with
$b(1\wedge\iota)\simeq dx$.
The morphism-of-triangles axiom (TR3), applied to the rotated sign
triangle and $H\to C^0\to C^1\to\Sigma H$, gives a map
$c'':S^V\wedge\Sigma G/K_+\longrightarrow C^1$ satisfying
\[
 \begin{cases}
 c''(1\wedge\delta)=\Delta b,\\
 \Sigma\partial\,c''=-\Sigma x\,\Sigma(1\wedge p).
 \end{cases}
\]
By the induction--restriction adjunction,
\[
[S^V\wedge G/K_+,H]^G
\cong [S^{V|_K},H|_K]^K,
\qquad
[S^V\wedge\Sigma G/K_+,C^1]^G
\cong [S^{(V+\alpha)|_K},C^1|_K]^K,
\]
where $\alpha|_K$ is the trivial one-dimensional
representation. Under these adjunction isomorphisms,
precomposition with $1\wedge p$ gives $\res_K^G$,
while precomposition with $1\wedge\delta$ gives
$\tr_K^G$, up to the chosen orientation sign.
Hence
\[
 \tr_K^G(c'')=\pm\Delta_{V+\alpha}^G(b)
 \ \Longrightarrow\ c''=\pm c,\qquad
 \partial_*c=\pm\res_K^G(x).
\]
By \eqref{eq:local-height-table}, $2\pi_V^G(C^0)=0$; injectivity of
$d_*$ gives $2x=0$. Therefore
\[
 \res_K^G(x)=\partial_*c=[c].
\]
\end{proof}

\leadbreak
We next give the corresponding transfer formula.

\begin{proposition}[Vertex-to-edge transfer]
\label{prop:two-natural-shifted-transfer}
Let $V\in\mathbb Z^n\times(\mathbb Z+\tfrac12)$ and
$x\in\ker\Delta_V^K=\pi_V^K(H)^{\mathrm{vertex}}$.
Choose a vertex array $b\in\pi_{V+1-\alpha}^G(C^0)$ with
$\res_K^G(b)=x$.
There is a unique edge array
$c\in\prod_{i<n}\pi_{V+1}^G(H^{t_i})$ such that
\[
 \iota_* c=\Delta_{V+1-\alpha}^G(b).
\]
Then
\begin{equation*}
 \tr_K^G(x)
 =[c]\in\coker\Delta_{V+1}^G
       =\pi_V^G(H)^{\mathrm{edge}}.
\end{equation*}
The class $[c]$ is independent of the chosen vertex lift $b$.
\end{proposition}
\begin{proof}
For $E=C^0,C^1$, \eqref{eq:local-height-table} gives
$\pi_{V+1}^K(E)=\pi_V^G(E)=0$. The sign exact sequence becomes
\[
 0\longrightarrow\pi_{V+1}^G(E)
 \xrightarrow{\iota_*}\pi_{V+1-\alpha}^G(E)
 \xrightarrow{\res_K^G}\pi_V^K(E)\longrightarrow0.
\]
Choose $b$ by surjectivity for $C^0$. Then
\[
 \res_K^G\Delta_{V+1-\alpha}^G(b)=\Delta_V^K(x)=0
 \ \Longrightarrow\
 \exists!\,c:\iota_*c=\Delta_{V+1-\alpha}^G(b).
\]
For another lift $b'$, exactness gives
\[
 b'-b=\iota_*e,\quad e\in\pi_{V+1}^G(C^0)
 \ \Longrightarrow\ c'-c=\Delta_{V+1}^G(e).
\]
Thus $[c]$ is independent of $b$.

Since $(V+1-\alpha)|_K=V|_K$, choose representatives
$x:S^{V+1-\alpha}\wedge G/K_+\to H$ and
$b:S^{V+1-\alpha}\to C^0$ with $b(1\wedge p)\simeq dx$.
Apply TR3 to $S^{V+1-\alpha}$ smashed with the sign triangle and
$H\xrightarrow{d}C^0\xrightarrow{\Delta}C^1
\xrightarrow{\Sigma\partial}\Sigma H$. The resulting map
$c'':S^{V+1}\to C^1$ satisfies
\[
 c''(1\wedge\iota)=\Delta b,\qquad
 \Sigma\partial\,c''=\Sigma x(1\wedge\delta).
\]
Injectivity of $\iota_*$ gives $c''=c$.
The induction--restriction adjunction gives
$\partial_*c=\pm\tr_K^G(x)$.
Since $V+1\in D_n^{(1)}$, \eqref{eq:local-height-table} gives
$2\pi_{V+1}^G(C^1)=0$.
Thus $2\partial_*c=0$, so
\[
 \tr_K^G(x)=\partial_*c=[c].
\]
\end{proof}

This gives the following isomorphism.

\begin{corollary}\label{cor:cross-transfer-quotient}
Under Proposition~\ref{prop:two-natural-shifted-transfer}, transfer
induces the isomorphism
\[
\begin{gathered}
 \frac{\pi_V^K(H)^{\mathrm{vertex}}}
 {\operatorname{im}\!\left(\res_K^G:
      \pi_{V+1-\alpha}^G(H)^{\mathrm{vertex}}
           \to\pi_V^K(H)^{\mathrm{vertex}}\right)}\\[4pt]
 \xrightarrow{\ \cong\ }
 \ker\!\left(\iota_*:
     \pi_V^G(H)^{\mathrm{edge}}
        \to\pi_{V-\alpha}^G(H)^{\mathrm{edge}}\right).
\end{gathered}
\]
\end{corollary}
\begin{proof}
Apply Proposition~\ref{prop:two-half-coset-kernel}.
\end{proof}
We also have the following result.

\begin{proposition}
\label{thm:two-half-coset-transfer}\label{prop:two-shifted-natural-res}
The edge-to-vertex maps are zero.
For $V\in(\mathbb Z+\tfrac12)^n\times\mathbb Z$,
\[
 \tr_K^G:\pi_V^K(H)^{\mathrm{edge}}
 \longrightarrow\pi_V^G(H)^{\mathrm{vertex}},\qquad x\longmapsto0,
\]
and for $V\in\mathbb Z^n\times(\mathbb Z+\tfrac12)$,
\[
 \res_K^G:\pi_V^G(H)^{\mathrm{edge}}
 \longrightarrow\pi_V^K(H)^{\mathrm{vertex}},\qquad y\longmapsto0.
\]
\end{proposition}

\begin{proof}
For $J=G,K$, the vertex map $d_*:\pi_V^J(H)\to\pi_V^J(C^0)$
vanishes on edge elements and is injective on vertex elements.
And we have
\[
 d_*\tr_K^G(x)=\tr_K^G d_*(x)=0,
 \qquad
 d_*\res_K^G(y)=\res_K^G d_*(y)=0.
\]
Injectivity of $d_*$ in the target degrees proves both statements.
\end{proof}

\begin{theorem}[Vertex-edge maps for $p=2$]\label{thm:cross-maps}
For every adjacent inclusion $C_{2^{m-1}}<C_{2^m}$,
vertex-to-edge restriction and transfer are given by
Propositions~\ref{prop:two-natural-shifted-res} and
\ref{prop:two-natural-shifted-transfer}, respectively.
The edge-to-vertex maps are zero by
Proposition~\ref{thm:two-half-coset-transfer}.
\end{theorem}
\begin{proof}
Apply the three propositions to
$C_{2^{m-1}}<C_{2^m}$ in degree $V|_{C_{2^m}}$.
Here $\alpha$ is the sign representation of $C_{2^m}$, so
$\alpha|_{C_{2^{m-1}}}=1$.
\end{proof}

We now extend the computation from $RO_\mathrm{sel}(G)$ to $RO(G)$ using $v_D$.

\begin{corollary}\label{cor:array-character-transport}
For $V\in RO(G)$, use the decomposition
$V=V_0+D$ with $V_0\in RO_{\mathrm{sel}}(G)$
from Section~\ref{sec:oriented-models}
and $v_D$ as in Lemma~\ref{lem:normalized-character-units}.
Then, for $y\in\pi_{V_0|_K}^K(H)$, we have
\[
 \begin{gathered}
 \pi_{V_0|_L}^L(H)\xrightarrow{\ \res_L^G(v_D)\cdot-\ }
 \pi_{V|_L}^L(H),\\
 \tr_K^L\bigl(\res_K^G(v_D)y\bigr)
       =\res_L^G(v_D)\tr_K^L(y).
 \end{gathered}
\]
These isomorphisms are compatible with products, restriction, and Weyl action.
\end{corollary}
\begin{proof}
Apply Lemma~\ref{lem:normalized-character-units} and Frobenius reciprocity.
\end{proof}

\section{Assemble to cyclic group}\label{sec:sylow-assembly}
In this section, we compute the
$RO(G)$-graded Green functor of $H\underline{\mathbb Z}$
for arbitrary finite cyclic groups $G$.
More precisely, motivated by Corollary~\ref{cor:cyclic-sylow-reduction},
we assemble the computations for
$C_{p^n}\leq G$.

Let $N=\prod_{p\mid N}p^{e_p}>1$, $G=C_N$, and
$P_p=C_{p^{e_p}}$. All primes below divide $N$.
For $K\leq G$, put $K_p=K\cap P_p$, so $K=\prod_{p\mid N}K_p$.
Section~\ref{sec:reconstruction} computes $\pi_V^{K_p}(H)$ and the
Green functor structure maps for each $P_p$. At $e$,
\[
 \pi_V^e(H)=\begin{cases}\mathbb Z,&\dim V=0,\\0,&\dim V\ne0.\end{cases}
\]
For $V=V^+-V^-\in RO(G)$, define the orientation character
$\varepsilon_V\colon G\to\{1,-1\}$ by
\[
\varepsilon_V(g)
=\det(g|V^+)\det(g|V^-)^{-1}.
\]
Then, for every $g\in G$,
\[
c_g|_{\pi_V^e(H)}
=\varepsilon_V(g)\operatorname{id},
\qquad
\varepsilon_{U+V}=\varepsilon_U\varepsilon_V.
\]

We first compute the Weyl actions.

\begin{lemma}\label{lem:complement-action-cyclic}\label{lem:complementary-weyl}
For $L\leq P_p$ and $h\in G$ of order prime to $p$,
\[
 c_h|\pi_V^L(H)=\varepsilon_V(h)\operatorname{id}.
\]
For $g\in G$, write $g_p\in P_p$ for its $p$-primary component. Then
\begin{equation}\label{eq:cyclic-total-weyl-action}
 c_g|\pi_V^{K_p}(H)=\varepsilon_V(gg_p^{-1})c_{g_p}.
\end{equation}
\end{lemma}
\begin{proof}
Since $\underline{\mathbb Z}$ is the constant Mackey functor,
its Weyl action is trivial. Thus $c_h$ is determined by
the action of $h$ on the representation sphere. On a rotation plane $\lambda$, a
path of rotations gives $h|_{S^\lambda}\simeq_L\operatorname{id}$;
$h|_{\mathbb R}=1$. If $N$ is even, the sign line $\sigma_N$ satisfies
\[
 \begin{aligned}
 p=2&\ \Longrightarrow\ h|_{\sigma_N}=1
       \quad\text{(since $h$ has odd order)},\\
 p\ne2&\ \Longrightarrow\ \Res_L^G\sigma_N=\mathbb R.
 \end{aligned}
\]
If $h|_{\sigma_N}=-1$, then $p\ne2$ and the induced map on
$S^{\Res_L^G\sigma_N}\cong S^1$ has degree $-1$.
Smash products and dual spheres give
$c_h=\varepsilon_V(h)\operatorname{id}$.
Finally, $g=g_p(gg_p^{-1})$ gives
\eqref{eq:cyclic-total-weyl-action}.
\end{proof}

Since $\varepsilon_V$ is trivial on odd-order subgroups,
\[
 \varepsilon_V|_{\prod_{q\ne p}K_q}\ne1
 \quad\Longleftrightarrow\quad
 p\ne2\ \text{and}\ \varepsilon_V|_K\ne1.
\]
We now assemble the results for $C_{p^n}$.

\begin{theorem}[Cyclic Green-functor assembly]
\label{thm:integral-cyclic-reconstruction}\label{thm:integral-cyclic-sylow}
\label{thm:global-main}\label{cor:arbitrary-cyclic-final}
For $K\leq C_N$ and $V\in RO(C_N)$, restriction induces an isomorphism
\[
 \pi_V^K(H)\xrightarrow{\cong}
 \left\{(x_p)_p\in\prod_{p\mid N}\pi_V^{K_p}(H):
 \begin{array}{l}
 \res_e^{K_p}x_p=\res_e^{K_q}x_q\quad(p,q\mid N),\\
 p\ne2,\ \varepsilon_V|_K\ne1\ \Longrightarrow\ x_p=0
 \end{array}\right\},
 \qquad x\longmapsto(\res_{K_p}^Kx)_p.
\]
Choose integers $b_p$ with
\begin{equation*}
 \sum_{p\mid N}b_p[K:K_p]=1.
\end{equation*}
The inverse is
\begin{equation*}
 (x_p)_p\longmapsto\sum_{p\mid N}b_p\tr_{K_p}^K(x_p).
\end{equation*}
\end{theorem}
\begin{proof}
\begin{equation}\label{eq:cyclic-invariants}
 \left(\pi_V^{K_p}(H)\right)^{\prod_{q\ne p}K_q}=
 \begin{cases}
 \pi_V^{K_p}(H),&p=2\ \text{or}\ \varepsilon_V|_K=1,\\
 0,&p\ne2\ \text{and}\ \varepsilon_V|_K\ne1.
 \end{cases}
\end{equation}
The first case follows from Lemma~\ref{lem:complement-action-cyclic}.
In the second case, every invariant element
$z\in \left(\pi_V^{K_p}(H)\right)^{ \prod_{q\ne p}K_q}$
satisfies $2z=0$.
Torsion-freeness of $\pi_V^e(H)$ and
Proposition~\ref{prop:cohomological-hz} give
\[
 \begin{gathered}
 2z=0\ \Longrightarrow\ \res_e^{K_p}z=0
 \ \Longrightarrow\ |K_p|z=\tr_e^{K_p}\res_e^{K_p}z=0,\\
 \gcd(2,|K_p|)=1\ \Longrightarrow\ z=0.
 \end{gathered}
\]
For $x\in\pi_V^K(H)$, put $x_p=\res_{K_p}^Kx$.
By transitivity,
$\res_e^{K_p}x_p=\res_e^Kx$ for every $p$.
Since $c_kx_p=x_p$ for every $k\in K$,
\eqref{eq:cyclic-invariants} gives $x_p=0$
whenever $p\ne2$ and $\varepsilon_V|_K\ne1$.

Since
$[K:K_p]=\prod_{q\ne p}|K_q|$, $\gcd_p[K:K_p]=1$,
choose $b_p\in\mathbb Z$ with $\sum_p b_p[K:K_p]=1$. For
$x\in\pi_V^K(H)$, cohomologicality gives
\[
 \sum_p b_p\tr_{K_p}^K\res_{K_p}^Kx
 =\left(\sum_p b_p[K:K_p]\right)x=x.
\]
For a compatible family $(x_p)_p$ satisfying the theorem's conditions,
\eqref{eq:cyclic-invariants} and
$c_k|_{\pi_V^{K_p}(H)}=\operatorname{id}$ for $k\in K_p$ give
$c_kx_p=x_p\;(k\in K)$.
Since $K_p\cap K_q=e$ for $p\ne q$, the double-coset formula and
cohomologicality give
\[
 \res_{K_p}^K\tr_{K_q}^K(x_q)
 =\begin{cases}
 [K:K_p]x_p,&q=p,\\
 [K:K_pK_q]\tr_e^{K_p}\res_e^{K_q}x_q=[K:K_q]x_p,&q\ne p.
 \end{cases}
\]
Here $\res_e^{K_q}x_q=\res_e^{K_p}x_p$. Hence
\[
 \res_{K_p}^K
       \left(\sum_q b_q\tr_{K_q}^Kx_q\right)
 =\sum_q b_q[K:K_q]x_p=x_p.
\]
This proves surjectivity and the third formula in the theorem.
\end{proof}

We now express the Green functor operations in terms of $(x_p)_p$.

\begin{proposition}\label{prop:cyclic-green-operations}
Under Theorem~\ref{thm:integral-cyclic-reconstruction},
addition is entrywise. For $L\leq K$, the unit, products,
restrictions, transfers, and Weyl actions are
\begin{align*}
 1_K&=(1_{K_p})_p,\\
 (xy)_p&=x_py_p,\\
 (\res_L^Kx)_p
 &=\res_{L_p}^{K_p}(x_p),\\
 (\tr_L^Ky)_p
 &=\begin{cases}
 \displaystyle\left(\prod_{q\ne p}[K_q:L_q]\right)
   \tr_{L_p}^{K_p}(y_p),&p=2\ \text{or}\ \varepsilon_V|_K=1,\\
 0,&p\ne2\ \text{and}\ \varepsilon_V|_K\ne1,
 \end{cases}\\
 (c_gx)_p&=\varepsilon_V(gg_p^{-1})c_{g_p}(x_p).
\end{align*}
The domains and codomains of above formulas are
\[
 \begin{gathered}
 \pi_U^K(H)\times\pi_V^K(H)\xrightarrow{\ \cdot\ }\pi_{U+V}^K(H),\qquad
 \res_L^K:\pi_V^K(H)\longrightarrow\pi_V^L(H),\\
 \tr_L^K:\pi_V^L(H)\longrightarrow\pi_V^K(H),\qquad
 c_g:\pi_V^K(H)\longrightarrow\pi_V^K(H).
 \end{gathered}
\]
\end{proposition}
\begin{proof}
Multiplicativity and transitivity of restriction give
the first three formulas in the proposition;
the fifth formula follows from
\eqref{eq:cyclic-total-weyl-action}. For transfer, if $p\ne2$ and
$\varepsilon_V|_K\ne1$, then $(\tr_L^K y)_p=0$ by
\eqref{eq:cyclic-invariants}. Otherwise,
\[
 K_p\backslash K/L\cong\prod_{q\ne p}K_q/L_q,
 \qquad K_p\cap L=L_p.
\]
Choose representatives $t\in\prod_{q\ne p}K_q$
for the double cosets in $K_p\backslash K/L$.
For every such $t$,
$c_t=\operatorname{id}$ on $\pi_V^{L_p}(H)$. Thus
\[
 \res_{K_p}^K\tr_L^K(y)
 =\sum_{[t]\in\prod_{q\ne p}K_q/L_q}\tr_{L_p}^{K_p}
        c_t\res_{L_p}^L(y)
 =\left(\prod_{q\ne p}[K_q:L_q]\right)
   \tr_{L_p}^{K_p}(y_p).
\]
\end{proof}

\begin{remark}
For $N=1$, the only nonzero group is $\pi_0^e(H)=\mathbb Z$.
\end{remark}
\section{Graded norms}\label{sec:tambara}

The preceding sections determine the $RO(G)$-graded Green
functor of $H\underline{\mathbb Z}$.
We now compute its norms. The explicit formula of norm is given in Theorem~\ref{thm:tambara-natural-kernel} and Theorem~\ref{thm:tambara-shifted-exact}. Also we get some the norm identities and prove an abstract result motivated from our computation.

\subsection{Induction of degrees}

Set $RO_{\mathrm{sel}}(e)=RO(e)=\mathbb Z\{1\}$.
Let $K=C_{p^k}\leq L=C_{p^h}$, where $0\leq k\leq h$,
and let $V\in RO_{\mathrm{sel}}(K)$.
Recall from Definition~\ref{def:degree-coordinates}
that the coordinates of a natural degree are half its
fixed-point dimensions along the subgroup chain.
We use the same convention for the source degree $V$
and the induced target degree $\Ind_K^L V$.

Write $\mathbf a=(a_0,\ldots,a_k)$ for the coordinates of $V$
and $\mathbf A=(A_0,\ldots,A_h)$ for those of $\Ind_K^L V$,
where
\[
\begin{aligned}
2a_j&=\dim(V^{C_{p^j}})
&& (0\leq j\leq k),\\
2A_j&=\dim\bigl((\Ind_K^L V)^{C_{p^j}}\bigr)
&& (0\leq j\leq h).
\end{aligned}
\]
These coordinates are defined by fixed-point dimensions
even when $\Ind_K^L V$ does not lie in
$RO_{\mathrm{sel}}(L)$.
The following proposition expresses the target coordinates
$\mathbf A$ in terms of the source coordinates $\mathbf a$.

\begin{proposition}
\label{prop:tambara-induced-degree}
The target coordinates are
\begin{equation*}
 A_j=
 \begin{cases}
 p^{h-k}a_j,&0\leq j\leq k,\\[2mm]
 p^{h-j}a_k,&k\leq j\leq h.
 \end{cases}
\end{equation*}
In particular, for the inclusion
$K=C_{p^k}<G=C_{p^{k+1}}$,
\begin{equation}\label{eq:tambara-adjacent-degree}
 I(\mathbf a)
 =(pa_0,\ldots,pa_k,a_k).
\end{equation}
In particular, for $p=2$, if the source degree has coordinates $\mathbf a\in(\mathbb Z+\tfrac12)^k\times\mathbb Z$, then $I(\mathbf a)\in\mathbb Z^{k+2}$.
\end{proposition}

\begin{proof}
For $H=C_{p^j}\leq L$, the double-coset formula gives
\[
\dim\bigl((\Ind_K^L V)^H\bigr)
=[L:KH]\dim(V^{K\cap H}),
\]
since $L$ is abelian.
If $0\leq j\leq k$, then $KH=K$ and $K\cap H=H$,
so $A_j=p^{h-k}a_j$.
If $k\leq j\leq h$, then $KH=H$ and $K\cap H=K$,
so $A_j=p^{h-j}a_k$.
Taking $h=k+1$ gives \ref{eq:tambara-adjacent-degree}.
For $p=2$, the final assertion follows from
$2a_j\in\mathbb Z$ for $0\leq j<k$
and $a_k\in\mathbb Z$.
\end{proof}

Also recall that for $D\in RO(K)$ as in Lemma~\ref{lem:normalized-character-units},
Section~\ref{sec:oriented-models} give
\begin{equation*}
 N_K^L(v_D)=v_{\Ind_K^L D}.
\end{equation*}
In general, $\Ind_K^L(RO_{\mathrm{sel}}(K))\not\subseteq
RO_{\mathrm{sel}}(L)$; hence the norm formulas below in valuation
coordinates should use $v_D$ to transport from the corresponding selected degree to $\Ind_K^L V$.

\subsection{Norms are compatible with the Tate compression diagram}

To compute the norms using Tate compression, we first show that the norm is compatible with the Tate compression diagram. We use a more catgorical language.

\begin{definition}
Let $G$ be a finite group.
An $N_\infty$ operad is a $G$-operad $\mathcal O$ such that
$\mathcal O(0)$ is $G$-contractible and, for each $d\geq0$,
$\mathcal O(d)$ is $\Sigma_d$-free and universal for a family of
subgroups of $G\times\Sigma_d$ containing every $K\times\{1\}$ with $K\leq G$
\cite[Definition~3.7]{BlumbergHill}.
For a finite $K$-set $T$, let $\sigma_T:K\to\Sigma_{|T|}$ be its action
and $\Gamma_T=\{(k,\sigma_T(k)):k\in K\}$.
The set $T$ is admissible if $\mathcal O(|T|)^{\Gamma_T}\simeq *$.
Write $\operatorname{Alg}_{\mathcal O}$ for the category of
$\mathcal O$-algebras in genuine $G$-spectra.
\end{definition}

We use the $N_\infty$ operad \cite{BlumbergHill} for
$N_{C_{p^{n-1}}}^{C_{p^n}}$ on the terms of the Tate compression diagram.

For $n\ge2$, let $\mathcal I$ be the incidence category of
Section~\ref{sec:incidence-category}, with $n$ edges. Write
$\mathcal I_<$, the truncated version, to be its full subcategory with objects
\[
 v_0,\ldots,v_{n-1},\ e_0,\ldots,e_{n-2},
 \qquad v_i\longrightarrow e_i\longleftarrow v_{i+1}
 \quad(0\le i<n-1).
\]
The diagram on $\mathcal I_<$ has values
$\mathcal R(v_i)=H^{h_i}$ and $\mathcal R(e_i)=H^{t_i}$,
with arrows $\ell_i$ and $\rho_i$.

The following lemma essentially says that norm is compatible with the Tate compression diagram of $H$, we prove it in $N_\infty$ operad language.

\begin{lemma}
\label{lem:tambara-local-admissibility}
Let $n\geq2$, $G=C_{p^n}$, and $K=C_{p^{n-1}}$.
Consider the truncated Tate compression diagram
\[
\mathcal R:\quad
H^{h_0}\xrightarrow{\ell_0}H^{t_0}
\xleftarrow{\rho_0}H^{h_1}
\longrightarrow\cdots\longrightarrow
H^{t_{n-2}}\xleftarrow{\rho_{n-2}}H^{h_{n-1}}.
\]
Then there is an $N_\infty$ operad $\mathcal O$ with
$\mathcal O(p)^{\Gamma_{G/K}}\simeq *$
such that $\mathcal R$ and its cone from $H$ take values in
$\operatorname{Alg}_{\mathcal O}$.

Moreover, there is an equivalence of cones
\[
\begin{aligned}
\bigl(\res_K^G H\to\res_K^G\mathcal R\bigr)
\simeq
\Bigl(
\res_K^G H\to
\bigl[
&(\res_K^G H)^{h_0}
\xrightarrow{\ell_0}(\res_K^G H)^{t_0}
\xleftarrow{\rho_0}(\res_K^G H)^{h_1}\\
&\longrightarrow\cdots\longrightarrow
(\res_K^G H)^{t_{n-2}}
\xleftarrow{\rho_{n-2}}
(\res_K^G H)^{\Phi_{n-1}}
\bigr]
\Bigr).
\end{aligned}
\]
In particular, the restricted cone induces an equivalence
\[
\res_K^G H\xrightarrow{\ \simeq\ }
\holim_{\mathcal I_<}\res_K^G\mathcal R.
\]
\end{lemma}

\begin{proof}
We use orthogonal $G$-spectra on a complete universe,
a commutative model of $H$, and derived algebra categories
and mapping spaces.
Let $\mathcal E$ be an ordinary $E_\infty$ operad with trivial
$G$-action. By \cite[Proposition~5.6 and Corollary~5.7]{Carrick},
a cofibrant model
\[
\mathcal O\simeq\operatorname{Map}(E\mathcal F_{n-1},\mathcal E)
\]
is an $N_\infty$ operad.
For $L\leq G$ and a finite $L$-set $T$, since
$E\mathcal F_{n-1}\simeq E(G/K)$ and
$\mathcal E(|T|)$ is $\Sigma_{|T|}$-free,
\[
\mathcal O(|T|)^{\Gamma_T}\simeq
\begin{cases}
*,&\sigma_T|_{L\cap K}=1,\\
\varnothing,&\sigma_T|_{L\cap K}\ne1.
\end{cases}
\]
Thus $N_K^G$ is the only nonidentity admissible norm.
The map from $\mathcal O$ to the commutative operad gives
$H$ an $\mathcal O$-algebra structure.

Write $S_{-1}=\mathrm{id}$ and
$S_s(-)=\widetilde E\mathcal F_s\wedge(-)$ for
$0\leq s\leq n-2$.
If $Z\in\Sp^K$ and
$\res_K^G\widetilde E\mathcal F_s\wedge Z\simeq0$, then
$\Phi^QZ\simeq0 \;(Q\leq K,\ Q\notin\mathcal F_s)$.
For $J\leq G$ with $J\notin\mathcal F_s$, we have
$J\cap K\notin\mathcal F_s$ and
\[
\Phi^J N_K^GZ
\simeq\bigwedge_{J\backslash G/K}\Phi^{J\cap K}Z
\simeq0
\]
on underlying spectra. Hence $S_sN_K^GZ\simeq0$.
By \cite[Theorem~5.2]{Carrick}, $S_s$ and its unit lift to
$\operatorname{Alg}_{\mathcal O}$.

For a family $\mathcal F$ and
$Y\in\operatorname{Alg}_{\mathcal O}$, the diagonal of
$E\mathcal F$ makes
$\eta_Y:Y\longrightarrow F((E\mathcal F)_+,Y)$
a map in $\operatorname{Alg}_{\mathcal O}$.
The diagonal and projections of
$E\mathcal F\times E\mathcal F$ give
\[
\begin{gathered}
F((E\mathcal F)_+,F((E\mathcal F)_+,Y))
\simeq F((E\mathcal F)_+,Y),\\
F((E\mathcal F)_+,\eta_Y)
\simeq\eta_{F((E\mathcal F)_+,Y)}.
\end{gathered}
\]
Hence $F((E\mathcal F)_+,-)$ is a localization of
$\operatorname{Alg}_{\mathcal O}$.
Consequently,
\[
\begin{aligned}
H^{\Phi_i}&=S_{i-1}H &&(0\leq i<n),\\
H^{h_i}&=F((E\mathcal F_i)_+,H^{\Phi_i}) &&(0\leq i<n),\\
H^{t_i}&=S_iH^{h_i} &&(0\leq i<n-1)
\end{aligned}
\]
are $\mathcal O$-algebras.
For $0\leq i<n-1$, the localization unit
$\ell_i:H^{h_i}\to H^{t_i}$ and the map
\[
\begin{aligned}
\theta_i:\quad H^{\Phi_{i+1}}
&\simeq S_iH^{\Phi_i}\\
&\xrightarrow{\ S_i\eta_{H^{\Phi_i}}\ }
S_iF((E\mathcal F_i)_+,H^{\Phi_i})=H^{t_i}
\end{aligned}
\]
are $\mathcal O$-algebra maps, where
$\eta_{H^{\Phi_i}}$ is the $\mathcal F_i$-completion unit.

By Corollary~\ref{cor:hz-module-descent}, for $0\leq i<n-1$,
\[
H^{t_i}\xrightarrow{\ \simeq\ }
F((E\mathcal F_{i+1})_+,H^{t_i}).
\]
Precomposition with $\mathcal F_{i+1}$-completion map gives
\[
\operatorname{Map}_{\operatorname{Alg}_{\mathcal O}}
(H^{h_{i+1}},H^{t_i})
\xrightarrow[\ \simeq\ ]{\ -\circ\eta_{H^{\Phi_{i+1}}}\ }
\operatorname{Map}_{\operatorname{Alg}_{\mathcal O}}
(H^{\Phi_{i+1}},H^{t_i}).
\]
Thus there is a contractibly unique factorization in
$\operatorname{Alg}_{\mathcal O}$:
\[
\begin{gathered}
H^{\Phi_{i+1}}
\xrightarrow{\ \eta_{H^{\Phi_{i+1}}}\ }
H^{h_{i+1}}\xrightarrow{\ \rho_i\ }H^{t_i},\\
\rho_i\eta_{H^{\Phi_{i+1}}}\simeq\theta_i.
\end{gathered}
\]
Naturality gives the cone from $H$.
The underlying diagram is $\mathcal R$ by
Definition~\ref{def:tate-compression-diagram}.

Finally,
$\mathcal F_i|_K=\{C_{p^j}:0\leq j\leq i\} \;(0\leq i<n-1)$, $\res_K^G E\mathcal F_{n-1}\simeq *$.

For $0\leq i<n-1$, there are natural equivalences
\[
\begin{aligned}
\res_K^G(\widetilde E\mathcal F_i\wedge X)
&\simeq
\widetilde E(\mathcal F_i|_K)\wedge\res_K^G X,\\
\res_K^G F_G((E\mathcal F_i)_+,X)
&\simeq
F_K((E(\mathcal F_i|_K))_+,\res_K^G X),
\end{aligned}
\]
compatible with the unit maps of localization and completion.
Hence, for $0\leq i<n-1$,
\[
\begin{aligned}
\res_K^G H^{h_i}
&\simeq(\res_K^G H)^{h_i},\\
\res_K^G H^{t_i}
&\simeq(\res_K^G H)^{t_i},\\
\res_K^G H^{h_{n-1}}
&\simeq(\res_K^G H)^{\Phi_{n-1}}.
\end{aligned}
\]

For $0\leq i<n-1$, restriction gives
\[
(\res_K^G\rho_i)\circ
(\res_K^G\eta_{H^{\Phi_{i+1}}})
\simeq\res_K^G\theta_i.
\]
By the contractible uniqueness of the factorization,
$\res_K^G\rho_i$ and its factorization homotopy agree
with those over $K$.
Naturality of the unit maps of localization and completion
gives the corresponding equivalence of cones.

Corollary~\ref{cor:hz-module-descent}, applied to $K$, gives
\[
\res_K^G H\xrightarrow{\ \simeq\ }
\holim_{\mathcal I_<}\res_K^G\mathcal R.
\]
\end{proof}

\subsection{Norm of Euler and orientation classes}

Let $K=C_{p^{n-1}}<G=C_{p^n}$, $n\geq2$.
For an actual $K$-representation $W$, applying the indexed smash product
to the Euler map $S^0\to S^W$ gives
\begin{equation}\label{eq:tambara-HHR-euler-orientation}
 N_K^G(a_W)=a_{\Ind_K^G W}.
\end{equation}
If $W$ is orientable and either $p$ is odd or $\dim W$ is even,
we also have
\begin{equation}\label{eq:tambara-orientation-identity}
 u_{\Ind_K^G W}
 =u_{\Ind_K^G(\dim W)}N_K^G(u_W)
 \qquad\text{\cite[Lemma~9.9.12]{HHRbook}}.
\end{equation}
The degree of both sides is
$[G:K]\dim W-\Ind_K^G W$.

Let
$E\in\{H^{h_i}:0\le i<n\}\cup\{H^{t_i}:0\le i<n-1\}$
be a $G$-spectrum in Lemma~\ref{lem:tambara-local-admissibility}. 
We have the following.

\begin{proposition}[Norms on $H^{h_i}$ and $H^{t_i}$]
\label{lem:tambara-cellular-generator}
We have
\begin{equation*}
 N_K^G(\prod_j a_{W_j}^{e_j}\prod_k u_{Q_k}^{f_k})=
 \prod_j a_{\Ind_K^G W_j}^{e_j}
 \prod_k
 \left(u_{\Ind_K^G Q_k}
       u_{\Ind_K^G(\dim Q_k)}^{-1}\right)^{f_k}.
\end{equation*}
\end{proposition}
\begin{proof}
Apply \eqref{eq:tambara-HHR-euler-orientation} and
\eqref{eq:tambara-orientation-identity} to each factor, using
multiplicativity and preservation of units.
\end{proof}

\begin{remark}
For $m\in\mathbb Z$, scalar multiplication gives
\[
 N_K^G(mg)=m^{[G:K]}N_K^G(g).
\]
The units $v_D$ of Lemma~\ref{lem:normalized-character-units}
transport the degrees to the selected degree.
\end{remark}

Since our computation using the valuation language, we have to translate Proposition~\ref{lem:tambara-cellular-generator} to valuation language. The following is primary but heavy-workloaded.

For the adjacent inclusion $K=C_{p^k}<G=C_{p^{k+1}}$, $k\ge1$,
and a natural degree $V\in RO(K)$, the norm maps are
\[
\begin{aligned}
 N_K^G:\pi_V^K(\res_K^G H^{h_i})
 &\longrightarrow\pi_{\Ind_K^G V}^G(H^{h_i})
 &&(0\le i\le k),\\
 N_K^G:\pi_V^K(\res_K^G H^{t_i})
 &\longrightarrow\pi_{\Ind_K^G V}^G(H^{t_i})
 &&(0\le i<k).
\end{aligned}
\]

\begin{definition}
\label{def:norm-valuation-notation}
For $L=C_{p^m}\in\{K,G\}$, use $D_m$ from
Definition~\ref{def:degree-coordinates} with $n=m$.
Write $\rho_L:D_m\to RO_{\mathrm{sel}}(L)$ for the map in
\eqref{eq:odd-forest-coordinate}--\eqref{eq:two-forest-coordinate}
formed over $L$.
For $\mathbf b\in D_m$ and
$E\in\{A_0,\ldots,A_m,T_0,\ldots,T_{m-1}\}$,
define, in each case listed in \eqref{eq:local-height-table},
\[
\begin{aligned}
\bigl(b_{E^L}(\mathbf b),e_{E^L}(\mathbf b)\bigr)
&:=\bigl(b_E(\mathbf b),e_E(\mathbf b)\bigr)\big|_{n=m},\\
c_{E^L}(\mathbf b)
&:=b_{E^L}(\mathbf b)+e_{E^L}(\mathbf b),\\
E^L(\mathbf b)
&:=\frac{p^{b_{E^L}(\mathbf b)}\mathbb Z}
        {p^{c_{E^L}(\mathbf b)}\mathbb Z}.
\end{aligned}
\]
Set $E^L(\mathbf b)=0$ in the remaining cases.
The products are those of Theorem~\ref{thm:local-forest-data}
for $L$.

Form the layers of $\res_L^G H$ over $L$ as in
\eqref{eq:layer-definitions}.
By Theorem~\ref{thm:local-forest-identification},
\[
E^L(\mathbf b)\cong
\begin{cases}
\pi_{\rho_L(\mathbf b)}^L\bigl((\res_L^G H)^{h_i}\bigr),
& E=A_i,\quad 0\leq i<m,\\
\pi_{\rho_L(\mathbf b)}^L\bigl((\res_L^G H)^{\Phi_m}\bigr),
& E=A_m,\\
\pi_{\rho_L(\mathbf b)}^L\bigl((\res_L^G H)^{t_i}\bigr),
& E=T_i,\quad 0\leq i<m.
\end{cases}
\]
For $E^L(\mathbf b)\ne0$, define $g_E^L(\mathbf b)$ by
$g_E^L(\mathbf b)\longleftrightarrow [p^{b_{E^L}(\mathbf b)}]$
under these isomorphisms.
We also use $[q]$ for the corresponding homotopy class when
$q\in p^{b_{E^L}(\mathbf b)}\mathbb Z$.
\end{definition}

For $V=\rho_K(\mathbf a)$ with $\mathbf a\in D_k$,
Lemma~\ref{lem:tambara-local-admissibility} gives the source
homotopy groups of the preceding norm maps as
$A_i^K(\mathbf a)$ and $T_i^K(\mathbf a)$.
In particular,
\[
\begin{aligned}
A_k^K(\mathbf a)
&\cong\pi_{\rho_K(\mathbf a)}^K
       \bigl((\res_K^G H)^{\Phi_k}\bigr)\\
&\cong\pi_{\rho_K(\mathbf a)}^K(\res_K^G H^{h_k}).
\end{aligned}
\]
The following is to translate Theorem~\ref{lem:tambara-cellular-generator} to valuation language.
\begin{theorem}
\label{thm:exact-weighted-local-norm}
Use the notation of Definition~\ref{def:norm-valuation-notation}.
Let $\mathbf a=(a_0,\ldots,a_k)\in D_k$, and let
$E=A_i$ $(0\leq i\leq k)$ or $E=T_i$ $(0\leq i<k)$
with $E^K(\mathbf a)\ne0$.
For $q\in p^{b_{E^K}(\mathbf a)}\mathbb Z$,
using the $v_D$ of Lemma~\ref{lem:normalized-character-units},
we have
\begin{equation}\label{eq:exact-weighted-cyclic-norm}
 v_{\Ind_K^G\rho_K(\mathbf a)-\rho_G(I\mathbf a)}^{-1}N_K^G([q])
 =[\eta_Ep^{a_k}q^p]\in E^G(I\mathbf a).
\end{equation}
Here $I\mathbf a=(pa_0,\ldots,pa_k,a_k)$ as in
\eqref{eq:tambara-adjacent-degree}. $\eta_E$ is $1$ for
odd $p$, and for $p=2$ is given by
\eqref{eq:norm-dyadic-edge-units}--\eqref{eq:norm-carry-and-unit}.
Moreover, let $N_K^G(0)=0$.
\end{theorem}

\begin{proof}[Proof of Theorem~\ref{thm:exact-weighted-local-norm}]
\noindent1. Verification for Euler and orientation classes.

For $p=2$, use $t_i(\mathbf a)$ from
Section~\ref{sec:local} and, for $0\leq i<k$, put
\begin{equation}\label{eq:norm-dyadic-edge-units}
 \eta_{T_i}(\mathbf a)=1+2^{k-i}t_i(\mathbf a).
\end{equation}
\begin{equation}\label{eq:norm-dyadic-vertex-units}
\begin{gathered}
 \eta_{A_0}(\mathbf a)=\eta_{T_0}(\mathbf a),\qquad
 \eta_{A_k}(\mathbf a)=\eta_{T_{k-1}}(\mathbf a),\\
 \eta_{A_i}(\mathbf a)=\begin{cases}
  \eta_{T_i}(\mathbf a),&t_i(\mathbf a)<0,\\
  \eta_{T_{i-1}}(\mathbf a),&t_i(\mathbf a)\geq0
 \end{cases}
 \qquad(1\leq i<k).
\end{gathered}
\end{equation}
For every $p$, use
\begin{equation}\label{eq:norm-carry-and-unit}
\begin{gathered}
 \eta_E(\mathbf a)\in
 \bigl(\mathbb Z/p^{e_{E^G}(I\mathbf a)}\bigr)^\times
 \quad(e_{E^G}(I\mathbf a)<\infty),\\
 \eta_E(\mathbf a)=1
 \quad(p\text{ odd or }e_{E^G}(I\mathbf a)=\infty).
\end{gathered}
\end{equation}
For $E=A_0$ with $a_0=0$, take the integer lift
$\eta_{A_0}(\mathbf a)=1$.

Suppose first that $\epsilon(\mathbf a)=0$.
For $0\leq j<k$, the complex-linear decomposition
\[
 \Ind_K^G\lambda_j
 =\bigoplus_{r=0}^{p-1}\lambda(p^j+rp^k)
\]
preserves orientations. For each of the two inputs below,
$\rho_K(\mathbf a)$ denotes its degree. Equations
\eqref{eq:tambara-HHR-euler-orientation},
\eqref{eq:tambara-orientation-identity}, and
\eqref{eq:character-euler-factor} give
\[
 v_{\Ind_K^G\rho_K(\mathbf a)-\rho_G(I\mathbf a)}^{-1}
 N_K^G([q])
 =\begin{cases}
 \displaystyle
 \left(\prod_{r=0}^{p-1}(1+rp^{k-j})\right)a_{\lambda_j}^p,
   &[q]=a_{\lambda_j},\\[2mm]
 u_{\lambda_j}^pu_{\lambda_k}^{-(p-1)},
   &[q]=u_{\lambda_j}.
 \end{cases}
\]
Here $u_{\lambda_k}$ is a unit in the target ring.
For $p=2$, $\lambda_{k-1}=2\sigma_{2^k}$ in the source and
$\lambda_k=2\sigma_{2^{k+1}}$ in the target.
In valuation language,
\[
\begin{aligned}
 a_{\lambda_j}&\longleftrightarrow[1],&
 u_{\lambda_j}&\longleftrightarrow[p^{k-j}]
 &&\text{in the source},\\
 a_{\lambda_j}&\longleftrightarrow[1],&
 u_{\lambda_j}&\longleftrightarrow[p^{k+1-j}]
 &&\text{in the target}.
\end{aligned}
\]
The target formula also holds for $j=k$.

For $a_{\lambda_j}$, the coordinates are
$a_\ell=-\mathbf1_{\ell\leq j}$ $(0\leq\ell\leq k)$,
so $q=1$ and $a_k=0$.
For odd $p$,
\[
 \prod_{r=0}^{p-1}(1+rp^{k-j})
 \equiv1\pmod{p^{k-j+1}},\qquad
 c_{E^G}(I\mathbf a)\leq k-j+1.
\]
For $p=2$, the product is $1+2^{k-j}$, and substitution in
\eqref{eq:local-height-table} and
\eqref{eq:norm-dyadic-edge-units}--\eqref{eq:norm-carry-and-unit}
gives
\[
 1+2^{k-j}-\eta_E(\mathbf a)
 \in2^{c_{E^G}(I\mathbf a)}\mathbb Z.
\]
Thus the formula holds for $a_{\lambda_j}$.

For $u_{\lambda_j}$, the coordinates are
$a_\ell=\mathbf1_{\ell>j}$ $(0\leq\ell\leq k)$,
so $q=p^{k-j}$ and $a_k=1$. We have
\[
 u_{\lambda_j}^pu_{\lambda_k}^{-(p-1)}
 \longleftrightarrow
 [p^{p(k+1-j)-(p-1)}]
 =[p^{p(k-j)+1}].
\]
For odd $p$, $\eta_E(\mathbf a)=1$.
For $p=2$ and finite $e_{E^G}(I\mathbf a)$, substitution in
the same formulas gives
\[
 (\eta_E(\mathbf a)-1)2^{2(k-j)+1}
 \in2^{c_{E^G}(I\mathbf a)}\mathbb Z.
\]
For $E=A_0$, $a_0=0$ and $\eta_{A_0}(\mathbf a)=1$.
Hence
\[
 [p^{p(k-j)+1}]
 =[\eta_E(\mathbf a)p^{a_k}q^p],
\]
which proves the formula for $u_{\lambda_j}$.

\par\medskip
\noindent2. Preservation under products.

Let $[q]\in E^K(\mathbf a)$ and $[r]\in E^K(\mathbf d)$
be monomials in these Euler and orientation classes, with
$\mathbf a,\mathbf d\in D_k^{(0)}$, and suppose that both satisfy
\eqref{eq:exact-weighted-cyclic-norm}.
Both $\rho_K(\mathbf a)$ and $\rho_K(\mathbf d)$ have even
fixed-point dimensions at every subgroup of $K$, so no
permutation sign occurs.
Since $I$, $\rho_K$, and $\rho_G$ are additive and
\eqref{eq:character-unit-coherence} holds, it remains to verify
in $E^G(I(\mathbf a+\mathbf d))$ that
\[
\begin{aligned}
 &[\eta_E(\mathbf a)\eta_E(\mathbf d)
       p^{a_k+d_k}(qr)^p]\\
 &\qquad=[\eta_E(\mathbf a+\mathbf d)
       p^{a_k+d_k}(qr)^p].
\end{aligned}
\]
By \eqref{eq:local-height-table},
\[
 p b_{E^K}(\mathbf a)+a_k=b_{E^G}(I\mathbf a),\qquad
 p^{a_k+d_k}(qr)^p
 \in p^{b_{E^G}(I(\mathbf a+\mathbf d))}\mathbb Z.
\]
For odd $p$,
$\eta_E(\mathbf a)=\eta_E(\mathbf d) =\eta_E(\mathbf a+\mathbf d)=1$.
For $p=2$, we have, for $0\leq j<k$,
\[
 (1+2^{k-j}a_j)(1+2^{k-j}d_j)
 \equiv1+2^{k-j}(a_j+d_j)
 \pmod{2^{k-j+1}}.
\]
This proves the required equality for $E=T_i$, $A_0$, and $A_k$,
and for $E=A_i$ when $a_i,d_i<0$ or $a_i,d_i\geq0$.
In the torsion-free case $E=A_0$ with $a_0+d_0=0$,
we have $a_0=d_0=0$, hence
$\eta_{A_0}(\mathbf a)=\eta_{A_0}(\mathbf d) =\eta_{A_0}(\mathbf a+\mathbf d)=1$.

It remains to consider $E=A_i$, $1\leq i<k$, with
$a_i<0\leq d_i$, after interchanging the factors if necessary.
If $d_i=0$, then
\[
 \eta_{A_i}(\mathbf a+\mathbf d)=\eta_{A_i}(\mathbf a),\qquad
 \eta_{A_i}(\mathbf d)\equiv1\pmod{2^{k-i+1}},\qquad
 e_{A_i^G}(I(\mathbf a+\mathbf d))=k-i+1.
\]
If $d_i>0$, then
\[
\begin{aligned}
 &b_{A_i^K}(\mathbf a)+b_{A_i^K}(\mathbf d)
       -b_{A_i^K}(\mathbf a+\mathbf d)\\
 &\qquad=a_i^++d_i^+-(a_i+d_i)^+
       =\min\{-a_i,d_i\}\geq1.
\end{aligned}
\]
Consequently,
\[
\begin{gathered}
 2^{a_k+d_k}(qr)^2
 \in2^{b_{A_i^G}(I(\mathbf a+\mathbf d))+2}\mathbb Z,\\
 \eta_{A_i}(\mathbf a)\eta_{A_i}(\mathbf d)
       -\eta_{A_i}(\mathbf a+\mathbf d)
 \in2^{k-i}\mathbb Z,\\
 e_{A_i^G}(I(\mathbf a+\mathbf d))\leq k-i+2.
\end{gathered}
\]
Thus the required equality holds in both cases.
By multiplicativity,
\[
\begin{aligned}
 &v_{\Ind_K^G\rho_K(\mathbf a+\mathbf d)
       -\rho_G(I(\mathbf a+\mathbf d))}^{-1}N_K^G([qr])\\
 &\qquad=[\eta_E(\mathbf a)\eta_E(\mathbf d)
                 p^{a_k+d_k}(qr)^p]\\
 &\qquad=[\eta_E(\mathbf a+\mathbf d)
                 p^{a_k+d_k}(qr)^p].
\end{aligned}
\]
Hence the product also satisfies
\eqref{eq:exact-weighted-cyclic-norm}.

\par\medskip
\noindent3. Preservation under inverses.

Suppose that the monomial $[q]$ is a unit and satisfies
\eqref{eq:exact-weighted-cyclic-norm}.
Choose a representative $q\in p^{\mathbb Z}$, so that
$[q]^{-1}=[q^{-1}]\in E^K(-\mathbf a)$.
For the inverse, it remains to verify in $E^G(0)$ that
\[
 [\eta_E(\mathbf a)p^{a_k}q^p]
 [\eta_E(-\mathbf a)p^{-a_k}q^{-p}]
 =[1].
\]
Taking $\mathbf d=-\mathbf a$ and $r=q^{-1}$ in the product
computation gives
\[
 [\eta_E(\mathbf a)\eta_E(-\mathbf a)]
 =[\eta_E(0)]=[1].
\]
By preservation of units and \eqref{eq:character-unit-coherence},
\[
\begin{aligned}
 &v_{\Ind_K^G\rho_K(-\mathbf a)-\rho_G(I(-\mathbf a))}^{-1}
     N_K^G([q^{-1}])\\
 &\qquad=
 \left(v_{\Ind_K^G\rho_K(\mathbf a)-\rho_G(I\mathbf a)}^{-1}
     N_K^G([q])\right)^{-1}\\
 &\qquad=[\eta_E(-\mathbf a)p^{-a_k}q^{-p}].
\end{aligned}
\]
Thus the inverse also satisfies
\eqref{eq:exact-weighted-cyclic-norm}.

By Propositions~\ref{prop:odd-fractional-lattices}
and~\ref{prop:two-fractional-lattices},
$g_E^K(\mathbf a)$ is a monomial in Euler and orientation
classes, with negative powers only of units. Steps 1--3 therefore give
\[
 v_{\Ind_K^G\rho_K(\mathbf a)-\rho_G(I\mathbf a)}^{-1}
 N_K^G(g_E^K(\mathbf a))
 \longleftrightarrow
 [\eta_E(\mathbf a)p^{p b_{E^K}(\mathbf a)+a_k}]
 \in E^G(I\mathbf a).
\]

Now let $p=2$ and $\epsilon(\mathbf a)=1$.
Proposition~\ref{prop:two-fractional-lattices},
\eqref{eq:local-height-table}, and
\eqref{eq:norm-dyadic-edge-units}--\eqref{eq:norm-carry-and-unit} give
\[
\begin{aligned}
 g_E^K(\mathbf a)
 &=a_{\sigma_{2^k}}
   g_E^K\bigl(\mathbf a+(\tfrac12,\ldots,\tfrac12,0)\bigr),\\
 b_{E^K}\bigl(\mathbf a+(\tfrac12,\ldots,\tfrac12,0)\bigr)
 &=b_{E^K}(\mathbf a),\\
 \eta_E\bigl(\mathbf a+(\tfrac12,\ldots,\tfrac12,0)\bigr)
 &=\eta_E(\mathbf a).
\end{aligned}
\]
Moreover,
\[
\begin{aligned}
 \rho_K\bigl(\mathbf a+(\tfrac12,\ldots,\tfrac12,0)\bigr)
 &=\rho_K(\mathbf a)+\sigma_{2^k},\\
 \rho_G\bigl(I(\mathbf a+(\tfrac12,\ldots,\tfrac12,0))\bigr)
 &=\rho_G(I\mathbf a)+\lambda_{k-1},\\
 \Ind_K^G\sigma_{2^k}&=\lambda_{k-1},\\
 N_K^G(a_{\sigma_{2^k}})
 &=a_{\lambda_{k-1}}\longleftrightarrow[1].
\end{aligned}
\]
Thus
\[
\begin{aligned}
 &\Ind_K^G\rho_K\bigl(\mathbf a+(\tfrac12,\ldots,\tfrac12,0)\bigr)
 -\rho_G\bigl(I(\mathbf a+(\tfrac12,\ldots,\tfrac12,0))\bigr)\\
 &\qquad=\Ind_K^G\rho_K(\mathbf a)-\rho_G(I\mathbf a).
\end{aligned}
\]
Since $\dim\rho_K\bigl(\mathbf a+(\tfrac12,\ldots,\tfrac12,0)\bigr)$
is even, multiplication by $a_{\sigma_{2^k}}$ introduces no
permutation sign. Multiplication by $a_{\lambda_{k-1}}$ is
$[q]\mapsto[q]$ in the target valuation language.
The generator formula therefore holds also when
$\epsilon(\mathbf a)=1$.

Since $q/p^{b_{E^K}(\mathbf a)}\in\mathbb Z$,
Proposition~\ref{lem:tambara-cellular-generator} gives
\[
\begin{aligned}
 &v_{\Ind_K^G\rho_K(\mathbf a)-\rho_G(I\mathbf a)}^{-1}
 N_K^G([q])\\
 &\quad=\left[
 \eta_E(\mathbf a)p^{p b_{E^K}(\mathbf a)+a_k}
 \left(\frac{q}{p^{b_{E^K}(\mathbf a)}}\right)^p
 \right]
 =[\eta_E(\mathbf a)p^{a_k}q^p].
\end{aligned}
\]
For $p=2$ and $\epsilon(\mathbf a)=1$,
substitution in \eqref{eq:local-height-table} gives
\[
 p b_{E^K}(\mathbf a)+a_k-b_{E^G}(I\mathbf a)
 =\begin{cases}
 k-1,&E=A_0,\\
 k-i-1+\mathbf1_{t_i(\mathbf a)\geq1},&E=A_i,\quad 1\leq i<k,\\
 k-i-1,&E=T_i,\quad 0\leq i<k,\\
 0,&E=A_k.
 \end{cases}
\]
These expressions are nonnegative, so
$[\eta_E(\mathbf a)p^{a_k}q^p]$ defines a class in $E^G(I\mathbf a)$.
For $p=2$ and $\epsilon(\mathbf a)=1$, we also have
\[
 e_{E^G}(I\mathbf a)
 -\bigl(p b_{E^K}(\mathbf a)+a_k-b_{E^G}(I\mathbf a)\bigr)=2,
\]
so only $\eta_E\bmod4$ affects the target class.
For $e_{E^G}(I\mathbf a)<\infty$, changing the integer lift
of $\eta_E$ does not change this class, since
\[
\begin{aligned}
 p^{e_{E^G}(I\mathbf a)+a_k+p b_{E^K}(\mathbf a)}\mathbb Z
 &=p^{c_{E^G}(I\mathbf a)
 +p b_{E^K}(\mathbf a)+a_k-b_{E^G}(I\mathbf a)}\mathbb Z\\
 &\subseteq p^{c_{E^G}(I\mathbf a)}\mathbb Z.
\end{aligned}
\]
The formula is independent of the representative $q$, since
the left-hand side is the norm of $[q]$.
\end{proof}
\subsection{The norm on vertex elements}

For $K=C_{p^{n-1}}<G=C_{p^n}$, $n\ge2$, and a natural degree
$V\in RO(K)$, the norm on vertex elements is
\[
 N_K^G:\pi_V^K(H)\longrightarrow\pi_{\Ind_K^G V}^G(H).
\]
For $0\le i<n$, its $i$th vertex entry is computed by the norm map
\[
 N_K^G:\pi_V^K(\res_K^G H^{h_i})
 \longrightarrow\pi_{\Ind_K^G V}^G(H^{h_i}).
\]
Theorem~\ref{thm:local-forest-identification} gives the homotopy groups
of $H^{h_i}$ and $H^{t_i}$ as $A_i$ and $T_i$, respectively.
We use their valuation language to describe the map.
Equation~\eqref{eq:exact-weighted-cyclic-norm} gives the first $n$
entries of the resulting vertex array $(b_0,\ldots,b_n)$.
The equation
$\rho_{n-1}(b_n)=\ell_{n-1}(b_{n-1})$
determines the last entry.

\begin{theorem}[Norm on vertex elements]
\label{thm:tambara-natural-kernel}
Let $V=\rho_K(\mathbf a)$ be a natural degree and let
$x=(x_0,\ldots,x_{n-1})\in\pi_V^K(H)\cong\ker\Delta^K_{\mathbf a}$.
For each nonzero $A_i^K(\mathbf a)$, choose
$q_i\in p^{b_{A_i^K}(\mathbf a)}\mathbb Z$ with $x_i=[q_i]$;
put $q_i=0$ when $A_i^K(\mathbf a)=0$. Then
\[
 v_{\Ind_K^G V-\rho_G(I\mathbf a)}^{-1}N_K^G(x)
 =(b_0,\ldots,b_n)
 \in\pi_{\rho_G(I\mathbf a)}^G(H)\cong\ker\Delta^G_{I\mathbf a},
\]
where
\[
 b_i=
 \begin{cases}
 [\eta_{A_i}(\mathbf a)p^{a_{n-1}}q_i^p]\in A_i^G(I\mathbf a),
   &0\le i<n,\ A_i^K(\mathbf a)\ne0,\\
 [q_{n-1}^p]\in\mathbb F_p,&i=n,\ a_{n-1}=0,\\
 0,&\text{otherwise}.
 \end{cases}
\]
Here $\eta_{A_i}$ is given by
\eqref{eq:norm-dyadic-edge-units}--\eqref{eq:norm-carry-and-unit},
with $k=n-1$.
\end{theorem}

\begin{proof}
By Lemma~\ref{lem:tambara-local-admissibility}, the maps $H\to H^{h_i}$
commute with the $G/K$-norm. ~\eqref{eq:exact-weighted-cyclic-norm}
therefore gives the first $n$ entries $b_0,\ldots,b_{n-1}$.
$\rho_{n-1}:A_n^G(I\mathbf a)\to T_{n-1}^G(I\mathbf a)$ is
\[
 \rho_{n-1}=\begin{cases}
  0\longrightarrow p^{a_{n-1}}\mathbb Z/p^{a_{n-1}+1}\mathbb Z,
       &a_{n-1}<0,\\
  \operatorname{id}_{p^{a_{n-1}}\mathbb Z/p^{a_{n-1}+1}\mathbb Z},
       &a_{n-1}\ge0.
 \end{cases}
\]
Hence $\rho_{n-1}(b_n)=\ell_{n-1}(b_{n-1})$ determines $b_n$ uniquely.
For $a_{n-1}=0$, the right side is $[q_{n-1}^p]$, since
$\eta_{A_{n-1}}\equiv1\pmod p$. For $a_{n-1}>0$,
we have $q_{n-1}\in p^{a_{n-1}}\mathbb Z$, so
\[
 p^{a_{n-1}}q_{n-1}^p\in p^{(p+1)a_{n-1}}\mathbb Z
 \subseteq p^{a_{n-1}+1}\mathbb Z,
\]
and $b_n=0$. For $a_{n-1}<0$, $A_n^G(I\mathbf a)=0$.
Thus we get the proof.
\end{proof}

\subsection{The norm on $F^1$}
In this subsection we get the norm on edge elements.

Since the degree of edge elements is shifted, we first introduce the following notations.
Let $K=C_{p^{n-1}}<G=C_{p^n}$, $n\ge2$, let $V$ be a shifted degree,
and write $V+1=\rho_K(\mathbf a)$. Define
\begin{equation}\label{eq:tambara-shifted-companion}
 J(\mathbf a)_j=
 \begin{cases}
 pa_j-(p-1)/2,&0\le j<n,\\
 a_{n-1},&j=n,
 \end{cases}
\end{equation}
and
\begin{equation}\label{eq:tambara-J-degree}
 D_V=\Ind_K^GV+1-\rho_G(J(\mathbf a)).
\end{equation}
Proposition~\ref{prop:tambara-induced-degree} gives
$\dim(D_V^L)=0$ for every $L\le G$.  Write
\[
 a_K^G=a_{\mathbb R[G/K]-1}
 \in\pi_{1-\mathbb R[G/K]}^G(\mathbb S).
\]
Use the same notation for its images under the unit maps to
$H$, $H^{h_i}$, and $H^{t_i}$.
For each $E=H^{h_i}$ $(0\le i<n)$ or $H^{t_i}$ $(0\le i<n-1)$,
write
\[
 N_K^G:\pi_{V+1}^K(\res_K^GE)
 \longrightarrow\pi_{\Ind_K^G(V+1)}^G(E).
\]

To apply Lemma~\ref{lem:tambara-local-admissibility} to edge elements,
we first compare the cokernels of the incidence maps for
$\mathcal I$ and $\mathcal I_<$.

Let $K=C_{p^{n-1}}<G=C_{p^n}$, $n\ge2$, and let
$V+1=\rho_K(\mathbf a)$ with $V$ a shifted degree. Use $J(\mathbf a)$
from \eqref{eq:tambara-shifted-companion}. The truncated incidence map for
$\mathcal I_<$ is
\[
\begin{aligned}
 \Delta^{G,<}_{J\mathbf a}:\prod_{i=0}^{n-1}A_i^G(J\mathbf a)
 &\longrightarrow\prod_{i=0}^{n-2}T_i^G(J\mathbf a),\\
 (x_i)_{i<n}&\longmapsto
 (\ell_i(x_i)-\rho_i(x_{i+1}))_{i<n-1}.
\end{aligned}
\]
Let
\[
\begin{aligned}
 \operatorname{pr}^0:\prod_{i=0}^{n}A_i^G(J\mathbf a)
 &\longrightarrow\prod_{i=0}^{n-1}A_i^G(J\mathbf a),\\
 \operatorname{pr}^1:\prod_{i=0}^{n-1}T_i^G(J\mathbf a)
 &\longrightarrow\prod_{i=0}^{n-2}T_i^G(J\mathbf a)
\end{aligned}
\]
drop the last entry.

\begin{lemma}
\label{lem:tambara-terminal-elimination}
Then
$\operatorname{pr}^1\Delta^G_{J\mathbf a}
 =\Delta^{G,<}_{J\mathbf a}\operatorname{pr}^0$,
and $\operatorname{pr}^1$ induces an isomorphism
\begin{equation*}
\coker\Delta^G_{J\mathbf a}
 \xrightarrow{\ \cong\ }
 \coker\Delta^{G,<}_{J\mathbf a}.
\end{equation*}
Explicitly,
\[
\begin{aligned}
 \overline{\operatorname{pr}}^1([z_0,\ldots,z_{n-1}])
     &=[z_0,\ldots,z_{n-2}],\\
 (\overline{\operatorname{pr}}^1)^{-1}([z_0,\ldots,z_{n-2}])
     &=[z_0,\ldots,z_{n-2},0].
\end{aligned}
\]
\end{lemma}
\begin{proof}
The maps $\operatorname{pr}^0$ and $\operatorname{pr}^1$ are surjective, and
$\operatorname{pr}^1(\operatorname{im}\Delta^G_{J\mathbf a}) =\operatorname{im}\Delta^{G,<}_{J\mathbf a}$.
For $t\in T_{n-1}^G(J\mathbf a)$, we show
$(0,\ldots,0,t)\in\operatorname{im}\Delta^G_{J\mathbf a}$.
If $a_{n-1}\ge0$,  $\rho_{n-1}$ is the identity on
$p^{a_{n-1}}\mathbb Z/p^{a_{n-1}+1}\mathbb Z$, so
\[
 \Delta^G_{J\mathbf a}(0,\ldots,0,-t)=(0,\ldots,0,t).
\]

If $a_{n-1}<0$, the map
$\ell_{n-1}:A_{n-1}^G(J\mathbf a)\to T_{n-1}^G(J\mathbf a)$
is an isomorphism. By \eqref{eq:local-comparison-displacements},
\[
 b_{A_{n-1}}(J\mathbf a)-b_{T_{n-2}}(J\mathbf a)
 =\begin{cases}
  -pa_{n-1}+(p-1)/2\ge2,&p\text{ odd},\\
  -2a_{n-1}\ge2,&p=2.
 \end{cases}
\]
$T_{n-2}^G(J\mathbf a)$ has length $2$ for odd $p$ and
length $1$ for $p=2$. Thus $\rho_{n-2}=0$ and
\[
 \Delta^G_{J\mathbf a}(0,\ldots,0,\ell_{n-1}^{-1}(t),0)
 =(0,\ldots,0,t).
\]
Hence $\ker\operatorname{pr}^1\subseteq\operatorname{im}\Delta^G_{J\mathbf a}$.
\end{proof}

The following identity describes the map on $\pi_1$ induced by the
norm, with the mapping spaces based at the zero maps.

\begin{lemma}
\label{lem:norm-loop-orientation}
Let $K<G$ have prime index $p$, let $R$ be a genuine $G$--spectrum
with $G/K$--indexed multiplication, and let $V\in RO(K)$.
Base both mapping spaces at the zero map. Applying $N_K^G$ to a
sphere map and then the indexed multiplication induces the operation
\[
 \pi_1\operatorname{Map}_{\mathrm{Sp}^K}(S^V,\res_K^G R)
 \longrightarrow
 \pi_1\operatorname{Map}_{\mathrm{Sp}^G}(S^{\Ind_K^G V},R)
\]
given, for $x\in\pi_{V+1}^K(\res_K^G R)$, by
\[
 x\longmapsto a_{\mathbb R[G/K]-1}\,N_K^G(x)
\]
in degree $\Ind_K^G V+1$.
\end{lemma}
\begin{proof}
By \cite[Theorem~3.1.1 and Section~3.2]{Balderrama},
\[
\begin{gathered}
 S^1
 \xrightarrow{\operatorname{id}_{S^1}\wedge a_{\mathbb R[G/K]-1}}
 S^1\wedge S^{\mathbb R[G/K]-1}
 \cong N_K^G(S^1),\\[4pt]
 \begin{aligned}
 \pi_1\bigl(
   \operatorname{Map}_{\mathrm{Sp}^K}(S^V,\res_K^G R),0
 \bigr)
 &\longrightarrow
 \pi_1\bigl(
   \operatorname{Map}_{\mathrm{Sp}^G}(S^{\Ind_K^G V},R),0
 \bigr),\\
 x&\longmapsto a_{\mathbb R[G/K]-1}N_K^G(x).
 \end{aligned}
\end{gathered}
\]

We next check the signs of the sphere isomorphisms.
Using
$S^V=S^{V+1}\wedge(S^1)^\vee$, the symmetry
\[
 \bigwedge_{G/K}\bigl(S^{V+1}\wedge(S^1)^\vee\bigr)
 \cong
 \left(\bigwedge_{G/K}S^{V+1}\right)
 \wedge
 \left(\bigwedge_{G/K}(S^1)^\vee\right)
\]
has underlying degree $(-1)^{\binom p2(\dim V+1)}$.
The duality isomorphism
\[
 \bigwedge_{G/K}(S^1)^\vee
 \xrightarrow{\cong}
 \left(\bigwedge_{G/K}S^1\right)^\vee
\]
has underlying degree $(-1)^{\binom p2}$, since duality reverses
the order of the $p$ odd sphere factors.
This change in the order of factors
\[
 \left(\bigwedge_{G/K}S^V\right)
 \wedge\left(\bigwedge_{G/K}S^1\right)
 \cong
 \bigwedge_{G/K}(S^V\wedge S^1)
\]
has underlying degree $(-1)^{\binom p2\dim V}$.
The product of these signs is
\[
 (-1)^{\binom p2(\dim V+1)}
 (-1)^{\binom p2}
 (-1)^{\binom p2\dim V}=1.\qedhere
\]
\end{proof}

We use the pointed-space functoriality of homotopy-limit spectral
sequences from \cite[Theorem~2.5.3]{Balderrama} to compute the induced
map on the edge cokernel.

\begin{lemma}
\label{lem:tambara-pointed-path}
Let $\mathcal I$ be the path defined in
Section~\ref{sec:incidence-category}.
Let $X,Y:\mathcal I\to\mathrm{Sp}$
be diagrams such that $\pi_0X(j)=\pi_0Y(j)=0$ for every object
$j\in\mathcal I$. Write $\ell_i,\rho_i$ for the corresponding maps for both $X$ and $Y$, so that the corresponding incidence maps $\Delta_X,\Delta_Y$ is $(x_i)\mapsto(\ell_i x_i-\rho_i x_{i+1})$.
Let
$f:\Omega^\infty X\longrightarrow\Omega^\infty Y$
be a natural transformation of diagrams of pointed spaces. Then
the map on $\pi_0\operatorname{holim}\Omega^\infty X$ is,
under the cokernel isomorphisms,
\[
\begin{gathered}
 \coker\!\left(
   \prod_i\pi_1X(v_i)\xrightarrow{\Delta_X}
   \prod_i\pi_1X(e_i)\right)
 \longrightarrow
 \coker\!\left(
   \prod_i\pi_1Y(v_i)\xrightarrow{\Delta_Y}
   \prod_i\pi_1Y(e_i)\right),\\
 [z_i]\longmapsto[\pi_1(f_{e_i})(z_i)].
\end{gathered}
\]
This map is a group homomorphism.
\end{lemma}
\begin{proof}
The fiber sequence
\[
\operatorname{holim}X\longrightarrow
\prod_i X(v_i)\xrightarrow{\Delta_X}\prod_i X(e_i)
\]
and $\pi_0X(j)=0$ give
$\pi_0\operatorname{holim}X\cong\coker\Delta_X$,
and similarly for $Y$.
A class $[z_i]$ is represented by zero maps at the vertices and
based loops $z_i$ at the edges.
Since $f$ is pointed, these loops map to
$\pi_1(f_{e_i})(z_i)$.
Moreover,
\[
\Delta_Y\circ\prod_i\pi_1(f_{v_i})
=
\left(\prod_i\pi_1(f_{e_i})\right)\circ\Delta_X.
\]
Thus the induced cokernel map is well-defined and additive.
\end{proof}

\begin{remark}
The descent filtration $F^\bullet$
from totalization on $\pi_0\operatorname{holim}X$ and
$\pi_0\operatorname{holim}Y$ has $F^2=0$.

There is no composable pair of nonidentity arrows in $\mathcal I$.
The normalized cosimplicial replacements of $X$ and $Y$ therefore
vanish above degree one, giving $F^2=0$.

\end{remark}
Since $F^2=0$, there is no higher-filtration ambiguity.
We can now compute the norm on edge elements.

\begin{theorem}[Norms of edge elements]
\label{thm:tambara-shifted-cokernel}
Let $K=C_{p^{n-1}}<G=C_{p^n}$, $n\ge2$, and let
$V+1=\rho_K(\mathbf a)$ with $V$ a shifted degree. Use $J(\mathbf a)$
and $D_V$ from \eqref{eq:tambara-shifted-companion}--\eqref{eq:tambara-J-degree}.
For $[z_i]\in\pi_V^K(H)=\coker\Delta^K_{\mathbf a}$,
with $z_i\in T_i^K(\mathbf a)$, the map $v_{D_V}^{-1}N_K^G$ is
\begin{equation}\label{eq:tambara-euler-power}
 [z_0,\ldots,z_{n-2}]
 \longmapsto
 \left[v_{D_V}^{-1}a_K^GN_K^G(z_0),\ldots,
       v_{D_V}^{-1}a_K^GN_K^G(z_{n-2}),0\right]
 \in\coker\Delta^G_{J(\mathbf a)}.
\end{equation}
The cokernel map is well-defined and additive.
\end{theorem}
\begin{proof}
Use the diagram $\mathcal R:\mathcal I_<\to\operatorname{Alg}_{\mathcal O}$
with $\mathcal R(v_i)=H^{h_i}$, $\mathcal R(e_i)=H^{t_i}$ and
arrows $\ell_i,\rho_i$ from Lemma~\ref{lem:tambara-local-admissibility}. Write
$R_<=\operatorname{holim}_{\operatorname{Alg}_{\mathcal O}} \mathcal R$.
Limits of operadic algebras are created on underlying spectra.
The compatible maps $H\to\mathcal R(j)$ induce an
$\mathcal O$-algebra map $q:H\to R_<$.
We applied Lemma~\ref{lem:tambara-local-admissibility}
and Theorem~\ref{thm:sharp-path-compression} to $K$, thus $\res_K^Gq$ is an equivalence.

For an object $j\in\mathcal I_<$ define the ordinary mapping
spectra
\[
 X(j)=\operatorname{map}_{\mathrm{Sp}^K}
           (S^V,\res_K^G\mathcal R(j)),\qquad
 Y(j)=\operatorname{map}_{\mathrm{Sp}^G}
           (S^{\Ind_K^GV},\mathcal R(j)).
\]
Choose one point of $\mathcal O(p)^{\Gamma_{G/K}}$ and use its
indexed multiplication $\mu_j$ at every object $j$. This gives a pointed map
\begin{equation*}
 \begin{split}
 f_j:\Omega^\infty X(j)&\longrightarrow\Omega^\infty Y(j),\\
 g&\longmapsto
 \left(S^{\Ind_K^GV}\cong N_K^GS^V
 \xrightarrow{N_K^Gg}N_K^G\res_K^G\mathcal R(j)
 \xrightarrow{\mu_j}\mathcal R(j)\right).
 \end{split}
\end{equation*}
Since $G/K\ne\varnothing$, $f_j(0)=0$.
By Lemma~\ref{lem:tambara-local-admissibility}, the maps $\ell_i$ and
$\rho_i$ are $\mathcal O$-algebra maps. Hence $(f_j)$ is a natural transformation
of diagrams of pointed spaces.

The operation on $\pi_1$ is
\[
 \pi_1f_j=a_K^GN_K^G:
 \pi_{V+1}^K\res_K^G\mathcal R(j)
 \longrightarrow
 \pi_{\Ind_K^GV+1}^G\mathcal R(j),
\]
by Lemma~\ref{lem:norm-loop-orientation}; compare
\cite[Theorems~3.1.1 and 3.2.1]{Balderrama} at the zero basepoint.

Since $V$ is a shifted degree,
\[
 \dim(\Ind_K^GV)^G=\dim V^K=2a_{n-1}-1\in2\mathbb Z+1.
\]
Theorem~\ref{thm:local-padic-identification} therefore gives $\pi_0X(j)=\pi_0Y(j)=0$ at every object of the
truncated Tate compression diagram. Lemma~\ref{lem:tambara-pointed-path} gives
the map on $\pi_V^K(H)$ as
\[
 [z_0,\ldots,z_{n-2}]
 \longmapsto
 [a_K^GN_K^G(z_0),\ldots,a_K^GN_K^G(z_{n-2})].
\]
The target cokernel computes
$\pi_{\Ind_K^GV}^G(R_<)$.

Since $R_<$ is the limit in $\operatorname{Alg}_{\mathcal O}$,
its indexed multiplication $\mu_<$ satisfies
\[
\pi_j\mu_<=\mu_jN_K^G(\res_K^G\pi_j)
\]
for every $\pi_j:R_<\to\mathcal R(j)$.
Thus $(f_j)$ induces the norm
\[
N_{K,R_<}^G:
\pi_V^K(\res_K^G R_<)
\longrightarrow
\pi_{\Ind_K^GV}^G(R_<).
\]
Since $q:H\to R_<$ is an $\mathcal O$-algebra map,
\[
N_{K,R_<}^G\bigl((\res_K^Gq)_*[z]\bigr)
=q_*N_{K,H}^G([z]).
\]
As $\res_K^Gq$ is an equivalence, the truncated
cokernel class computed above is
$q_*N_{K,H}^G([z])$.

Multiplication by $v_{D_V}^{-1}$ transports
$\Ind_K^GV+1$ to $\rho_G(J(\mathbf a))$.
After this degree transport, $q_*$ is the quotient map
$\overline{\operatorname{pr}}^1:\coker\Delta^G_{J\mathbf a}\to
\coker\Delta^{G,<}_{J\mathbf a}$ induced by the chain projection
$(\operatorname{pr}^0,\operatorname{pr}^1)$ of
Lemma~\ref{lem:tambara-terminal-elimination}. The inverse of
$\overline{\operatorname{pr}}^1$ gives
\eqref{eq:tambara-euler-power}. Additivity follows from
Lemma~\ref{lem:tambara-pointed-path}.
\end{proof}

\begin{remark}
For $0\le i<n-1$, the $i$th entry is given by the composite
\[
 T_i^K(\mathbf a)
 \xrightarrow{\ N_K^G\ }\pi_{\Ind_K^G(V+1)}^G(H^{t_i})
 \xrightarrow{\ a_K^G\ }\pi_{\Ind_K^GV+1}^G(H^{t_i})
 \xrightarrow{\ v_{D_V}^{-1}\ }T_i^G(J(\mathbf a)).
\]
\end{remark}

\begin{remark}[The coefficient of $a_K^G$]
\label{rem:tambara-loop-euler-coefficient}
For odd $p$, use the standard complex orientations in
\[
 \mathbb R[G/K]-1
 \cong\bigoplus_{r=1}^{(p-1)/2}\lambda(rp^{n-1}).
\]
Write
\[
 D'=\sum_{r=1}^{(p-1)/2}\lambda(rp^{n-1})
 -\frac{p-1}{2}\lambda(p^{n-1}).
\]
Equation~\eqref{eq:character-euler-factor} gives
\[
 a_K^G=\left(\frac{p-1}{2}\right)!\,v_{-D'}a_{\lambda(p^{n-1})}^{(p-1)/2}.
\]
In the ordered real permutation basis, with $0\leq j<p$, the Fourier columns
\[
 1,\ \cos(2\pi j/p),\ \sin(2\pi j/p),\ldots,
 \cos(\pi(p-1)j/p),\ \sin(\pi(p-1)j/p)
\]
have determinant $p^{p/2}/2^{(p-1)/2}>0$. Indeed the complex Fourier
Vandermonde determinant has phase $(-1)^{(p-1)/2}i^{(p-1)/2}$ and absolute value
$p^{p/2}$; arranging the nonzero characters in conjugate pairs is an
even permutation, and each conversion from a conjugate pair to
cosine and sine has determinant $i/2$. The invariant diagonal line,
followed by the standard
complex orientations of the nontrivial pairs, therefore agrees with
the permutation orientation.

\end{remark}

\leadbreak
This gives the following corollary.

\begin{corollary}[Norms preserve the descent filtration]
\label{cor:tambara-shifted-visible}
For cyclic $p$-groups $e<K\le L$ and $V\in RO(K)$,
\[
 N_K^L(F^1\pi_V^K(H))\subseteq F^1\pi_{\Ind_K^LV}^L(H).
\]
On each shifted degree the adjacent norm is additive.
\end{corollary}

\begin{proof}
Theorem~\ref{thm:tambara-shifted-cokernel} proves additivity.  Every inclusion of cyclic
$p$-groups factors as a chain into adjacent inclusions, so we norm preserve filtration on every subgroup pairs.
\end{proof}

We now give the explicit norm formula for edge elements, it's a primary but heavy-workloaded caculation in valuation language. Surprisingly, the eventual formula is clean.

\begin{theorem}[Explicit norm formula in shifted degrees]
\label{thm:tambara-shifted-exact}
Let $K=C_{p^{n-1}}<G=C_{p^n}$, $n\ge2$, and let $V$
be a shifted degree with $V+1=\rho_K(\mathbf a)$.
Use $J(\mathbf a)$ and $D_V$ from
\eqref{eq:tambara-shifted-companion}--\eqref{eq:tambara-J-degree}.
For
$[q_0,\ldots,q_{n-2}]\in\coker\Delta^K_{\mathbf a}$, $q_i\in p^{b_{T_i^K}(\mathbf a)}\mathbb Z$,
we have
\begin{equation}\label{eq:tambara-shifted-exact-scalar}
N_K^G([q_0,\ldots,q_{n-2}])
=v_{D_V}[p^{a_{n-1}}q_0^p,\ldots,
p^{a_{n-1}}q_{n-2}^p,0].
\end{equation}
The entries on the two sides denote their classes in
$T_i^K(\mathbf a)$ and $T_i^G(J(\mathbf a))$,
respectively.
\end{theorem}

\begin{proof}
For $0\le i\le n-2$, Theorem~\ref{thm:exact-weighted-local-norm}
and Remark~\ref{rem:tambara-loop-euler-coefficient} give
\[
 v_{D_V}^{-1}a_K^GN_K^G([q_i])=
 \begin{cases}
 \displaystyle\left[\left(\frac{p-1}{2}\right)!
    p^{a_{n-1}}q_i^p\right],&p\text{ odd},\\[2mm]
 [\eta_{T_i}(\mathbf a)2^{a_{n-1}}q_i^2],&p=2.
 \end{cases}
\]
For odd $p$, \eqref{eq:local-height-table} gives
\[
\begin{aligned}
 b_{T_i^K}(\mathbf a)
   &=S_i(\mathbf a),\\
 b_{T_i^G}(J(\mathbf a))
   &=pS_i(\mathbf a)+a_{n-1}
     -\frac{(p-1)(n-i-1)}2,\\
 c_{T_i^G}(J(\mathbf a))
   &=b_{T_i^G}(J(\mathbf a))+n-i.
\end{aligned}
\]
If $p=3$, then
\[
\begin{gathered}
 \left(\frac{p-1}{2}\right)!=1,\qquad
 v_{D_V}^{-1}a_K^GN_K^G([q_i])=[p^{a_{n-1}}q_i^p],\\
 v_3(3^{a_{n-1}}q_i^3)
 \ge a_{n-1}+3b_{T_i^K}(\mathbf a)
 =b_{T_i^G}(J(\mathbf a))+n-i-1
 =c_{T_i^G}(J(\mathbf a))-1.
\end{gathered}
\]
If $p\ge5$, then $(p-1)(n-i-1)/2\ge n-i$, so
\[
\begin{aligned}
 v_p(p^{a_{n-1}}q_i^p)
 &\ge a_{n-1}+p\,b_{T_i^K}(\mathbf a)\\
 &=b_{T_i^G}(J(\mathbf a))
   +\frac{(p-1)(n-i-1)}2\\
 &\ge c_{T_i^G}(J(\mathbf a)),
\end{aligned}
\]
and hence
\[
 v_{D_V}^{-1}a_K^GN_K^G([q_i])
 =\left[\left(\frac{p-1}{2}\right)!p^{a_{n-1}}q_i^p\right]
 =[p^{a_{n-1}}q_i^p]=0.
\]
For $p=2$, \eqref{eq:local-height-table} and
\eqref{eq:norm-dyadic-edge-units} give
\[
\begin{aligned}
 b_{T_i^K}(\mathbf a)
  &=S_i(\mathbf a)+\frac{\epsilon(\mathbf a)(n-i-2)}2,\\
 b_{T_i^G}(J(\mathbf a))
  &=2S_i(\mathbf a)+a_{n-1},\\
 c_{T_i^G}(J(\mathbf a))
  &=b_{T_i^G}(J(\mathbf a))+1,\\
 \eta_{T_i}(\mathbf a)
  &=1+2^{n-1-i}t_i(\mathbf a)\equiv1\pmod2.
\end{aligned}
\]
Therefore
\[
\begin{aligned}
 v_2(2^{a_{n-1}}q_i^2)
 &\ge a_{n-1}+2b_{T_i^K}(\mathbf a)\\
 &=b_{T_i^G}(J(\mathbf a))
    +\epsilon(\mathbf a)(n-i-2)\\
 &\ge b_{T_i^G}(J(\mathbf a)),\\
 (\eta_{T_i}(\mathbf a)-1)2^{a_{n-1}}q_i^2
 &\in 2^{c_{T_i^G}(J(\mathbf a))}\mathbb Z,
\end{aligned}
\]
so
\[
 v_{D_V}^{-1}a_K^GN_K^G([q_i])
 =[\eta_{T_i}(\mathbf a)2^{a_{n-1}}q_i^2]
 =[2^{a_{n-1}}q_i^2].
\]
Thus, for every $p$ and $0\le i\le n-2$,
\[
 v_{D_V}^{-1}a_K^GN_K^G([q_i])=[p^{a_{n-1}}q_i^p].
\]
By \eqref{eq:tambara-euler-power},
\[
 N_K^G([q_0,\ldots,q_{n-2}])
 =v_{D_V}[p^{a_{n-1}}q_0^p,\ldots,
 p^{a_{n-1}}q_{n-2}^p,0].
\]
Theorem~\ref{thm:tambara-shifted-cokernel} gives independence
of the representatives $q_i$.
\end{proof}

The following extend the norm form $RO_\mathrm{sel}(K)$ to $RO(K)$.

\begin{proposition}
\label{prop:norm-full-degree-coherence}
Let $K\leq L$ be cyclic $p$-groups, $V\in RO(K)$, and
$V_0$ its selected representative. Put $D=V-V_0$.
For $x\in\pi_V^K(H)$,
\[
N_K^L(x)
=v_{\Ind_K^L D}N_K^L(v_D^{-1}x).
\]
Thus the norms on $RO_{\mathrm{sel}}(K)$ determine
all $RO(K)$-graded norms when $K=C_{p^n}$.
\end{proposition}
\begin{proof}
This follows from \eqref{eq:character-unit-coherence}
and multiplicativity.
\end{proof}

Square-zero ideals can have nonzero norms
\cite[Example~14.13]{StricklandTambara}.
The following corollary determines their possible indices
for $F^1\underline\pi_\star H$.

\begin{corollary}[Vanishing of norms]
\label{cor:shifted-norm-exponent-vanishing}
Let $K<L$ be cyclic $p$-groups and let
$m\in F^1\pi_\star^K(H)$ be homogeneous. Then
\[
pN_K^L(m)=0,
\qquad
\left(p\ge5\ \text{or}\ [L:K]\ge p^2\right)
\Longrightarrow N_K^L(m)=0.
\]
If $J<K<L$ and $m\in F^1\pi_\star^J(H)$, then
\begin{equation*}
N_K^L\!\left(N_J^K(m)\right)=N_J^L(m)=0.
\end{equation*}
Nonzero norms on $F^1$ occur only when
$p\in\{2,3\}$ and $[L:K]=p$, and both cases occur.
\end{corollary}

\begin{proof}
For $C_{p^{n-1}}<C_{p^n}$, $n\ge2$, the proof of
Theorem~\ref{thm:tambara-shifted-exact} give
\[
 p[p^{a_{n-1}}q_i^p]=0\quad(p=2,3),\qquad
 [p^{a_{n-1}}q_i^p]=0\quad(p\ge5).
\]
By \eqref{eq:tambara-shifted-exact-scalar},
$pN_K^L(m)=0$, and $N_K^L(m)=0$ for $p\ge5$.

For $p=2,3$, consider
$C_{p^{n-1}}<C_{p^n}<C_{p^{n+1}}$, $n\ge2$.
Put $V=\rho_{C_{p^{n-1}}}(\mathbf a)-1$ and choose
$q_i\in p^{b_{T_i^{C_{p^{n-1}}}}(\mathbf a)}\mathbb Z$ representing a class
\[
 [q_0,\ldots,q_{n-2}]\in\coker\Delta^{C_{p^{n-1}}}_{\mathbf a}.
\]
By \eqref{eq:tambara-shifted-exact-scalar},
\[
 v_{D_V}^{-1}N_{C_{p^{n-1}}}^{C_{p^n}}([q_0,\ldots,q_{n-2}])
 =[p^{a_{n-1}}q_0^p,\ldots,p^{a_{n-1}}q_{n-2}^p,0].
\]
Since $N_{C_{p^n}}^{C_{p^{n+1}}}(v_{D_V})$ is a unit,
\eqref{eq:tambara-shifted-exact-scalar} for
$C_{p^n}<C_{p^{n+1}}$ reduces the second norm to
\[
 [p^{(p+1)a_{n-1}}q_0^{p^2},\ldots,
 p^{(p+1)a_{n-1}}q_{n-2}^{p^2},0,0]
 \in\coker\Delta^{C_{p^{n+1}}}_{J(J(\mathbf a))}.
\]
For $0\le i\le n-2$,
\eqref{eq:tambara-shifted-companion} and
\eqref{eq:local-height-table} give
\[
\begin{aligned}
 S_i(J(\mathbf a))
  &=pS_i(\mathbf a)+a_{n-1}-\frac{p-1}{2}(n-i-1),\\
 c_{T_i^{C_{p^{n+1}}}}(J(J(\mathbf a)))
  &=p^2S_i(\mathbf a)+(p+1)a_{n-1}
    -\frac{p(p-1)}{2}(n-i-1)+1,\\
 v_p(q_i)-S_i(\mathbf a)
  &\ge\begin{cases}
     \epsilon(\mathbf a)(n-i-2)/2,&p=2,\\
     0,&p=3.
   \end{cases}
\end{aligned}
\]
Consequently,
\[
\begin{aligned}
 &v_p\bigl(p^{(p+1)a_{n-1}}q_i^{p^2}\bigr)
   -c_{T_i^{C_{p^{n+1}}}}(J(J(\mathbf a)))\\
 &\qquad\ge
 \begin{cases}
 (1+2\epsilon(\mathbf a))(n-i-2),&p=2,\\
 3(n-i-1)-1,&p=3
 \end{cases}
 \ge0.
\end{aligned}
\]
Thus every entry of the second norm is zero. Norm transitivity gives
the second formula in the corollary and $N_K^L(m)=0$ for
$[L:K]\ge p^2$. Proposition~\ref{prop:norm-full-degree-coherence}
extends these conclusions to all $RO$-degrees.

For $p=2,3$, take $C_p<C_{p^2}$ and $\mathbf a=(1,-1)$.
Then
\[
 A_0^{C_p}(\mathbf a)=A_1^{C_p}(\mathbf a)=0,\qquad
 \coker\Delta^{C_p}_{\mathbf a}=p^{-1}\mathbb Z/\mathbb Z,
\]
and \eqref{eq:tambara-shifted-exact-scalar} gives
\[
\begin{aligned}
 p=2:\quad&
 \coker\Delta^{C_4}_{J(\mathbf a)}
 \cong2^{-3}\mathbb Z/2^{-2}\mathbb Z,
 &v_{D_V}^{-1}N_{C_2}^{C_4}([2^{-1}])&=[2^{-3}],\\
 p=3:\quad&
 \coker\Delta^{C_9}_{J(\mathbf a)}
 \cong3^{-5}\mathbb Z/3^{-3}\mathbb Z,
 &v_{D_V}^{-1}N_{C_3}^{C_9}([3^{-1}])&=[3^{-4}].
\end{aligned}
\]
Both classes are nonzero, since
$\ell_0(A_0^{C_{p^2}}(J(\mathbf a)))=
\rho_0(A_1^{C_{p^2}}(J(\mathbf a)))=0$ for $p=2,3$.
\end{proof}

\subsection{Norm identities}
\label{sec:norm-identities}
This subsection we get some identities of norm.

Use the ordered bases and representation sphere models of
Section~\ref{sec:oriented-models} and
\eqref{eq:oriented-model-coefficient-transport}.
For $K\le L\le C_N$, $V\in RO(K)$ and $g\in L$, conjugation is
\[
 c_g:\pi_V^K(H)\longrightarrow\pi_{{}^gV}^K(H),
\]
where the $K$-action on ${}^gV$ is obtained from that on $V$
by $k\mapsto g^{-1}kg$.
For $U\in RO(L)$ and $V=\Res_K^L U$, the Weyl action of
$gK\in L/K$ on $\pi_V^K(H)$ is
\[
 x\longmapsto c_g(x)\circ S^{g^{-1}},
\]
where $g:{}^gV\xrightarrow{\cong}V$ is induced by the $L$-action on $U$.

For products and induction,
let $J\le K\le L$, $U\in RO(J)$, $V,W\in RO(K)$,
$x\in\pi_V^K(H)$, $y\in\pi_W^K(H)$ and $z\in\pi_U^J(H)$.
The graded Tambara identities \cite[Theorem~5.16]{AngeltveitBohmann} give
\[
 N_K^L(xy)=N_K^L(x)N_K^L(y),\qquad
 N_J^L(z)=N_K^L(N_J^K(z)),
\]
using \eqref{eq:oriented-model-coefficient-transport} for the
$L$-linear isomorphisms
$\Ind_K^L(V\oplus W)\cong\Ind_K^L V\oplus\Ind_K^L W$
and $\Ind_K^L\Ind_J^K U\cong\Ind_J^L U$, respectively.
For $d=[L:K]$, the first isomorphism has determinant sign
\[
 (-1)^{\binom d2\dim(V)\dim(W)}:
\]
the required permutation interchanges $\binom d2$ pairs consisting
of a $V$-summand and a $W$-summand.

For $J\le K\le L$ and $U\in RO(J)$, put $a=[K:J]$, $b=[L:K]$.
Iterated induction orders $(r,s)$ first by $r$, then by $s$,
where $0\le r<b$ and $0\le s<a$; direct induction orders by $r+bs$.
Thus $\Ind_K^L\Ind_J^K U\cong\Ind_J^L U$ has determinant sign
\begin{equation}\label{eq:norm-direct-iterated-sign}
 \chi_{a,b}(U)=(-1)^{\dim(U)\binom a2\binom b2}.
\end{equation}
Indeed, an inversion occurs exactly when $r<r'$ and $s>s'$.
Each of the $\binom a2\binom b2$ inversions contributes
$(-1)^{\dim U}$. Duality gives both signs for virtual representations.
For $U_1=(\Ind_J^K U)_0$ and $z\in\pi_U^J(H)$,
\[
\begin{aligned}
 &v_{\Ind_J^L U-(\Ind_J^L U)_0}^{-1}N_J^L(z)\\
 &\quad=\chi_{a,b}(U)\,
 v_{\Ind_K^L U_1-(\Ind_K^L U_1)_0}^{-1}
 N_K^L\bigl(v_{\Ind_J^K U-U_1}^{-1}N_J^K(z)\bigr).
\end{aligned}
\]
This follows from norm transitivity,
\eqref{eq:character-unit-coherence} and
\eqref{eq:oriented-model-determinant}.

For restrictions, 
let $T$ be the ordered representatives of $L/K$ and $P\le L$.
Order the $P$-orbits by their smallest coset indices and choose the
corresponding $g_O\in T$. Within each orbit, order
$s\in P/(P\cap K)$ by powers of the chosen generator of $P$.
Write
$sg_O=t k_{O,s}$, $t\in T$, $k_{O,s}\in K$,
and let $\pi$ be the permutation from the ordered labels $(O,s)$ to $t$.
The resulting $P$-linear isomorphism
\[
 \kappa_P:\bigoplus_O
 \Ind_{P\cap K}^{P}({}^{g_O}\!\Res_{P\cap K}^K V)
 \xrightarrow{\cong}\Res_P^L\Ind_K^L V
\]
is given, for actual $V$, by
$s\otimes v\mapsto t\otimes k_{O,s}v$ on the $O$-summand,
where $sg_O=tk_{O,s}$.
For virtual $V$, use the sphere models and duality of
Section~\ref{sec:oriented-models}.
Its determinant sign is
\begin{equation}\label{eq:norm-orbit-determinant}
 (\sgn\pi)^{\dim V}\prod_{O,s}\varepsilon_V(k_{O,s}).
\end{equation}
The block permutation contributes $(\sgn\pi)^{\dim V}$;
the action on each copy of $V$ contributes $\varepsilon_V(k_{O,s})$.

For $x\in\pi_V^K(H)$, the orbit decomposition
\[
 \Res_P^L(L/K)\cong
 \coprod_{[g]\in P\backslash L/K}P/(P\cap K)
\]
and the graded norm--restriction identity give
\begin{equation}\label{eq:tambara-norm-double-coset}
 \res_P^L N_K^L(x)=(\kappa_P)_*\left(
 \prod_{[g]\in P\backslash L/K}
 N_{P\cap K}^{P}(c_g\res_{P\cap K}^K x)\right).
\end{equation}
The map $(\kappa_P)_*$ is given by
\eqref{eq:oriented-model-coefficient-transport}, using
\eqref{eq:norm-orbit-determinant} for $\sgn\det(\kappa_P)$.

For the norm of a sum, let $x,y\in\pi_V^K(H)$ and put
$x_0=x$, $x_1=y$ and $T=L/K$.
For $t=\gamma_L^jK\in T$ with $0\leq j<[L:K]$,
the notation $ {}^tV$ means $ {}^{\gamma_L^j}V$.
The action $(\ell f)(t)=f(\ell^{-1}t)$ on $\{0,1\}^{T}$ has stabilizers
$H_f=\{\ell\in L:\ell f=f\}$.
Since $L$ is abelian, $K\le H_f$, and $f$ is constant on each
$H_f$-orbit $O$. Use the orders specified above with $P=H_f$ to form
\[
 \kappa_f:\bigoplus_{O\in H_f\backslash T}
 \Ind_K^{H_f}({}^{g_O}V)
 \xrightarrow{\cong}\Res_{H_f}^L\Ind_K^L V.
\]
Then
\begin{equation}\label{eq:norm-explicit-exponential}
 N_K^L(x+y)=\sum_{[f]\in L\backslash\{0,1\}^{T}}
 \tr_{H_f}^L\left((\kappa_f)_*\left(
 \prod_{O\in H_f\backslash T}
 N_K^{H_f}(c_{g_O}x_{f(g_OK)})\right)\right),
\end{equation}
where $(\kappa_f)_*$ is given by
\eqref{eq:oriented-model-coefficient-transport}, using
\eqref{eq:norm-orbit-determinant} with $P=H_f$
for $\sgn\det(\kappa_f)$.
Apply graded Tambara distributivity
\cite[Section~5.4 and Theorem~5.16]{AngeltveitBohmann}
to $\{0,1\}^T$, using $\kappa_f$ for the representation degrees.
Compare \cite[Theorem~2.5]{HillMazur} in degree zero.

\subsection{Assembly of norm to cyclic case}
\label{sec:norm-sylow-assembly}
The last mission is to assenbly the norm from $C_{p^n}$-equivariant case to $C_N$-equivariant case.

Let $K\leq L\leq C_N$, $V\in RO(K)$ and $x\in\pi_V^K(H)$.
Write $W=\Ind_K^L V$ and assume $L\ne e$. For every prime $p\mid |L|$, let
$P=P_p(L)$ and $K_p=K\cap P$. 
Use the conjugation maps $c_g$ and the representation isomorphisms
$\kappa_P$ as in Section~\ref{sec:norm-identities}.

\begin{theorem}[Assembly of norm to cyclic case]
\label{thm:prime-local-cyclic-norms}
The product of the restriction maps
\begin{equation*}
 \pi_W^L(H)\longrightarrow
 \prod_{p\mid |L|}\pi_{\Res_P^L W}^{P}(H),
 \qquad
 y\longmapsto(\res_P^Ly)_p,
\end{equation*}
is injective.
Choose integers $b_p$ satisfying
$\sum_{p\mid|L|}b_p[L:P_p(L)]=1$.
Let $y_p$ be defined in Proposition~\ref{prop:sylow-norm-coordinates}.
Then $N_K^L(x)$ is the unique $y\in\pi_W^L(H)$ such that
$\res_P^Ly=y_p$ for every $p\mid|L|$, and
\begin{equation*}
\quad
 N_K^L(x)=\sum_{p\mid|L|}b_p\tr_{P_p(L)}^L(y_p).
 \quad
\end{equation*}
\end{theorem}
\begin{proof}
Since $\sum_{p\mid|L|}b_p[L:P]=1$,
Proposition~\ref{prop:cohomological-hz} gives, for any $y\in\pi_W^L(H)$,
\[
 y=\sum_{p\mid|L|}b_p\tr_P^L\res_P^L(y).
\]
Hence $\prod_p\res_P^L$ is injective.

Equation~\eqref{eq:tambara-norm-double-coset} gives
$y_p=\res_P^LN_K^L(x)$.
Substituting $y=N_K^L(x)$ in
$y=\sum_{p\mid|L|}b_p\tr_P^L\res_P^L(y)$ proves the formula for
$N_K^L(x)$, independently of the choice of the integers $b_p$.
The family $(y_p)_p$ therefore satisfies the compatibility conditions of
Theorem~\ref{thm:integral-cyclic-reconstruction}.
\end{proof}

We have the following.

\begin{proposition}\label{prop:sylow-norm-coordinates}
For every prime $p\mid|L|$,
\begin{equation*}
 \res_P^L N_K^L(x)=
 y_p:=(\kappa_P)_*\left(
 \prod_{[g]\in P\backslash L/K}
 N_{K_p}^{P}
       (c_g\res_{K_p}^{K}x)\right).
\end{equation*}
The map $(\kappa_P)_*$ is given by
\eqref{eq:oriented-model-coefficient-transport}, using
\eqref{eq:norm-orbit-determinant} for $\sgn\det(\kappa_P)$.
\end{proposition}
\begin{proof}
Apply \eqref{eq:tambara-norm-double-coset} with $P=P_p(L)$.
Theorems~\ref{thm:tambara-natural-kernel} and~\ref{thm:tambara-shifted-exact}
compute $N_{C_{p^{n-1}}}^{C_{p^n}}$ for $n\ge2$.
For $m\in\pi_0^e(H)=\mathbb Z$, use $N_e^{C_p}(m)=m^p$.
Norm transitivity gives $N_{K_p}^P$ by iteration, with the signs
in \eqref{eq:norm-direct-iterated-sign}.
\end{proof}

\begin{corollary}[Cyclic norms]\label{cor:cyclic-norm-determination}
The $RO$-graded Tambara functor of $H$ for $C_N$ is determined by
the corresponding functors for $C_{p^n}\leq C_N$ and the Weyl actions.
\end{corollary}
\begin{proof}
Apply Theorem~\ref{thm:prime-local-cyclic-norms} and
Proposition~\ref{prop:sylow-norm-coordinates}.
\end{proof}

\subsection{Completion}\label{sec:sheaf}
In this subsection, we prove that norms preserve the descent filtration
and extend Corollary~\ref{cor:shifted-norm-exponent-vanishing}
to arbitrary finite cyclic groups.
We then prove a norm-of-sum formula for square-zero Green ideals
preserved by norms in $RO$-graded Tambara functors
for finite abelian groups.

Let $G=C_N$.
For $P=C_{p^r}\leq G$, $r\geq1$, and $V\in RO(P)$,
define as in \eqref{eq:vertex-tensor-descent}
\[
\begin{aligned}
 B_P&=\left(\prod_{i=0}^{r-1}(\res_P^G H)^{h_i}\right)
 \times(\res_P^G H)^{\Phi_r},\\
 F^1\pi_V^P(H)&=
 \ker\bigl(\pi_V^P(H)\longrightarrow\pi_V^P(B_P)\bigr).
\end{aligned}
\]
Set $F^1\pi_*^e(H)=0$. For nontrivial $K\leq G$, write
\begin{equation}\label{eq:cyclic-descent-ideal}
 F^1\pi_V^K(H)=\bigcap_{p\mid|K|}
 (\res_{P_p(K)}^K)^{-1}
 F^1\pi_{\Res_{P_p(K)}^K V}^{P_p(K)}(H).
\end{equation}

\begin{theorem}[Cyclic descent filtration]
\label{thm:first-order-normed-descent}
The filtration
$0=F^2\underline\pi_\star H\subset
 F^1\underline\pi_\star H\subset
 F^0\underline\pi_\star H=\underline\pi_\star H$
is multiplicative. For $K\leq L\leq G$ and $V\in RO(K)$,
\[
 (F^1)^2=0,\qquad
 N_K^L(F^1\pi_V^K(H))\subseteq
 F^1\pi_{\Ind_K^L V}^L(H).
\]
Restrictions, transfers and Weyl actions also preserve $F^1$.
\end{theorem}
\begin{proof}
For cyclic $p$-groups, Theorems~\ref{thm:filtered-green-descent}
and \ref{thm:padic-homotopy-ring} show that $F^1$ is a square-zero
Green ideal, and Corollary~\ref{cor:tambara-shifted-visible}
shows that norms preserve $F^1$.
The homomorphisms induced by morphisms in the grading category
for $P$ of \cite[Section~3.1]{AngeltveitBohmann} are given by
precomposition with stable equivalences of representation spheres.
The maps $\pi_V^P(H)\to\pi_V^P(B_P)$ commute with this precomposition.
Hence these homomorphisms preserve $F^1$.
For general cyclic groups, \eqref{eq:cyclic-descent-ideal}
and Theorem~\ref{thm:sylow-reconstruction} show that $F^1$ is a
square-zero ideal preserved by restrictions and Weyl actions.
For $P=P_p(L)$, $K_p=K\cap P$, $W\in RO(L)$, and
$m\in F^1\pi_{\Res_K^L W}^K(H)$, the Mackey double-coset formula gives
\[
 \res_P^L\tr_K^L(m)=
 \sum_{[g]\in P\backslash L/K}
 \tr_{K_p}^P\bigl(c_g(\res_{K_p}^K m)\circ S^{g^{-1}}\bigr)\in F^1.
\]
Here $g:{}^g\!\Res_{K_p}^L W\xrightarrow{\cong}\Res_{K_p}^L W$
is induced by the $L$-action on $W$.
For $m\in F^1\pi_V^K(H)$,
\begin{equation*}
 \res_P^L N_K^L(m)=(\kappa_P)_*\left(
 \prod_{[g]\in P\backslash L/K}
 N_{K_p}^P(c_g\res_{K_p}^K m)\right),
 \qquad |P\backslash L/K|=[L:PK].
\end{equation*}
If $K_p=e$, each factor is zero. Otherwise each factor
lies in $F^1$, and the product vanishes when $[L:PK]\geq2$.
Thus transfers and norms preserve $F^1$ by
\eqref{eq:cyclic-descent-ideal}.
\end{proof}

\begin{corollary}[Vanishing of norm]
\label{cor:cyclic-proper-ideal-norms}
For $K<L$, $d=[L:K]$, and $m\in F^1\pi_V^K(H)$,
\begin{equation*}
 dN_K^L(m)=0,\qquad
 d\notin\{2,3\}\ \Longrightarrow\ N_K^L(m)=0.
\end{equation*}
For $J<K<L$ and $m\in F^1\pi_V^J(H)$, we have in degree $\Ind_J^L V$
\begin{equation*}
 N_K^LN_J^K(m)=N_J^L(m)=0.
\end{equation*}
Nonzero norms on $F^1$ occur for both $d=2$ and $d=3$.
\end{corollary}
\begin{proof}
Set $P=P_p(L)$ and $K_p=K\cap P$.
The proof of Theorem~\ref{thm:first-order-normed-descent} gives
$\res_P^L N_K^L(m)=0$ if $K_p=e$ or $[L:PK]\geq2$.
Otherwise $d=[P:K_p]>1$ is a power of $p$.
Corollary~\ref{cor:shifted-norm-exponent-vanishing} gives
$N_{K_p}^P(c_g\res_{K_p}^K m)=0$ unless $d=p\in\{2,3\}$,
and then $pN_{K_p}^P(c_g\res_{K_p}^K m)=0$. By Theorem~\ref{thm:sylow-reconstruction},
$dN_K^L(m)=0$, and $N_K^L(m)=0$ when $d\notin\{2,3\}$.
Norm transitivity and $[L:J]=[L:K][K:J]\geq4$ give
$N_K^LN_J^K(m)=N_J^L(m)=0$.
Examples for both exceptional indices are given in
Corollary~\ref{cor:shifted-norm-exponent-vanishing}.
\end{proof}
Let $G$ be a finite abelian group.
For an $RO$-graded Tambara functor $T$, write $T_K(V)$
for its value at $G/K$ in degree $V\in RO(K)$.
\begin{definition}
\label{def:norm-stable-square-zero}
For every $K\leq G$ and $V\in RO(K)$, let
$I_K(V)\subseteq T_K(V)$ be an additive subgroup.
Assume that every homomorphism $T_K(V)\to T_K(W)$
induced by a morphism in the grading category for $K$ of
\cite[Section~3.1]{AngeltveitBohmann}
carries $I_K(V)$ into $I_K(W)$.
For every inclusion $K\leq L$, assume the following, with
$U,V,W\in RO(K)$ in the first line, $V\in RO(L)$ in the second,
and $V\in RO(K)$, $g\in G$ in the third:
\begin{equation*}
\begin{gathered}
 T_K(U)I_K(V)\subseteq I_K(U+V),\qquad I_K(V)I_K(W)=0,\\
 \res_K^L I_L(V)\subseteq I_K(\Res_K^L V),\qquad
 \tr_K^L I_K(\Res_K^L V)\subseteq I_L(V),\\
 c_g I_K(V)\subseteq I_K({}^gV),\qquad
 N_K^L I_K(V)\subseteq I_L(\Ind_K^L V).
\end{gathered}
\end{equation*}
\end{definition}
For $K\leq Q\leq L$ and $V\in RO(K)$, ordered representatives
of $L/Q$ containing $1$ give the $Q$-linear isomorphism
\begin{equation*}
 \left(\bigoplus_{\substack{\bar g\in L/Q\\\bar g\ne\bar e}}
 {}^g\!\Ind_K^Q V\right)\oplus\Ind_K^Q V
 \xrightarrow{\cong}\Res_Q^L\Ind_K^L V.
\end{equation*}

Hill--Mazur give the ungraded norm-of-sum formula
\cite[Theorem~2.5]{HillMazur}.
The norms $N_K^L$ need not vanish on $I_K(V)$.

\begin{theorem}
\label{thm:finite-abelian-square-zero-norm}
Assume the hypotheses of Definition~\ref{def:norm-stable-square-zero}.
For
$K\leq L$, $V\in RO(K)$, $x\in T_K(V)$ and $m\in I_K(V)$, then
\begin{equation*}
\begin{aligned}
 N_K^L(x+m)={}&N_K^L(x)\\
 &+\sum_{K\leq Q\leq L}\tr_Q^L\left[
 \left(\prod_{\substack{\bar g\in L/Q\\\bar g\ne\bar e}}
 c_g N_K^Q(x)\right)N_K^Q(m)\right].
\end{aligned}
\end{equation*}
For each $K\leq Q\leq L$, apply the homomorphism induced by
the above $Q$-linear isomorphism to the product before
applying $\tr_Q^L$.
\end{theorem}
\begin{proof}
The norm-of-sum theorem \cite[Theorem~2.5]{HillMazur}
indexes terms by $L$-orbits of subsets $S\subseteq L/K$.
Let $Q=\operatorname{Stab}_L(S)$. The subset $S$ is a union
of $Q$-orbits. Its contribution transfers from $Q$ the product
of $c_gN_K^Q(m)$ on the orbits in $S$ and
$c_gN_K^Q(x)$ on the other orbits.
Because $I_Q^2=0$, only the empty subset and subsets consisting
of one $Q$-orbit survive. The empty subset contributes
$N_K^L(x)$; for fixed $Q$, the latter subsets form one
$L$-orbit, represented by $Q/K$. Its transfer is the $Q$-term
of the formula in the theorem.
\end{proof}

For finite cyclic $G$, take $T=\underline\pi_\star H$ and $I=F^1$.
Theorem~\ref{thm:first-order-normed-descent} and its proof verify
the hypotheses of Definition~\ref{def:norm-stable-square-zero},
including preservation of $F^1$ by the homomorphisms induced
by morphisms in the grading category.
Thus Theorem~\ref{thm:finite-abelian-square-zero-norm} shows that
$\underline\pi_\star H/F^1$ inherits an $RO$-graded Tambara
structure. Unlike the square-zero extension constructed in
\cite[Proposition~3.2]{HillAQ}, the norms $N_K^L$ with $K<L$
need not vanish on $F^1$.

\begin{corollary}[Restrictions of norms]
\label{cor:isotropy-support-law}
\label{cor:top-isotropy-concentration}
Let $G$ be a finite abelian group, and let $T,I$ satisfy
Definition~\ref{def:norm-stable-square-zero}.
Let $J,K\leq L\leq G$, $V\in RO(K)$ and $m\in I_K(V)$.
Choose double-coset representatives containing $1$ and let
\[
 \kappa_J:\bigoplus_{[g]\in J\backslash L/K}
 \Ind_{J\cap K}^{J}({}^g\!\Res_{J\cap K}^{K}V)
 \xrightarrow{\cong}\Res_J^L\Ind_K^L V
\]
be the corresponding $J$-linear isomorphism. Then
\[
\begin{aligned}
 JK<L&\ \Longrightarrow\ \res_J^L N_K^L(m)=0,\\
 JK=L&\ \Longrightarrow\ 
 \res_J^L N_K^L(m)=(\kappa_J)_*
 N_{J\cap K}^{J}(\res_{J\cap K}^K m).
\end{aligned}
\]
In particular, if $K<L$ are cyclic $p$-groups and
$m\in F^1\pi_V^K(H)$, then, for $J<L$,
\begin{equation*}
 \res_J^L N_K^L(m)=0.
\end{equation*}
\end{corollary}
\begin{proof}

The graded norm--restriction formula
\cite{AngeltveitBohmann} gives

\[
 \res_J^L N_K^L(m)=(\kappa_J)_*\left(
 \prod_{[g]\in J\backslash L/K}
 N_{J\cap K}^{J}(c_g\res_{J\cap K}^{K}m)\right).
\]
Each factor lies in $I_J$, and $|J\backslash L/K|=[L:JK]$.
For $JK<L$, the product vanishes by $I_J^2=0$;
for $JK=L$, the only representative is $g=1$.
For proper $K,J<L$ in a cyclic $p$-group, $JK<L$.
\end{proof}
\appendix
\section{Coefficient calculations and comparison maps}\label{app:local-calculations}
This appendix provides the detailed calculations used in the main text.
The main tool we use is the generalized homotopy fixed point spectral sequence.

\subsection{Localizations and the skeletal tower for $p$ odd}

Let $p$ be odd and $G=C_{p^n}$. Use the oriented planes
$\lambda_j$, the selected degree $RO_{\mathrm{sel}}(G)$,
the families $\mathcal F_i$ of Section~\ref{sec:local},
and the spectra in \eqref{eq:layer-definitions}. Since
$\widetilde E\mathcal F_i\simeq S^{\infty\lambda_i}$,
\begin{equation}\label{eq:layer-localizations}
 H^{\Phi_{i+1}}\simeq H^{\Phi_i}[a_{\lambda_i}^{-1}],
 \qquad H^{t_i}\simeq H^{h_i}[a_{\lambda_i}^{-1}].
\end{equation}
We use the presentation of $\pi_\star^G(H)$ in the positive cone
\cite[Proposition~5.7(5)]{Zeng}.

Write $Y^{(m)}$ for the $m$-skeleton of a $G$-CW complex $Y$.
The tower
$F(((E\mathcal F_i)_+)^{(m)},H^{\Phi_i})$, $m\geq0$, defines the
generalized homotopy fixed-point spectral sequence (HFPSS).
The Euler and orientation classes have degrees
$|a_{\lambda_j}|=-\lambda_j$, $|u_{\lambda_j}|=2-\lambda_j$.
The Euler relations are
\begin{equation*}
 p^{n-j}a_{\lambda_j}=0\qquad(0\leq j<n).
\end{equation*}
For $0\leq j<k<n$, the gold relations are
\begin{align*}
 a_{\lambda_k}u_{\lambda_j}
   &=p^{k-j}a_{\lambda_j}u_{\lambda_k}.
\end{align*}
For $1\leq i\leq n-1$, the equivalence
$\widetilde E\mathcal{F}_{i-1}\simeq S^{\infty\lambda_{i-1}}$ gives
\[
 H^{\Phi_i}\simeq H[a_{\lambda_{i-1}}^{-1}].
\]

We start with the odd-prime case.

\begin{theorem}\label{thm:odd-main}
For $p$ odd and $1\leq i<n$, we have an $RO_{\mathrm{sel}}(G)$-graded ring isomorphism
\begin{equation}
\pi_\star^{G}H^{h_i}\cong
\frac{
 (\mathbb{Z}/p^{n-i+1})[
  \{a_{\lambda_j}\}_{i\leq j\leq n-1},
  \{a_{\lambda_j}^{\pm1}\}_{0\leq j<i},
  \{u_{\lambda_j}^{\pm1}\}_{i\leq j\leq n-1},
  \{u_{\lambda_j}\}_{0\leq j<i}]
}{\left(\begin{gathered}
 p^{n-j}a_{\lambda_j}\quad(i\leq j<n),\\
 a_{\lambda_k}u_{\lambda_j}-p^{k-j}a_{\lambda_j}u_{\lambda_k}
 \quad(0\leq j<k<n)
\end{gathered}\right)},                                      \label{eq:odd-main-ring}
\end{equation}
where empty index sets are omitted.
\end{theorem}

\begin{remark}\label{prop:odd-index}
Lemma~\ref{lem:odd-normal} gives
$\operatorname{char}\pi_\star^G(H^{h_i})=p^{n-i+1}$.
\end{remark}

\subsection{The generalized HFPSS for $p$ odd}

Fix $n\geq2$ and $1\leq i<n$. Put
$r=n-i$, $K=C_{p^i}$ and $W=G/K=C_{p^r}$.
Write $\gamma$ for the image in $W$ of the chosen generator of $G$.
The model $E\mathcal F_i\simeq S(\infty\lambda_i)$ has one $G/K$-cell
in each dimension. Put $K_{-1}=*$ and
$K_m=((E\mathcal F_i)_+)^{(m)}$ for $m\geq0$. Then
$K_m/K_{m-1}\simeq G/K_+\wedge S^m$.

An overline denotes the image under
$\pi_\star^G(H)\to\pi_\star^G(H^{\Phi_i})\xrightarrow{\res_K^G}
\pi_\star^K(H^{\Phi_i})$.
We write bidegrees as $(V,s)$, with representation degree $V$,
filtration $s$, and total degree $V-s$. Define
\begin{align*}
 L_i^{(p)}&:=\mathbb{F}_p[
  \overline{a}_{\lambda_0}^{\pm1},\ldots,
  \overline{a}_{\lambda_{i-1}}^{\pm1},
  \overline{u}_{\lambda_i}^{\pm1},\ldots,
  \overline{u}_{\lambda_{n-1}}^{\pm1}],\\
 q&:=\frac{\overline{u}_{\lambda_{i-1}}}
          {\overline{a}_{\lambda_{i-1}}},
 \qquad |q|=(2,0).
\end{align*}
We first consider the spectral sequence.

\begin{proposition}[The spectral sequence]
\label{prop:odd-E2}
For $p$ odd and $1\leq i<n$, the tower $F(K_m,H^{\Phi_i})$
gives a multiplicative spectral sequence
\begin{equation}
 E_2^{V,s}=H^s\!\left(W;\pi_V^{K}(H^{\Phi_i})\right)
 \Longrightarrow \pi_{V-s}^{G}H^{h_i},
 \qquad |d_k|=(k-1,k).                                        \label{eq:odd-hfpss}
\end{equation}
\begin{equation*}
 \pi_\star^{K}(H^{\Phi_i})\cong L_i^{(p)}[q],                    
\end{equation*}
The $W$-action on $\pi_\star^K(H^{\Phi_i})$ is trivial.
Write $\Lambda(s)$ for the exterior algebra over $\mathbb F_p$
on $s$. Then
\[
 H^*(W;\mathbb{F}_p)=\mathbb{F}_p[t]\otimes\Lambda(s),
 \qquad |s|=(0,1),\quad |t|=(0,2),
\]
where $c_1$ is the first Chern class and
$t=c_1(\lambda_i)\bmod p$ for the $W$-representation
$\lambda_i$. Hence
\begin{equation*}
 E_2=L_i^{(p)}[q,t]\otimes\Lambda(s).                         
\end{equation*}
\end{proposition}
\begin{proof}
The cellular filtration gives \eqref{eq:odd-hfpss}. Put $K'=C_{p^{i-1}}$. Since
$\mathcal{F}_{i-1}\cap K$ is the family $\mathcal P_K$ of proper subgroups of
$K$, one has
\[
 \res_K^G(H^{\Phi_i})\simeq\widetilde E\mathcal P_K\wedge\res_K^G H.
\]
Moreover $E\mathcal P_K\simeq S(\infty\lambda_{i-1}|_K)$ is the inflation of
a free contractible $K/K'\cong C_p$-space.  Every cell has isotropy $K'$.
The cellular homology spectral sequence for
$(E\mathcal P_K)_+\wedge\res_K^G H$ has a single nonzero row,
given by $\pi_0^{K'}H=\mathbb Z$ with trivial $K/K'$-action.
Its cellular chain complex at the $K/K$-orbit is
\[
 C_m=\mathbb Z\quad(m\geq0),\qquad
 d_{2a}=p\ (a\geq1),\qquad d_{2a+1}=0\ (a\geq0).
\]
Thus
\[
 \pi_d^{K}((E\mathcal P_K)_+\wedge\res_K^G H)
 \cong H_d(BC_p;\mathbb Z).
\]
The collapse $(E\mathcal P_K)_+\to S^0$ induces
\[
 \pi_0^K((E\mathcal P_K)_+\wedge\res_K^G H)
 \longrightarrow \pi_0^K(H),\qquad
 \mathbb Z\xrightarrow{\ \tr_{K'}^K=p\ }\mathbb Z.
\]
The cofiber sequence
\[
 (E\mathcal P_K)_+\wedge\res_K^G H\longrightarrow\res_K^G H
 \longrightarrow\widetilde E\mathcal P_K\wedge\res_K^G H
\]
gives
\begin{equation*}
 \pi_d^{K}(H^{\Phi_i})\cong
 \begin{cases}
  \mathbb F_p,&d=0,2,4,\ldots,\\
  0,&d<0\text{ or }d\text{ odd}.
 \end{cases}                                                   
\end{equation*}
Localize the presentation in
\cite[Proposition~5.7(5)]{Zeng} at $a_{\lambda_{i-1}|_K}$.
The Euler relation becomes $p=0$, and
$q=\frac{\overline u_{\lambda_{i-1}}} {\overline a_{\lambda_{i-1}}} \in\pi_2^K(H^{\Phi_i})$.
The integer-graded subring is $\mathbb F_p[q]$,
with $q^m$ generating degree $2m$.

For $j<i$, one has $(\lambda_j|_K)^K=0$. The Euler cofiber sequence
\[
 S(\lambda_j|_K)_+\longrightarrow S^0
 \xrightarrow{a_{\lambda_j}}S^{\lambda_j|_K}
\]
has first term with proper $K$-isotropy, so
\[
 S(\lambda_j|_K)_+\wedge\widetilde E\mathcal P_K
 \wedge\res_K^G H\simeq0.
\]
Thus $\overline a_{\lambda_j}$ is a unit. For $k\geq i$,
$\lambda_k|_K\cong2$ as oriented representations, so
$\overline u_{\lambda_k}$ is a unit.

For $V=d+\sum_jc_j\lambda_j\in RO_{\mathrm{sel}}(G)$,
multiplication by
$\prod_{j<i}\overline a_{\lambda_j}^{c_j}
 \prod_{j\geq i}\overline u_{\lambda_j}^{c_j}$ gives an isomorphism
$\pi_V^K(H^{\Phi_i})\cong
\pi_{d+2\sum_{j\geq i}c_j}^K(H^{\Phi_i})$.
Since $1,\lambda_0,\ldots,\lambda_{n-1}$ form a
$\mathbb Z$-basis of $RO_{\mathrm{sel}}(G)$, these exponents
are uniquely determined by $V$. Together with the integer-graded
subring $\mathbb F_p[q]$ computed above, this gives the
monomial $\mathbb F_p$-basis of $L_i^{(p)}[q]$.

Every generator is a product or inverse of restrictions of
classes in $\pi_\star^G(H)$, so the $W$-action is trivial.
The periodic resolution of $W=C_{p^r}$ gives
$E_2=L_i^{(p)}[q,t]\otimes\Lambda(s)$.
\end{proof}

The cellular chains give the periodic free $\mathbb Z[W]$-resolution
$P_*\to\mathbb Z$:
\begin{equation*}
 \cdots\longrightarrow\mathbb{Z}[W]
 \xrightarrow{N}\mathbb{Z}[W]
 \xrightarrow{1-\gamma}\mathbb{Z}[W]
 \xrightarrow{N}\mathbb{Z}[W]
 \xrightarrow{1-\gamma}\mathbb{Z}[W],
 \quad N=1+\gamma+\cdots+\gamma^{p^r-1}.                     
\end{equation*}
The complex $\operatorname{Hom}_{\mathbb Z[W]}
(P_*,\pi_V^K(H^{\Phi_i}))$ computes $H^s(W;\pi_V^K(H^{\Phi_i}))$.

Since $\lambda_i|_K\cong2$, multiplication by $u_{\lambda_i}$
gives an equivalence
\[
 \Sigma^{2-\lambda_i}F(G/K_+\wedge S^m,H^{\Phi_i})
 \xrightarrow{\ u_{\lambda_i}\ }
 F(G/K_+\wedge S^m,H^{\Phi_i}).
\]
Induction over the cells gives compatible equivalences on
$F(K_m,H^{\Phi_i})$. Thus $x=a_{\lambda_i}u_{\lambda_i}^{-1}$
defines compatible classes in $\pi_{-2}^G F(K_m,H^{\Phi_i})$.

Define
\[
 \Lambda_i^{(p)}:=\mathbb{Z}[
  a_{\lambda_0}^{\pm1},\ldots,a_{\lambda_{i-1}}^{\pm1},
  u_{\lambda_i}^{\pm1},\ldots,u_{\lambda_{n-1}}^{\pm1}]
\]
and, with formal variables $|q|=2$ and $|x|=-2$,
\begin{equation}
 B_{r,i}^{(p)}:=
 \frac{\Lambda_i^{(p)}[q,x]}
 {(p^{r+1},\;qx-p,\;p^rx)}.                                  \label{eq:odd-normalized}
\end{equation}

We need the following lemma.

\begin{lemma}\label{lem:odd-normal}
Let $\mathcal M_i^{(p)}$ be the $\mathbb Z$-basis of
$\Lambda_i^{(p)}$ consisting of Laurent monomials.
For $B_{r,i}^{(p)}$ in \eqref{eq:odd-normalized}, we have
\begin{equation}
 B_{r,i}^{(p)}\cong
 \bigoplus_{w\in\mathcal M_i^{(p)}}
 \left(
  \bigoplus_{A\geq0}\mathbb{Z}/p^{r+1}\{wq^A\}
  \oplus
  \bigoplus_{B\geq1}\mathbb{Z}/p^r\{wx^B\}
 \right).                                                     \label{eq:odd-normal-sum}
\end{equation}
\end{lemma}

\begin{proof}
Let
\[
 M_r:=
 \bigoplus_{A\geq0}\mathbb Z/p^{r+1}\{Q_A\}
 \oplus
 \bigoplus_{B\geq1}\mathbb Z/p^r\{X_B\}.
\]
Define commuting endomorphisms $q$ and $x$ of $M_r$ by
\[
 qQ_A=Q_{A+1},\qquad
 xQ_0=X_1,\qquad xQ_A=pQ_{A-1}\ (A\geq1),
\]
\[
 xX_B=X_{B+1},\qquad
 qX_1=pQ_0,\qquad qX_B=pX_{B-1}\ (B\geq2).
\]
These maps satisfy
\[
 qx=xq=p,\qquad p^{r+1}=0,\qquad p^rx=0.
\]
Hence $M_r$ is a module over
$\mathbb Z[q,x]/(p^{r+1},qx-p,p^rx)$, generated by $Q_0$.
Since $qx=p$ reduces every mixed monomial to a multiple of $q^A$ or
$x^B$, the map
\[
 M_r\longrightarrow
 \mathbb Z[q,x]/(p^{r+1},qx-p,p^rx),
 \qquad Q_A\mapsto q^A,\quad X_B\mapsto x^B,
\]
is surjective, with inverse $f\mapsto fQ_0$.
Tensoring with $\Lambda_i^{(p)}$ gives \eqref{eq:odd-normal-sum}.
\end{proof}

In the setting of Proposition~\ref{prop:odd-E2}, let
\[
 X:=H[\{a_{\lambda_j}^{-1}\}_{j<i},
          \{u_{\lambda_j}^{-1}\}_{j\geq i}].
\]
Put $q=u_{\lambda_{i-1}}/a_{\lambda_{i-1}}$ and
$x=a_{\lambda_i}/u_{\lambda_i}$.

\leadbreak
We also need the following result.

\begin{lemma}[Localization and completion]\label{lem:odd-weighted}
The maps
\[
 X\xrightarrow{\simeq}F((E\mathcal F_i)_+,X)
   \xleftarrow{\simeq}H^{h_i}
\]
are equivalences. The left map is induced by
$(E\mathcal F_i)_+\to S^0$, and the right map is
$F((E\mathcal F_i)_+,H^{\Phi_i}\to X)$.
Moreover,
\begin{equation}\label{eq:odd-localized-core}
 \pi_\star^{G}X\cong
 \frac{\mathbb Z[
 \{a_{\lambda_j}^{\pm1}\}_{j<i},
 \{u_{\lambda_j}^{\pm1}\}_{j\geq i},q,x]}
 {(p^{r+1},\ qx-p,\ p^rx)}.
\end{equation}
\end{lemma}
\begin{proof}
The positive-cone presentation \cite[Proposition~5.7(5)]{Zeng} gives
\[
 \bigoplus_{d\in\mathbb Z,\ b_j\geq0}
 \pi_{d-\sum_jb_j\lambda_j}^{G}H
 =\frac{\mathbb Z[a_{\lambda_j},u_{\lambda_j}\mid0\leq j<n]}
 {\bigl(p^{n-j}a_{\lambda_j},\
 a_{\lambda_k}u_{\lambda_j}
       -p^{k-j}a_{\lambda_j}u_{\lambda_k}\ (j<k)\bigr)}.
\]
In particular, inverting $a_{\lambda_{i-1}}$ imposes $p^{r+1}=0$.
For $V=d+\sum_jc_j\lambda_j$ and integers $k_j\geq0$,
a fraction with denominator
$\prod_{j<i}a_{\lambda_j}^{k_j} \prod_{j\geq i}u_{\lambda_j}^{k_j}$
has degree $V$ when its numerator has degree
$d+2\sum_{j\geq i}k_j+\sum_j(c_j-k_j)\lambda_j$.
The multi-indices $(k_0,\ldots,k_{n-1})$ satisfying
$k_j\geq\max(c_j,0)$ form a cofinal subset of
$\mathbb Z_{\geq0}^n$ with its coordinatewise order.
For these indices the numerator degree lies in the positive cone.
Thus the localized presentation computes $\pi_\star^G X$.
For $j<i$ and $k\geq i$, the gold relations then give
\[
 u_{\lambda_j}=p^{i-1-j}a_{\lambda_j}q,\qquad
 a_{\lambda_k}=p^{k-i}u_{\lambda_k}x.
\]
Thus the relations reduce to $qx=p$, $p^{r+1}=0$, and $p^rx=0$, giving
\eqref{eq:odd-localized-core} with the normal forms of
Lemma~\ref{lem:odd-normal}.

The orientation classes $u_{\lambda_j}$ for $j\geq i$
restrict to units at $K$. For $j<i$, $S(\lambda_j)$ has
isotropy in $\mathcal F_{i-1}$, so the Euler cofiber sequence
for $S(\lambda_j)_+$ makes $a_{\lambda_j}$ a unit in $H^{\Phi_i}$.
Induction on the $G/K$-cells gives
\[
 F(K_m,H^{\Phi_i})\xrightarrow{\simeq}F(K_m,X).
\]
For $h\geq1$, dualizing the Euler
cofiber sequence for $K_{2h-1}=S(h\lambda_i)_+$ gives
\[
 F(K_{2h-1},X)\simeq
 \operatorname{cofib}\bigl(\Sigma^{-2h}X\xrightarrow{x^h}X\bigr).
\]
The transition maps on the fibers $\Sigma^{-2h}X$
are multiplication by $x$. Write $V=d+\sum_jc_j\lambda_j$.
The corresponding Laurent monomial is
\[
 w=\prod_{j<i}a_{\lambda_j}^{-c_j}
   \prod_{j\geq i}u_{\lambda_j}^{-c_j}
   \in\mathcal M_i^{(p)}.
\]
For $h\gg0$, each nonzero group $\pi_{V+2h}^G X$ is cyclic,
generated by $wq^A$, where
$2A=d+2h+2\sum_{j\geq i}c_j\geq0$.
If this integer is odd, the group is zero.
Multiplication by $x$ sends $wq^{A+1}$ to $pwq^A$.
Thus, for sufficiently large $h$, the composite of $r+1$
transition homomorphisms
\[
 \pi_{V+2(h+r+1)}^G X\xrightarrow{\ x^{r+1}\ }
 \pi_{V+2h}^G X
\]
is zero. The same argument applies to $V+1$. Hence
\[
 \varprojlim_h\pi_{V+2h}^G X=
 \varprojlim_h^{1}\pi_{V+2h}^G X=0,
\]
and similarly for $V+1$. The Milnor exact sequence gives
\[
 \pi_V^G X\xrightarrow{\cong}
 \pi_V^G\!\left(\holim_hF(K_{2h-1},X)\right).
\]
For $J=C_{p^m}$ with $i<m\leq n$, restriction sends
$\lambda_j$ to $\lambda_{p^m}(p^j)$ for $j<m$,
and to the trivial oriented plane for $j\geq m$.
Applying the same localization argument to $\res_J^G H$ gives
\[
 \bigoplus_{d\in\mathbb Z}\pi_d^J(X)\cong
 \mathbb Z[q,x]/(p^{m-i+1},qx-p,p^{m-i}x).
\]
Here $d$ denotes the degree $d\cdot1\in RO(J)$.
Fix $d\in\mathbb Z$. For large enough $h$, every composite of $m-i+1$
transition maps starting at $\pi_d^J(\Sigma^{-2h}X)$ is zero.
If $J\leq K$, then $E\mathcal F_i|_J$ is equivariantly contractible.
Hence the completion map is an isomorphism on every $\pi_d^J$, giving
the genuine $G$-equivalences
\[
 X\simeq\holim_hF(K_{2h-1},X)
 \simeq F((E\mathcal F_i)_+,X)\simeq H^{h_i}.
\]
\end{proof}

\begin{definition}\label{def:appendix-skeletal-filtration}
For $s\geq0$, the skeletal filtration is
\[
 F^s\pi_\star^G H^{h_i}
 =\ker\bigl(\pi_\star^G H^{h_i}\to
       \pi_\star^G F(K_{s-1},H^{\Phi_i})\bigr).
\]
\end{definition}

The Euler cofiber sequence computes this filtration.

\begin{proposition}[The skeletal filtration for odd primes]\label{prop:odd-skeletal-filtration}
For $\dim V$ even and $h\geq0$,
\begin{equation}\label{eq:odd-x-bockstein-tower}
 F^{2h}\pi_V^{G}H^{h_i}
 =\operatorname{im}\bigl(x^h:\pi_{V+2h}^{G}X\longrightarrow
                              \pi_V^{G}X\bigr),
 \qquad F^{2h+1}=F^{2h+2}.
\end{equation}
In every fixed $RO_{\mathrm{sel}}(G)$-degree this filtration is finite.
\end{proposition}
\begin{proof}
The cofiber sequence gives the formula for $F^{2h}$.
For $\dim V$ even, the cell in dimension $2h+1$ gives
\[
 \pi_V^G F(G/K_+\wedge S^{2h+1},H^{\Phi_i})
 \cong\pi_{V+2h+1}^K(H^{\Phi_i})=0.
\]
The corresponding long exact sequence gives
$F^{2h+1}=F^{2h+2}$. Finally, for fixed $V$ and sufficiently
large $h$, the image of $x^h$ is zero by
\eqref{eq:odd-normal-sum}. Thus the filtration is finite.
\end{proof}

\leadbreak
We now prove strong convergence.

\begin{proposition}
\label{prop:odd-convergence}
The spectral sequence \eqref{eq:odd-hfpss} converges strongly.  In each
fixed $RO_{\mathrm{sel}}(G)$-degree, the induced filtration is finite, exhaustive,
complete, and separated.
\end{proposition}

\begin{proof}
Lemma~\ref{lem:odd-weighted} shows that, for each $V$,
the inverse system $\{\pi_{V+2h}^G X\}_h$ is pro-zero:
every stage receives a zero transition homomorphism from
some later stage. Its limit and first derived inverse limit
therefore vanish, giving conditional convergence
\cite[Definition~5.10]{Boardman}.
Each group $E_2^{V,s}$ is finite: $V$ determines the
Laurent monomial and the power of $q$.
For each fixed $(V,s)$, the cycle subgroups of $E_2^{V,s}$
form a decreasing sequence and therefore stabilize.
Their first derived inverse limit is zero.
Since $E_2^{V,s}=0$ for $s<0$,
\cite[Theorem~7.1]{Boardman} gives strong convergence.
Equation~\eqref{eq:odd-x-bockstein-tower} gives the finite,
exhaustive, complete, separated filtration.
\end{proof}

Here $F^s$ denotes the skeletal filtration, and
$\operatorname{gr}^s\pi_D^G H^{h_i} =F^s\pi_D^G H^{h_i}/F^{s+1}\pi_D^G H^{h_i} \cong E_\infty^{D+s,s}$.
For $x\in F^s$, write $c\Longrightarrow x$ or $x\longleftrightarrow c$
when $\operatorname{gr}^s(x)=c\in E_\infty^{D+s,s}$.

We need one more lemma.

\begin{lemma}\label{lem:odd-even-cell}
In \eqref{eq:odd-hfpss}, for $m\geq1$,
\[
 \operatorname{gr}^{2m}(x^m)=t^m\ne0.
\]
\end{lemma}

\begin{proof}
The Euler cofiber sequence for $K_{2m-1}=S(m\lambda_i)_+$
gives $x^m\in F^{2m}$. Since $t=c_1(\lambda_i)\bmod p$,
the mod-$p$ Euler class of the associated bundle
$EW\times_W m\lambda_i\to BW$ is $t^m$.
Lemma~\ref{lem:odd-weighted} gives
\[
 F^{2m}\pi_{-2m}^G H^{h_i}=(\mathbb Z/p^r)\{x^m\},\qquad
 F^{2m+1}\pi_{-2m}^G H^{h_i}=p(\mathbb Z/p^r)\{x^m\}.
\]
Thus $t^m$ survives and detects $x^m$.
\end{proof}

We can now compute the differentials.

\begin{theorem}[Differentials]\label{thm:odd-differentials}
The only nonzero differential on the algebra generators in
\eqref{eq:odd-hfpss} is
\begin{equation*}
 d_{2r+1}(s)=q^rt^{r+1}.
\end{equation*}
All other nonzero differentials follow from
$L_i^{(p)}[q,t]$-linearity and the Leibniz rule.
\end{theorem}
\begin{proof}
Multiplication by $q$ and by the Laurent generators commutes
with the maps in the tower $F(K_m,H^{\Phi_i})$.
Lemma~\ref{lem:odd-even-cell} gives $\operatorname{gr}^2(x)=t$.
Thus all outgoing differentials vanish on $L_i^{(p)}[q,t]$.
Since $\pi_d^K(H^{\Phi_i})=0$ for odd integers $d$,
the first nonzero differential on $s$ has the following form
for some $m\geq1$ and $c\in\mathbb F_p^\times$:
\[
 d_{2m+1}(s)=c\,q^mt^{m+1}.
\]
The only monomial in bidegree $(2m,2m+2)$ is $q^mt^{m+1}$.
Even differentials have zero target, and $s$ has no incoming differential
because its filtration is one. Strong convergence and
$\pi_{-1}^G X=0$ therefore force a nonzero outgoing differential.

By \eqref{eq:odd-localized-core} and
\eqref{eq:odd-x-bockstein-tower},
\[
 \pi_0^{G}H^{h_i}=\mathbb Z/p^{r+1},\qquad
 F^{2h}\pi_0^{G}H^{h_i}=p^h(\mathbb Z/p^{r+1}).
\]
Thus the associated graded group of $\pi_0^G H^{h_i}$
has $r+1$ nonzero graded pieces, represented by
$1,qt,\ldots,q^rt^r$.
The differential $d_{2m+1}(s)=c q^mt^{m+1}$ is injective on
$L_i^{(p)}[q,t]s$ and leaves $m+1$ nonzero graded pieces
in total degree zero. Hence $m=r$. Choose $s$ so that $c=1$.
After this differential, $L_i^{(p)}[q,t]s$ contributes zero
to the next page, and every remaining class is permanent.
\end{proof}

\begin{remark}
No element of $L_i^{(p)}[q,t]$ supports a nonzero
differential, although multiples of $q^rt^{r+1}$ are hit by the displayed
differential.
\end{remark}

\subsection{The extensions for $p$ odd}

We now compute the associated graded ring.

\begin{proposition}\label{prop:odd-Einf}
For \eqref{eq:odd-hfpss},
\begin{equation}
 E_\infty\cong
 \frac{L_i^{(p)}[q,t]}{(q^rt^{r+1})}.                          \label{eq:odd-Einf}
\end{equation}
\end{proposition}

\begin{proof}
Take homology using the formula in Theorem~\ref{thm:odd-differentials}.
\end{proof}

\leadbreak
We have the following proposition.

\begin{proposition}\label{prop:odd-detections}
The filtration-zero detections in \eqref{eq:odd-hfpss} are
\[
 \overline a_{\lambda_j}^{\pm1}\Longrightarrow a_{\lambda_j}^{\pm1}
 \quad(j<i),\qquad
 \overline u_{\lambda_j}^{\pm1}\Longrightarrow u_{\lambda_j}^{\pm1}
 \quad(j\ge i).
\]
Moreover,
\begin{equation*}
 t\Longrightarrow x:=\frac{a_{\lambda_i}}{u_{\lambda_i}},
 \qquad
 q\Longrightarrow
 \frac{u_{\lambda_{i-1}}}{a_{\lambda_{i-1}}}.                 
\end{equation*}
\end{proposition}

\begin{proof}
The edge homomorphism
\[
 \pi_V^G H^{h_i}\longrightarrow
 H^0\!\left(W;\pi_V^K H^{\Phi_i}\right)
 \subseteq\pi_V^K H^{\Phi_i}
\]
is induced by $K_0=(G/K)_+\hookrightarrow(E\mathcal F_i)_+$.
The images of the Euler and orientation classes under
$H^{\Phi_i}\to H^{h_i}$ map under this edge homomorphism
to their restrictions at $K$.
The classes $a_{\lambda_j}$ are invertible for $j<i$
by Euler localization, and $u_{\lambda_j}$ are invertible
for $j\geq i$. The definition of $q$ and
Lemma~\ref{lem:odd-even-cell} give the remaining detections.
\end{proof}

For $w\in\mathcal M_i^{(p)}$, let $\overline w\in L_i^{(p)}$
be its image under the map $\Lambda_i^{(p)}\to L_i^{(p)}$
that reduces integers modulo $p$ and sends
$a_{\lambda_j}^{\pm1}$ to $\overline a_{\lambda_j}^{\pm1}$
for $j<i$, and $u_{\lambda_j}^{\pm1}$ to
$\overline u_{\lambda_j}^{\pm1}$ for $j\geq i$.
The corresponding associated-graded classes are given below,
in filtrations $2m$ and $2(B+m)$, respectively.
\begin{equation*}
 \begin{aligned}
 p^mwq^A&\longleftrightarrow\overline wq^{A+m}t^m,
 &&A\geq0,\quad0\leq m\leq r,\\
 p^mwx^B&\longleftrightarrow\overline wq^mt^{B+m},
 &&B\geq1,\quad0\leq m<r.
 \end{aligned}                                                
\end{equation*}

We now determine the hidden extensions. For $h\geq0$, give $B_{r,i}^{(p)}$
the filtration
\[
 F^{2h}B_{r,i}^{(p)}=x^hB_{r,i}^{(p)},\qquad
 F^{2h+1}B_{r,i}^{(p)}=F^{2h+2}B_{r,i}^{(p)}.
\]

\begin{theorem}
\label{thm:odd-extensions}
The assignments
\[
\begin{aligned}
 a_{\lambda_j}^{\pm1}&\longmapsto a_{\lambda_j}^{\pm1}
 &&(j<i),\\
 u_{\lambda_j}^{\pm1}&\longmapsto u_{\lambda_j}^{\pm1}
 &&(j\geq i),\\
 q&\longmapsto u_{\lambda_{i-1}}/a_{\lambda_{i-1}},\\
 x&\longmapsto a_{\lambda_i}/u_{\lambda_i}
\end{aligned}
\]
define an isomorphism of filtered rings
\[
 B_{r,i}^{(p)}\xrightarrow{\cong}\pi_\star^{G}H^{h_i}.
\]
As an algebra over $\Lambda_i^{(p)}$, its defining relations
can be written as
\begin{equation*}
 qx=p,\qquad p^rx=0.
\end{equation*}
\end{theorem}
\begin{proof}
The map is the isomorphism \eqref{eq:odd-localized-core}
followed by the completion equivalence of
Lemma~\ref{lem:odd-weighted}.
The relation $qx=p$ is the gold relation
$a_{\lambda_i}u_{\lambda_{i-1}}
=p\,a_{\lambda_{i-1}}u_{\lambda_i}$ after division by
$a_{\lambda_{i-1}}u_{\lambda_i}$.
The relation $p^rx=0$ follows from
$p^ra_{\lambda_i}=0$; multiplying by $q$ gives $p^{r+1}=0$.

The normal forms of Lemma~\ref{lem:odd-normal} and
\eqref{eq:odd-x-bockstein-tower} give the detections of
$p^mwq^A$ and $p^mwx^B$ displayed above. Hence
\[
 \operatorname{gr}B_{r,i}^{(p)}
 =\frac{(\Lambda_i^{(p)}/p)[q,t]}{(q^rt^{r+1})}.
\]
\end{proof}

This gives the additive extensions.

\begin{corollary}\label{cor:odd-additive-extensions}
The additive extensions from \eqref{eq:odd-Einf} to
$\pi_\star^G H^{h_i}$ are determined by
\begin{align*}
 p(wq^A)&=wq^{A+1}x &&(A\geq0),\\
 p(wx^B)&=wqx^{B+1} &&(B\geq1),
\end{align*}
for $w\in\mathcal M_i^{(p)}$.
\end{corollary}
\begin{proof}
Multiplying $qx=p$ by $wq^A$ or $wx^B$ gives these formulas.
Together with the orders in \eqref{eq:odd-normal-sum},
they determine the cyclic summands and their filtrations.
\end{proof}

\begin{proof}[Proof of Theorem~\ref{thm:odd-main}]
For $j<i$ and $k\geq i$, the adjacent gold relations give
\begin{align*}
 u_{\lambda_j}
   &=p^{i-1-j}\frac{a_{\lambda_j}}{a_{\lambda_{i-1}}}
      u_{\lambda_{i-1}},                                     \\
 a_{\lambda_k}
   &=p^{k-i}\frac{u_{\lambda_k}}{u_{\lambda_i}}
      a_{\lambda_i}.                                         
\end{align*}
With $q=u_{\lambda_{i-1}}/a_{\lambda_{i-1}}$ and
$x=a_{\lambda_i}/u_{\lambda_i}$, these identities recover every generator
in \eqref{eq:odd-main-ring} from \eqref{eq:odd-normalized}.
After substitution, every gold and Euler relation follows from
$qx=p$, $p^{r+1}=0$ and $p^rx=0$.
Theorem~\ref{thm:odd-extensions} gives \eqref{eq:odd-main-ring}.
\end{proof}

This gives the following corollary.

\begin{corollary}[Ordinary HFPSS for $p$ odd]\label{cor:odd-ordinary}
For $G=C_{p^n}$ with $p$ odd,
\begin{equation*}
 \pi_\star^{G}H^{h_0}
 \cong
 \frac{\mathbb{Z}[
  a_{\lambda_0},\ldots,a_{\lambda_{n-1}},
  u_{\lambda_0}^{\pm1},\ldots,u_{\lambda_{n-1}}^{\pm1}]}
 {(p^{n-j}a_{\lambda_j}\ (0\leq j\leq n-1),\
    \ \text{all gold relations})}.                            
\end{equation*}
\end{corollary}

\begin{proof}
For every prime $p$, choose $EG$ with one $G/e$-cell in each
dimension and put $K_m=((EG)_+)^{(m)}$. Since
\[
 \pi_V^G F(G/e_+\wedge S^m,H)
 =\pi_{\dim V+m}(H\mathbb Z),
\]
the map $\pi_V^G F(K_m,H)\to\pi_V^G F(K_{m-1},H)$ is an
isomorphism for $m>1-\dim V$. The same holds in degree $V+1$.
Thus $\lim_m^1\pi_{V+1}^G F(K_m,H)=0$ and
$\pi_V^G\operatorname{fib}(H^{h_0}\to F(K_m,H))=0$ for $m\gg0$.
The ordinary HFPSS converges strongly with a finite skeletal
filtration in each degree.

For odd $p$, \cite[Proposition~6.1]{Zeng} gives
\[
 \pi_\star^G H^{h_0}
 \cong
 \mathbb Z[a_{\lambda_0},u_{\lambda_0}^{\pm1},\ldots,
 u_{\lambda_{n-1}}^{\pm1}]/(p^n a_{\lambda_0}).
\]
The ordinary HFPSS has
$E_2^{V,s}=H^s(G;\pi_V^eH)$ and collapses at $E_2$.
Here $H^*(G;\mathbb Z)=\mathbb Z[t]/(p^nt)$, $|t|=2$,
and $t$ detects $a_{\lambda_0}/u_{\lambda_0}$.
Since $E_2^{V+s,s}\ne0$ implies $s=-\dim V$,
there are no hidden extensions.
For $0<j<n$, the gold relations give
\[
 a_{\lambda_j}
 =p^j a_{\lambda_0}u_{\lambda_j}u_{\lambda_0}^{-1},
\]
which gives the stated presentation.
\end{proof}

\subsection{Localizations and the skeletal tower for $p=2$}

Let $p=2$, with the sign representation $\alpha$ and oriented planes
from Section~\ref{sec:local}. Thus
$\lambda_{n-1}=2\alpha$, $a_{\lambda_{n-1}}=a_\alpha^2$, $u_{\lambda_{n-1}}=u_{2\alpha}$.
The family $\lambda_0,\ldots,\lambda_{n-2}$ is empty for $n=1$.
The degrees are
$|a_\alpha|=-\alpha$, $|a_{\lambda_j}|=-\lambda_j$, $|u_{2\alpha}|=2-2\alpha$, $|u_{\lambda_j}|=2-\lambda_j$.
The Euler classes $a_\alpha$ and $a_{\lambda_j}$ satisfy
\begin{equation}
 2a_\alpha=0,
 \qquad
 2^{n-j}a_{\lambda_j}=0\quad(0\leq j\leq n-2).               \label{eq:two-euler}
\end{equation}
For $0\leq j<k\leq n-1$, the gold relations are
\begin{align*}
 a_{\lambda_k}u_{\lambda_j}
   &=2^{k-j}a_{\lambda_j}u_{\lambda_k}.
\end{align*}
For $k=n-1$ and $0\leq j\leq n-2$, these relations become
\begin{equation}
 a_\alpha^2u_{\lambda_j}
   =2^{n-1-j}a_{\lambda_j}u_{2\alpha}.                        \label{eq:two-top-gold}
\end{equation}
We next consider the case $p=2$.

\begin{theorem}[The ring $\pi_\star^G(H^{h_i})$ for $p=2$]\label{thm:two-main}
For $p=2$ and $1\leq i<n$, we have an $RO_{\mathrm{sel}}(G)$-graded ring isomorphism
\begin{equation}
\pi_\star^{G}H^{h_i}\cong
\frac{
 \begin{gathered}
  (\mathbb{Z}/2^{n-i+1})[
   a_\alpha,
   \{a_{\lambda_j}\}_{i\leq j\leq n-2},
   \{a_{\lambda_j}^{\pm1}\}_{0\leq j<i},\\
   u_{2\alpha}^{\pm1},
   \{u_{\lambda_j}^{\pm1}\}_{i\leq j\leq n-2},
   \{u_{\lambda_j}\}_{0\leq j<i}]
 \end{gathered}
}{\left(\begin{gathered}
 2a_\alpha,\quad 2^{n-j}a_{\lambda_j}\quad(i\leq j\leq n-2),\\
 a_{\lambda_k}u_{\lambda_j}-2^{k-j}a_{\lambda_j}u_{\lambda_k}
 \quad(0\leq j<k\leq n-1)
\end{gathered}\right)}.                                      \label{eq:two-main-ring}
\end{equation}
\end{theorem}

\subsection{The generalized HFPSS and its differentials for $p=2$}

Assume $n\geq2$ and fix $1\leq i\leq n-1$. Put $r=n-i$,
$K=C_{2^i}$ and $W=G/K=C_{2^r}$, with $\gamma$ the image of the
chosen generator in $W$. Put
$K_{-1}=*$ and $K_m=((E\mathcal F_i)_+)^{(m)}$ for $m\geq0$.
An overline denotes the image under
$\pi_\star^G(H)\to\pi_\star^G(H^{\Phi_i})\xrightarrow{\res_K^G}
\pi_\star^K(H^{\Phi_i})$.
Since $\alpha|_K\cong1$ as an oriented representation, the resulting
sphere isomorphism gives a unit in $\pi_{\alpha-1}^K(H)$; let
$e_\alpha\in\pi_{\alpha-1}^K(H^{\Phi_i})$ be its image.
Define
\begin{align*}
 L_i^{(2)}&:=\mathbb{F}_2[
  \overline{a}_{\lambda_0}^{\pm1},\ldots,
  \overline{a}_{\lambda_{i-1}}^{\pm1},
  \overline{u}_{\lambda_i}^{\pm1},\ldots,
  \overline{u}_{\lambda_{n-2}}^{\pm1}],\\
 q&:=\frac{\overline{u}_{\lambda_{i-1}}}
          {\overline{a}_{\lambda_{i-1}}},
 \qquad |q|=(2,0),
 \qquad |e_\alpha|=(\alpha-1,0).
\end{align*}

For $r=1$, put
\[
 H^*(W;\mathbb{F}_2)=\mathbb{F}_2[b],\qquad |b|=(0,1),\qquad t=b^2.
\]
For $r\geq2$,
\[
 H^*(W;\mathbb{F}_2)=\mathbb{F}_2[t]\otimes\Lambda(s),
 \qquad |s|=(0,1),\quad |t|=(0,2).
\]

We first compute the $E_2$-term.

\begin{proposition}[The $E_2$-term for $p=2$]\label{prop:two-E2}
For $p=2$ and $1\leq i<n$, the tower $F(K_m,H^{\Phi_i})$ above
gives a multiplicative spectral sequence
\begin{equation}
 E_2^{V,s}=H^s\!\left(W;\pi_V^{K}(H^{\Phi_i})\right)
 \Longrightarrow\pi_{V-s}^{G}H^{h_i},
 \qquad |d_k|=(k-1,k),                                        \label{eq:two-hfpss}
\end{equation}
with
\[
 \pi_\star^{K}(H^{\Phi_i})\cong L_i^{(2)}[q,e_\alpha^{\pm1}].
\]
Consequently
\begin{equation*}
 E_2=\begin{cases}
  L_i^{(2)}[q,e_\alpha^{\pm1},b],&r=1,\\[1mm]
  L_i^{(2)}[q,e_\alpha^{\pm1},t]\otimes\Lambda(s),&r\geq2.
 \end{cases}                                                   
\end{equation*}
\end{proposition}
\begin{proof}
Let $\mathcal P_K$ be the family of proper subgroups of $K$,
and let $\alpha_K$ be the sign representation of $K$.
Since $\mathcal F_{i-1}\cap K=\mathcal P_K$,
\[
 \res_K^G H^{\Phi_i}
 \simeq\widetilde E\mathcal P_K\wedge\res_K^G H.
\]
By \cite[Proposition~3.18]{HHRNonexistence},
$a_{\alpha_K}$ is invertible and
\[
 \bigoplus_{d\in\mathbb Z}\pi_d^K(H^{\Phi_i})
 \cong\pi_*(\Phi^K\res_K^G H)
 \cong\mathbb F_2[q].
\]
Here $\lambda_{i-1}|_K=2\alpha_K$ and
$q=u_{2\alpha_K}a_{\alpha_K}^{-2}
=\overline u_{\lambda_{i-1}}/\overline a_{\lambda_{i-1}}$,
with $|q|=2$.
The gold relations give
$\overline u_{\lambda_j}=0$ for $0\leq j<i-1$.

For $j<i$, $(\lambda_j|_K)^K=0$, so $S(\lambda_j|_K)$ has only proper
$K$-isotropy. Smashing its Euler cofiber sequence with
$\widetilde E\mathcal P_K\wedge\res_K^G H$ makes
$\overline a_{\lambda_j}$ invertible. For $j\geq i$,
$\lambda_j|_K\cong2$, so $\overline u_{\lambda_j}$ is invertible.
The oriented isomorphism $\alpha|_K\cong1$ also gives
\[
 e_\alpha^2=\overline u_{2\alpha}^{-1}.
\]
These units reduce every selected degree to an integer degree.
The ring $\mathbb F_2[q]$ therefore gives
\[
 \pi_\star^{K}(H^{\Phi_i})\cong L_i^{(2)}[q,e_\alpha^{\pm1}].
\]

The classes $\overline a_{\lambda_j}$ for $j<i$,
$\overline u_{\lambda_j}$ for $i\leq j\leq n-2$, and $q$
are restrictions of classes in $\pi_\star^G(H^{\Phi_i})$.
Hence $W$ fixes these classes and the inverses occurring in $L_i^{(2)}$.
On $e_\alpha$,
$\gamma e_\alpha=-e_\alpha$,
because $\gamma$ reverses the orientation of $\alpha$.
Since $2=0$ in $L_i^{(2)}[q,e_\alpha^{\pm1}]$,
the $W$-action on this ring is trivial.
The cellular filtration of
$E\mathcal F_i\simeq S(\infty\lambda_i)$ therefore gives
\eqref{eq:two-hfpss}; the periodic resolution of $W$ gives
the last formula in the proposition.
\end{proof}

Define
\[
 \Lambda_i^{(2)}:=\mathbb{Z}[
  a_{\lambda_0}^{\pm1},\ldots,a_{\lambda_{i-1}}^{\pm1},
  u_{\lambda_i}^{\pm1},\ldots,u_{\lambda_{n-2}}^{\pm1},
  u_{2\alpha}^{\pm1}]
\]
and, with formal variables $|q|=2$, $|x|=-2$, $|z|=\alpha-2$,
\begin{equation}
 B_{r,i}^{(2)}:=
 \frac{\Lambda_i^{(2)}[q,x,z]}
 {(2^{r+1},\;qx-2,\;2^rx,\;2z,\;
   z^2-2^{r-1}xu_{2\alpha}^{-1})}.                            \label{eq:two-normalized}
\end{equation}
We have the following lemma.

\begin{lemma}\label{lem:two-normal}
Let $\mathcal M_i^{(2)}$ be the Laurent-monomial basis of
$\Lambda_i^{(2)}$. For $B_{r,i}^{(2)}$ in \eqref{eq:two-normalized}, we have
\begin{equation*}
\begin{split}
 B_{r,i}^{(2)}\cong\bigoplus_{w\in\mathcal M_i^{(2)}}\Bigg(&
 \bigoplus_{A\geq0}\mathbb{Z}/2^{r+1}\{wq^A\}
 \oplus\bigoplus_{B\geq1}\mathbb{Z}/2^r\{wx^B\}\\
 &{}\oplus\bigoplus_{A\geq0}\mathbb{F}_2\{wzq^A\}
 \oplus\bigoplus_{B\geq1}\mathbb{F}_2\{wzx^B\}\Bigg).
\end{split}
\end{equation*}
\end{lemma}

\begin{proof}
Put
\[
 E_r:=\frac{\Lambda_i^{(2)}[q,x]}
 {(2^{r+1},qx-2,2^rx)}.
\]
The proof of Lemma~\ref{lem:odd-normal} applies with $p=2$ and gives
the additive summands of $E_r$ generated by $wq^A$ and $wx^B$.
Since
\[
 E_r/(2)\cong(\Lambda_i^{(2)}/2)[q,x]/(qx),\qquad
 2\bigl(2^{r-1}xu_{2\alpha}^{-1}\bigr)=0,
\]
for $a,a'\in E_r$ and $b,b'\in E_r/(2)$, the multiplication
\[
 (a+bz)(a'+b'z)
 =aa'+2^{r-1}xu_{2\alpha}^{-1}bb'+(ab'+a'b)z
\]
on $E_r\oplus(E_r/(2))z$ is independent of the lifts of $b,b'$.
It is associative, commutative, and unital, since
$2^{r-1}xu_{2\alpha}^{-1}$ is central and annihilated by
$\ker(E_r\to E_r/(2))=(2)$.
Reduction by $2z=0$ and $z^2=2^{r-1}xu_{2\alpha}^{-1}$ gives
mutually inverse maps
\[
 B_{r,i}^{(2)}\cong E_r\oplus(E_r/(2))z.
\]
The additive decomposition of $E_r$ and the monomial
$\mathbb F_2$-basis of $E_r/(2)$ give the formula in the lemma.
\end{proof}

In the setting of Proposition~\ref{prop:two-E2}, put
\[
 X=H^{\Phi_i}[u_{2\alpha}^{-1},
                  \{u_{\lambda_j}^{-1}\}_{i\leq j\leq n-2}],
 \qquad
 q=\frac{u_{\lambda_{i-1}}}{a_{\lambda_{i-1}}},\quad
 x=\frac{a_{\lambda_i}}{u_{\lambda_i}},\quad
 z=\frac{a_\alpha}{u_{2\alpha}}.
\]

We also need the following lemma.

\begin{lemma}\label{lem:two-weighted}
Then
\[
 \pi_\star^{G}X\cong B_{r,i}^{(2)},
\]
and, for $h\geq0$,
\begin{equation*}
 F(K_{2h-1},H^{\Phi_i})\simeq
 \operatorname{cofib}\bigl(\Sigma^{-2h}X\xrightarrow{x^h}X\bigr),
 \qquad X\xrightarrow{\simeq}H^{h_i}.                         
\end{equation*}
The skeletal filtration satisfies
\begin{equation}
 F^{2h}\pi_V^{G}H^{h_i}
   =x^h\pi_{V+2h}^{G}X.                                    \label{eq:two-even-filtration}
\end{equation}
\end{lemma}
\begin{proof}
By \cite[Theorem~5.4]{SikoraYan}, the subring of $\pi_\star^G H$
in the positive cone is generated by the Euler and orientation classes,
subject to
\eqref{eq:two-euler} and the gold relations.  Inverting
$a_{\lambda_j}$ for $j<i$, $u_{\lambda_j}$ for $i\leq j\leq n-2$, and
$u_{2\alpha}$ moves every selected degree into this cone. For
$V=d+c_\alpha\alpha+\sum_{j=0}^{n-2}c_j\lambda_j$, the numerator
degree after using $k_j$ copies of the indicated $\lambda_j$-class
and $h$ copies of $u_{2\alpha}$ is
$d+2h+2\sum_{j\geq i}k_j +(c_\alpha-2h)\alpha+\sum_j(c_j-k_j)\lambda_j$.
The tuples of localization exponents satisfying
$k_j\geq\max(c_j,0)$ and
$h\geq\max(0,\lceil c_\alpha/2\rceil)$ form a cofinal set.
For these tuples, $c_\alpha-2h\leq0$ and $c_j-k_j\leq0$ for every $j$.
Equivariant homotopy groups commute with filtered colimits, so the
localized positive-cone presentation computes $\pi_\star^{G}X$.
Eliminating the remaining Euler and orientation generators by the gold
relations gives
\[
 qx=2,\qquad 2^{r+1}=0,\qquad 2^rx=0,\qquad 2z=0,
 \qquad z^2=2^{r-1}xu_{2\alpha}^{-1}.
\]
These are the relations of \eqref{eq:two-normalized}, with normal forms
given by Lemma~\ref{lem:two-normal}.

Every cell of $K_m$ has isotropy $K=C_{2^i}$, and the inverted
orientation classes restrict to units at $K$. Hence
\[
 F(K_m,H^{\Phi_i})\xrightarrow{\simeq}F(K_m,X).
\]
Dualizing $S(h\lambda_i)_+\to S^0\to S^{h\lambda_i}$ and using the
unit $u_{\lambda_i}$ on $X$ gives the first equivalence in the second formula of the lemma. The fiber transition map is
\[
 \Sigma^{-2h-2}X\xrightarrow{x}\Sigma^{-2h}X.
\]
For fixed $V$ and $h\gg0$, the normal forms for $\pi_{V+2h}^{G}X$
are $wq^A$ and $wzq^A$, with $A\geq1$. Since
\[
 xq^A=2q^{A-1},\qquad xzq^A=0\quad(A\geq1),
\]
for $h\gg0$, every composite of $r+1$ transition maps is zero.
Thus the fiber homotopy groups form a pro-zero tower in each degree,
and $\lim=\lim^1=0$.
For $J=C_{2^m}$ with $i<m<n$,
$\alpha|_J=1$, $\lambda_{m-1}|_J=2\alpha_J$, $\lambda_j|_J=2\;(j\geq m)$.
Here $\alpha_J$ is the sign representation of $J$.
Localizing the positive-cone presentation for $J$
at the restrictions of the classes inverted in $X$
gives the integer-graded subring
\[
 \bigoplus_{d\in\mathbb Z}\pi_d^J(\res_J^G X)=
 \mathbb Z[q,x]/(2^{m-i+1},qx-2,2^{m-i}x).
\]
The class $z$ restricts to zero since $\alpha|_J=1$.
In sufficiently high integer degrees, $xq^A=2q^{A-1}$.
Thus the fiber homotopy groups at $J$ also form a pro-zero tower.
For $J\subseteq K$, the space $E\mathcal F_i|_J$
is $J$-equivariantly contractible and
the inverted orientations are already units. The completion map
$X\to H^{h_i}$ therefore induces isomorphisms on $\pi_d^J$ for every
$J\leq G$ and $d\in\mathbb Z$, so $X\simeq H^{h_i}$.
The long exact sequence of the cofiber in the second formula of the lemma gives \eqref{eq:two-even-filtration}.
For fixed $V$, Lemma~\ref{lem:two-normal} gives
$x^h\pi_{V+2h}^G X=0$ for $h\gg0$, so the skeletal filtration is finite.
\end{proof}

We now prove strong convergence.

\begin{proposition}
\label{prop:two-convergence}
The spectral sequence \eqref{eq:two-hfpss} converges strongly, and its
filtration is finite, exhaustive, complete, and separated in each fixed
$RO_{\mathrm{sel}}(G)$-degree.
\end{proposition}

\begin{proof}
Lemma~\ref{lem:two-weighted} gives degreewise pro-zero fiber homotopy
groups and a finite skeletal filtration. Since each $E_2^{V,s}$ is finite
by Proposition~\ref{prop:two-E2}, the proof of
Proposition~\ref{prop:odd-convergence} applies.
\end{proof}

Put
\[
 y=\begin{cases}
  e_\alpha b,&r=1,\\
  e_\alpha s,&r\geq2.
 \end{cases}
\]

We have the following proposition.

\begin{proposition}
\label{prop:two-detections}
In filtration zero of \eqref{eq:two-hfpss}, the detections in
$\pi_\star^G X\cong\pi_\star^G H^{h_i}$ are
\[
 \begin{aligned}
 \overline a_{\lambda_j}^{\pm1}&\Longrightarrow a_{\lambda_j}^{\pm1}
 &&(j<i),\\
 \overline u_{\lambda_j}^{\pm1}&\Longrightarrow u_{\lambda_j}^{\pm1}
 &&(i\leq j\leq n-2),\\
 e_\alpha^{\pm2}&\Longrightarrow u_{2\alpha}^{\mp1}.
 \end{aligned}
\]
Moreover,
\begin{equation*}
 t\Longrightarrow x=\frac{a_{\lambda_i}}{u_{\lambda_i}},
 \qquad y\Longrightarrow z=u_{2\alpha}^{-1}a_\alpha.         
\end{equation*}
\end{proposition}

\begin{proof}
Lemma~\ref{lem:two-weighted} lifts the generators of $L_i^{(2)}$, together
with $e_\alpha^{\pm2}$, to the corresponding Euler and orientation
units. The Euler cofiber sequence shows that $t$ detects $x$:
in degree $-2$,
$F^2=\mathbb Z/2^r\{x\}$ and $F^4=2\mathbb Z/2^r\{x\}$,
while the odd filtration quotient between them is zero because
$E_2^{1,3}=0$. Since $\alpha|_K\cong1$, the restriction of $z$ lies in
$\pi_{-1}^{K}\res_K^G(H^{\Phi_i})=0$; hence $z\in F^1$.
Lemmas~\ref{lem:two-normal} and \ref{lem:two-weighted} give
\[
 \pi_{\alpha-2}^{G}X=\mathbb F_2\{z\},\qquad
 F^2\pi_{\alpha-2}^{G}X
 =x\pi_\alpha^GX=x\mathbb F_2\{qz\}=0,
\]
since $xqz=2z=0$. Strong convergence
(Proposition~\ref{prop:two-convergence}) therefore gives
\[
 E_\infty^{\alpha-1,1}
 \cong F^1\pi_{\alpha-2}^GX/F^2\pi_{\alpha-2}^GX
 =\mathbb F_2\{[z]\}.
\]
By Proposition~\ref{prop:two-E2},
$\pi_{\alpha-1}^K(H^{\Phi_i})=\mathbb F_2\{e_\alpha\}$, so
$E_2^{\alpha-1,1}=\mathbb F_2\{y\}$.
No $d_k$ enters filtration one, since $1-k<0$ for $k\ge2$.
Hence $y$ survives and detects $z$.
\end{proof}

We can now give all the differentials at $p=2$.

For $C_4$, Yan gives the two $d_3$-differentials below
in \cite[Proposition~3.5]{Yan}.

\begin{theorem}[Differentials for $p=2$]
\label{thm:two-differentials}
In \eqref{eq:two-hfpss}, if $r=1$, the differentials are determined by
\begin{equation*}
 d_3(b)=qb^4=qt^2,
 \qquad d_3(e_\alpha)=qe_\alpha b^3=qe_\alpha bt.            
\end{equation*}
If $r\geq2$, the differentials are determined by
\begin{equation*}
 d_3(e_\alpha)=qe_\alpha st,
 \qquad d_{2r+1}(s)=q^rt^{r+1}.                               
\end{equation*}
All other nonzero differentials follow from linearity and the Leibniz rule.
\end{theorem}
\begin{proof}
The class $q\in\pi_2^G X$ and
Proposition~\ref{prop:two-detections} show that
$q,t,e_\alpha^2$ and all elements of $L_i^{(2)}$ are permanent cycles.
For $V=d+c_\alpha\alpha+\sum_{j=0}^{n-2}c_j\lambda_j$,
$E_2^{V,s}=0$ unless $d\equiv c_\alpha\pmod2$.
Hence all even differentials vanish.

First assume $r\geq2$. By Proposition~\ref{prop:two-detections},
\[
 y=e_\alpha s\Longrightarrow z=a_\alpha u_{2\alpha}^{-1},\qquad
 F^2\pi_{\alpha-2}^{G}X=0.
\]
Thus $qty$, of total degree $\alpha-2$ and filtration $3$, must be hit.
The classes $q,t,y$ are permanent cycles, so $qty$ cannot support a differential.
An incoming $d_k$ has source filtration $3-k$; nonnegativity and parity
leave only $k=3$, with source $\mathbb F_2\{e_\alpha\}$. Hence
\[
 d_3(e_\alpha)=qty=qe_\alpha st.
\]
Also
\[
 F^4\pi_{-2}^{G}X/F^6\pi_{-2}^{G}X
   =2\mathbb Z/2^r\{x\}/4\mathbb Z/2^r\{x\}
   \cong\mathbb F_2,
\]
so $qt^2$ survives and $d_3(s)=0$.
Over $L_i^{(2)}[q,t,e_\alpha^{\pm2}]$, the $d_3$-homology has
module generators $1,s,y$, with
\[
 qty=0,\qquad s^2=sy=y^2=0.
\]
Indeed, every power of $e_\alpha$ is uniquely an even power times $1$
or $e_\alpha$. For $f\in L_i^{(2)}[q,t,e_\alpha^{\pm2}]$,
\[
 d_3(f)=d_3(sf)=d_3(e_\alpha sf)=0,\qquad
 d_3(e_\alpha f)=qtyf.
\]
The last map is injective, with image $(qt)y$.
By bidegree, the first possible nonzero differential on $s$ has the form
\[
 d_{2k+1}(s)=q^kt^{k+1}\quad(k\geq1).
\]
The skeletal filtration is
\[
 \pi_{-2}^{G}X=\mathbb Z/2^r\{x\},\qquad
 F^{2h}\pi_{-2}^{G}X=
 \begin{cases}
  2^{h-1}\mathbb Z/2^r\{x\},&1\leq h\leq r,\\
  0,&h\geq r+1.
 \end{cases}
\]
For $h\geq0$, $E_2^{2h-1,2h+1}=0$, so
$\operatorname{gr}^{2h+1}\pi_{-2}^G X=0$.
Thus $q^{h-1}t^h$ survives for $1\leq h\leq r$, whereas
$q^rt^{r+1}$ must be hit. The possible incoming differentials on
$q^rt^{r+1}$ have sources
$s q^{r-k}t^{r-k}$ on pages $2k+1$, $1\leq k\leq r$.
A differential with $k<r$ would also kill $q^kt^{k+1}$, contradicting
the nonzero filtration quotient above. Hence
\[
 d_{2r+1}(s)=q^rt^{r+1}.
\]
This differential is injective on the summand
$L_i^{(2)}[q,t,e_\alpha^{\pm2}]\{s\}$.

When $r=1$, the same filtration calculation gives
\[
 F^2\pi_{\alpha-2}^G X=0,\qquad F^4\pi_{-2}^G X=0.
\]
Thus $qty$ and $qt^2$ must be hit, respectively from $e_\alpha$ and $b$
on page three. Hence
\[
 d_3(e_\alpha)=qe_\alpha bt,\qquad d_3(b)=qt^2.
\]
For $A,\ell\geq0$ and $k\in\mathbb Z$, the differential on monomials is
\[
 d_3(q^A e_\alpha^k b^\ell)
  =(k+\ell)q^{A+1}e_\alpha^k b^{\ell+3},
\]
with the coefficient read modulo two. Thus the cycles are the
$L_i^{(2)}$-linear combinations of the monomials with $k+\ell$ even,
and the boundaries give
$qt^2=qty=0$. The resulting homology ring is generated by
$q,e_\alpha^{\pm2},t,y$ with
$y^2=e_\alpha^2t$. The differential $d_3$ is $L_i^{(2)}$-linear.
After these differentials, the surviving algebra is generated by
$q,t,e_\alpha^{\pm2},L_i^{(2)}$ and $y$, all of which are permanent cycles.
Thus all further differentials vanish. Linearity and the Leibniz rule
give the remaining differentials.
\end{proof}

\begin{remark}
The classes $q,t,e_\alpha^2$ and all generators of $L_i^{(2)}$ are
permanent cycles.
\end{remark}

\subsection{The extensions for $p=2$}

We now give the associated graded ring.

\begin{proposition}
\label{prop:two-Einf}
For \eqref{eq:two-hfpss} with $r\geq2$,
\begin{equation*}
 E_\infty\cong
 \frac{L_i^{(2)}[q,t,y,e_\alpha^{\pm2}]}
 {(q^rt^{r+1},\;qty,\;y^2)}.                                  
\end{equation*}
For $r=1$,
\begin{equation*}
 E_\infty\cong
 \frac{L_i^{(2)}[q,t,y,e_\alpha^{\pm2}]}
 {(qt^2,\;qty,\;y^2-e_\alpha^2t)}.                           
\end{equation*}
\end{proposition}

\begin{proof}
Take homology under the differentials of
Theorem~\ref{thm:two-differentials}.
\end{proof}

The relations in the abutment are
\begin{align}
 qx&=2,                                                        \notag\\
 2^rx&=0,                                                      \notag\\
 2z&=0,                                                        \notag\\
 z^2&=2^{r-1}xu_{2\alpha}^{-1}
     =q^{r-1}x^ru_{2\alpha}^{-1}.                            \label{eq:two-hidden-square}
\end{align}
Multiplying $2a_\alpha=0$ by $u_{2\alpha}^{-1}$ gives
the third formula above.
For $r\geq2$, the first equality in \eqref{eq:two-hidden-square} follows
by dividing \eqref{eq:two-top-gold}, with $j=i$, by
$u_{\lambda_i}u_{2\alpha}^2$.
For $r=1$, $i=n-1$ and $\lambda_{n-1}=2\alpha$, so
\[
 x=a_\alpha^2u_{2\alpha}^{-1},\qquad
 z=a_\alpha u_{2\alpha}^{-1},\qquad
 z^2=xu_{2\alpha}^{-1}.
\]
In both cases, $qx=2$ gives the second expression in
\eqref{eq:two-hidden-square}.

For $w\in\mathcal M_i^{(2)}$, let
$\overline w\in L_i^{(2)}[e_\alpha^{\pm2}]$ be the filtration-zero image of $w$. The leading terms are
\begin{align*}
 2^mwq^A&\longleftrightarrow\overline wq^{A+m}t^m,
 &&A\geq0,\quad0\leq m\leq r,\\
 2^mwx^B&\longleftrightarrow\overline wq^mt^{B+m},
 &&B\geq1,\quad0\leq m<r,\\
 wzq^A&\longleftrightarrow\overline wq^Ay,
 &&A\geq0,\\
 wzx^B&\longleftrightarrow\overline wt^By,
 &&B\geq1.
\end{align*}
Equation~\eqref{eq:two-even-filtration} gives the first two identities.
The class $q^Ay\ne0$ detects $zq^A$ in filtration one, and
$F^2\pi_{\alpha-2+2A}^G X=x\mathbb F_2\{zq^{A+1}\}=0$.
Multiplicativity gives $zx^B\in F^{2B+1}$, whereas
\[
 F^{2B+2}\pi_{\alpha-2-2B}^G X
   =x^{B+1}\mathbb F_2\{zq\}=0.
\]
Since $zx^B\ne0$, its filtration is exactly $2B+1$.
Since $w$ and $w^{-1}$ are filtration-zero units, multiplication by $w$
multiplies each leading term by $\overline w$.
For $r\geq2$, the nonzero class $z^2$ has filtration $2r>2$.
Hence the product $y^2$ is zero in filtration two of the
associated graded ring.

We next determine the hidden extensions.

For $h\geq0$, equip $B_{r,i}^{(2)}$ with the filtration
\[
 F^{2h}B_{r,i}^{(2)}=x^hB_{r,i}^{(2)},\qquad
 F^{2h+1}B_{r,i}^{(2)}
 =x^h(x,z)B_{r,i}^{(2)}.
\]

\begin{theorem}[Hidden extensions for $p=2$]
\label{thm:two-extensions}
For \eqref{eq:two-hfpss}, the homomorphism
\[
 B_{r,i}^{(2)}\longrightarrow\pi_\star^G H^{h_i}
\]
sending the Laurent generators to the classes with the same names and
\[
 q\longmapsto\frac{u_{\lambda_{i-1}}}{a_{\lambda_{i-1}}},
 \qquad
 x\longmapsto\frac{a_{\lambda_i}}{u_{\lambda_i}},
 \qquad
 z\longmapsto\frac{a_\alpha}{u_{2\alpha}}
\]
is an isomorphism of filtered rings.
\end{theorem}
\begin{proof}
Lemma~\ref{lem:two-weighted} gives
$\pi_\star^G X\cong B_{r,i}^{(2)}$ and $X\simeq H^{h_i}$.
The preceding leading terms show that the filtrations agree.
\end{proof}

The hidden extensions are as follows.

\begin{corollary}\label{cor:two-extension-relations}
The relations
\begin{align*}
 2(wq^A)&=wq^{A+1}x &&(A\geq0),\\
 2(wx^B)&=wqx^{B+1} &&(B\geq1),                               
\end{align*}
for $w\in\mathcal M_i^{(2)}$, together with \eqref{eq:two-hidden-square} and
the Euler and gold relations, give all additive and multiplicative hidden extensions.
\end{corollary}
\begin{proof}
Multiplying $qx=2$ by $wq^A$ or $wx^B$ gives the formulas in the corollary.
The normal forms of Lemma~\ref{lem:two-normal}, together with
\eqref{eq:two-hidden-square}, determine all products and additive relations.
\end{proof}

\begin{proof}[Proof of Theorem~\ref{thm:two-main}]
For $j<i$ and $k\geq i$, the gold relations give
\begin{align*}
 u_{\lambda_j}
  &=2^{i-1-j}\frac{a_{\lambda_j}}{a_{\lambda_{i-1}}}
     u_{\lambda_{i-1}},\\
 a_{\lambda_k}
  &=2^{k-i}\frac{u_{\lambda_k}}{u_{\lambda_i}}
     a_{\lambda_i}.
\end{align*}
These identities and $a_\alpha=zu_{2\alpha}$ recover all generators
of \eqref{eq:two-main-ring}. Conversely, the gold relations reduce to
$qx=2$ and \eqref{eq:two-hidden-square}, while the Euler relations reduce
to $2^{r+1}=0$, $2^rx=0$, and $2z=0$.
Theorem~\ref{thm:two-extensions} and
Corollary~\ref{cor:two-extension-relations} give \eqref{eq:two-main-ring}.
\end{proof}

This gives the following corollary.

\begin{corollary}\label{cor:two-ordinary}
For $G=C_{2^n}$,
\begin{equation*}
 \pi_\star^{G}H^{h_0}
 \cong
 \frac{\mathbb{Z}[
  a_\alpha,\{a_{\lambda_j}\}_{0\leq j\leq n-2},
  u_{2\alpha}^{\pm1},
  \{u_{\lambda_j}^{\pm1}\}_{0\leq j\leq n-2}]}
 {(2a_\alpha,\ 2^{n-j}a_{\lambda_j}\ (0\leq j\leq n-2),
   \text{all gold relations})}.                               
\end{equation*}
For $n=1$, the indexed generator families are empty.
\end{corollary}

\begin{proof}
The ordinary HFPSS $E_2^{V,s}=H^s(G;\pi_V^eH)$ uses the coefficient modules
$\mathbb Z$ and $\mathbb Z_{\mathrm{sgn}}$ and collapses at $E_2$.
Here a generator of $G$ acts by $-1$ on $\mathbb Z_{\mathrm{sgn}}$.
The finite-stage stabilization
argument of Corollary~\ref{cor:odd-ordinary} gives strong convergence.
The corresponding group cohomology groups are
\[
 H^*(C_{2^n};\mathbb{Z})=\mathbb{Z}[t]/(2^nt),
 \qquad
 H^*(C_{2^n};\mathbb{Z}_{\mathrm{sgn}})=\mathbb{Z}/2[t]\{s\},
\]
with $|t|=2$ and $|s|=1$. The class $t$ detects
$a_{\lambda_0}u_{\lambda_0}^{-1}$; $s$ in internal
degree $1-\alpha$ and filtration one detects $a_\alpha$.
The orientation classes $u_{2\alpha}$ and $u_{\lambda_j}$ are invertible.

If $n=1$, then $\lambda_0=2\alpha$, so $t$ detects
$a_\alpha^2u_{2\alpha}^{-1}$. The displayed presentation reduces to
\[
 \frac{\mathbb{Z}[a_\alpha,u_{2\alpha}^{\pm1}]}{(2a_\alpha)}.
\]
Hence assume $n\geq2$.

Equation~\eqref{eq:two-top-gold} with $j=0$ gives
\begin{equation*}
 a_\alpha^2u_{\lambda_0}
  =2^{n-1}a_{\lambda_0}u_{2\alpha}.                        
\end{equation*}
For $1\leq j\leq n-2$, the remaining Euler classes are determined by
\[
 a_{\lambda_j}=2^j a_{\lambda_0}\frac{u_{\lambda_j}}{u_{\lambda_0}}.
\]
Every monomial reduces to a scalar multiple of a monomial in
$a_{\lambda_0}$, $u_{\lambda_0}^{\pm1},\ldots, u_{\lambda_{n-2}}^{\pm1}$, $u_{2\alpha}^{\pm1}$,
or of $a_\alpha$ times such a monomial. Products of two terms containing
$a_\alpha$ reduce using
$a_\alpha^2=2^{n-1}a_{\lambda_0}u_{2\alpha}u_{\lambda_0}^{-1}$,
and the torsion relations are compatible because
\[
 2a_\alpha^2u_{\lambda_0}=2^na_{\lambda_0}u_{2\alpha}=0.
\]

For fixed $V\in RO_{\mathrm{sel}}(G)$,
$E_2^{V+s,s}\ne0\;\Longrightarrow\;s=-\dim V$.
Thus no hidden extensions remain, and the normal forms give
the formula in the corollary.
\end{proof}

\subsection{Valuation language}

Let $g_V$ generate a cyclic group in degree $V$, of order $p^{e_V}$;
$e_V=\infty$ denotes an infinite cyclic group. For $b_V\in\mathbb Z$,
\begin{equation}\label{eq:appendix-fractional-normalization}
 \mathbb Z/p^{e_V}\{g_V\}
 \xrightarrow{\ \cong\ }
 p^{b_V}\mathbb Z/p^{b_V+e_V}\mathbb Z,
 \qquad mg_V\longmapsto [mp^{b_V}].
\end{equation}
Here $p^\infty\mathbb Z=0$, and the fractional lattices lie in
$\mathbb Z[1/p]$. Suppose $g_0=1$, $b_0=0$, and, whenever
$g_V,g_W,g_{V+W}$ are defined,
\[
 b_V+b_W-b_{V+W}\geq0,\qquad
 g_Vg_W=p^{b_V+b_W-b_{V+W}}g_{V+W}.
\]
Under \eqref{eq:appendix-fractional-normalization}, products are
computed by
$(mp^{b_V})(np^{b_W})=mn\,p^{b_V+b_W}$
modulo $p^{b_{V+W}+e_{V+W}}\mathbb Z$. The generator orders give
\[
 e_{V+W}\leq\min(e_V,e_W)+b_V+b_W-b_{V+W},
\]
so the product is well-defined. If the target group is zero, the product is zero.

\subsubsection{Valuation language for $p$ odd}

Use $\rho_p$ and $S_i$ from Section~\ref{sec:local}.
We now describe the result in valuation language.

\begin{proposition}[Valuation language for $p$ odd]\label{prop:odd-fractional-lattices}
For $p$ odd, the rings $\pi_\star^G(H^{h_i})$,
$\pi_\star^G(H^{t_i})$, and $\pi_\star^G(H^{\Phi_n})$
are given by the valuation quotients in \eqref{eq:local-height-table}. Products are computed by multiplying
fractional representatives in $\mathbb Z[1/p]$ and reducing in the target quotient.
\end{proposition}
\begin{proof}
\emph{(1) Degree coordinates and support.}
Since $\dim\lambda_j^{C_{p^k}}=2\mathbf1_{k\leq j}$,
$\dim\rho_p(\mathbf a)^{C_{p^k}} =2a_n+2\sum_{j=k}^{n-1}(a_j-a_{j+1})=2a_k$
for $0\leq k\leq n$. Thus $\rho_p$ is injective. For $E=H^{h_i},H^{t_i},H^{\Phi_n}$ and
$V\in RO_{\mathrm{sel}}(G)\setminus\operatorname{im}\rho_p$,
the normal forms below show that
\[
 \pi_V^G(E)=0.
\]

\smallskip
\noindent\emph{(2) Generators and orders.}
For $0\leq i<n$, put
\[
 \Lambda_i^{(p)}=
 \mathbb Z[
  \{a_{\lambda_j}^{\pm1}\}_{0\leq j<i},
  \{u_{\lambda_k}^{\pm1}\}_{i\leq k\leq n-1}]
\]
and define
\begin{equation}\label{eq:odd-forest-w}
 w_i(\mathbf a):=
 \prod_{j=0}^{i-1}a_{\lambda_j}^{a_{j+1}-a_j}
 \prod_{k=i}^{n-1}u_{\lambda_k}^{a_{k+1}-a_k}.
\end{equation}
Empty products are one, and
$|w_i(\mathbf a)|=\rho_p(\mathbf a)-2a_i$.

For $1\leq i\leq n-1$, put
$q_i=\frac{u_{\lambda_{i-1}}}{a_{\lambda_{i-1}}}$, $x_i=\frac{a_{\lambda_i}}{u_{\lambda_i}}$.
Theorem~\ref{thm:odd-extensions} gives
\[
 \pi_\star^G(H^{h_i})\cong
 \frac{\Lambda_i^{(p)}[q_i,x_i]}
 {(p^{n-i+1},q_ix_i-p,p^{n-i}x_i)}.
\]
For $\mathbf a\in\mathbb Z^{n+1}$, let
\begin{equation}\label{eq:odd-forest-generator}
 g_{A_i}(\mathbf a)=w_i(\mathbf a)
 \begin{cases}
  q_i^{a_i},&a_i\geq0,\\
  x_i^{-a_i},&a_i<0.
 \end{cases}
\end{equation}
$|g_{A_i}(\mathbf a)|=\rho_p(\mathbf a)$, and
Lemma~\ref{lem:odd-normal} gives
\begin{equation*}
 \pi_{\rho_p(\mathbf a)}^G(H^{h_i})\cong
 \begin{cases}
  \mathbb Z/p^{n-i}\{g_{A_i}(\mathbf a)\},&a_i<0,\\
  \mathbb Z/p^{n-i+1}\{g_{A_i}(\mathbf a)\},&a_i\geq0.
 \end{cases}
\end{equation*}

For $i=0$, put $x_0=a_{\lambda_0}u_{\lambda_0}^{-1}$.  The presentation
in the proof of Corollary~\ref{cor:odd-ordinary} is
\[
 \pi_\star^G(H^{h_0})\cong
 \Lambda_0^{(p)}[x_0]/(p^nx_0).
\]
For $a_0\leq0$, put $g_{A_0}(\mathbf a)=w_0(\mathbf a)x_0^{-a_0}$.
Then
\[
 \pi_{\rho_p(\mathbf a)}^G(H^{h_0})\cong
 \begin{cases}
  \mathbb Z/p^n\{g_{A_0}(\mathbf a)\},&a_0<0,\\
  \mathbb Z\{g_{A_0}(\mathbf a)\},&a_0=0,\\
  0,&a_0>0.
 \end{cases}
\]

Since $a_{\lambda_i}=x_i u_{\lambda_i}$ with $u_{\lambda_i}$ invertible,
the preceding presentations give, for $0\leq i<n$,
\begin{equation*}
 \pi_\star^G(H^{t_i})
 \cong\pi_\star^G(H^{h_i})[a_{\lambda_i}^{-1}]
 \cong(\Lambda_i^{(p)}/p^{n-i})[x_i^{\pm1}].
\end{equation*}
The class $g_{T_i}(\mathbf a)=w_i(\mathbf a)x_i^{-a_i}$ has degree
$\rho_p(\mathbf a)$ and order $p^{n-i}$, giving $e_{T_i}=n-i$.

In $H^{\Phi_n}$, \eqref{eq:layer-localizations} makes
$a_{\lambda_j}$ invertible for $0\leq j<n$,
and $p a_{\lambda_{n-1}}=0$ gives $p=0$. Put
$q_n=u_{\lambda_{n-1}}a_{\lambda_{n-1}}^{-1}$, of degree $2$.
The cellular argument in the proof of
Proposition~\ref{prop:odd-E2} also applies with $K=G$
and gives the integer-graded ring $\mathbb F_p[q_n]$.
Multiplication by the invertible Euler classes gives
\begin{equation}\label{eq:odd-phin-ring}
 \pi_\star^G(H^{\Phi_n})
 \cong\mathbb F_p[
  a_{\lambda_0}^{\pm1},\ldots,a_{\lambda_{n-1}}^{\pm1},q_n].
\end{equation}
For $a_n\geq0$, the degree-$\rho_p(\mathbf a)$ group is generated by
\[
 g_{A_n}(\mathbf a)=
 \left(\prod_{j=0}^{n-1}
  a_{\lambda_j}^{a_{j+1}-a_j}\right)q_n^{a_n},
\]
which has order $p$. The group is zero for $a_n<0$.

\smallskip
\noindent\emph{(3) Valuation bounds and multiplication.}
For $1\leq i<n$,
\[
 b_{A_i}(\mathbf a)=S_i(\mathbf a)+a_i^+,
 \qquad e_{A_i}(\mathbf a)=n-i+\mathbf1_{a_i\geq0}.
\]
The relation $q_ix_i=p$ and additivity of $S_i$ give
\begin{equation*}
 g_{A_i}(\mathbf a)g_{A_i}(\mathbf d)
 =p^{a_i^++d_i^+-(a_i+d_i)^+}g_{A_i}(\mathbf a+\mathbf d).
\end{equation*}
Thus
\[
 g_{A_i}(\mathbf a)\longmapsto p^{b_{A_i}(\mathbf a)}
 \pmod {p^{c_{A_i}(\mathbf a)}}
\]
is multiplicative by \eqref{eq:appendix-fractional-normalization}.
For $i=0$ and $A_0(\mathbf a),A_0(\mathbf b)\ne0$,
$g_{A_0}(\mathbf a)g_{A_0}(\mathbf b)=g_{A_0}(\mathbf a+\mathbf b)$, and the
lower bound $S_0$ is additive.  The same argument applies
to $g_{T_i}$ with bound $S_i$ and to $g_{A_n}$ with bound $a_n$.
The degree-zero generator maps to $1$.
\end{proof}

\subsubsection{Valuation language for $p=2$}

Use $D_n^{(0)},D_n^{(1)},\rho_2,\epsilon,t_i$, and $S_i$ from
Section~\ref{sec:local}.
We also give the case $p=2$.

\begin{proposition}[Valuation language for $p=2$]\label{prop:two-fractional-lattices}
For $p=2$, the rings $\pi_\star^G(H^{h_i})$,
$\pi_\star^G(H^{t_i})$, and $\pi_\star^G(H^{\Phi_n})$
are given by the valuation quotients in \eqref{eq:local-height-table}. Products are computed by multiplying
fractional representatives in $\mathbb Z[1/2]$ and reducing in the target quotient.
\end{proposition}
\begin{proof}
\emph{(1) Degree coordinates and support.}
For $0\leq k\leq n$,
$\dim\rho_2(\mathbf a)^{C_{2^k}}=2a_k$.
Thus $\rho_2$ is injective. Moreover,
$2(a_{n-1}-a_n)\equiv\epsilon(\mathbf a)\pmod2$.
The normal forms below give the support in \eqref{eq:local-height-table}.

\smallskip
\noindent\emph{(2) Generators and orders.}
Put $z=a_\alpha u_{2\alpha}^{-1}$.  For
$\mathbf a\in D_n^{(\epsilon)}$, define the Laurent monomial
\begin{equation}\label{eq:two-forest-w}
 w_i(\mathbf a):=
 \prod_{j=0}^{i-1}a_{\lambda_j}^{a_{j+1}-a_j}
 \prod_{j=i}^{n-2}u_{\lambda_j}^{a_{j+1}-a_j}
 u_{2\alpha}^{a_n-a_{n-1}+\epsilon/2}.
\end{equation}
All exponents are integral, and
$|w_i(\mathbf a)z^\epsilon| =\rho_2(\mathbf a)-2t_i(\mathbf a)$.

For $1\leq i<n$, put
$q_i=\frac{u_{\lambda_{i-1}}}{a_{\lambda_{i-1}}}$, $x_i=\frac{a_{\lambda_i}}{u_{\lambda_i}}$.
Theorem~\ref{thm:two-extensions} gives
\[
 \pi_\star^G(H^{h_i})\cong
 \frac{\Lambda_i^{(2)}[q_i,x_i,z]}
 {(2^{n-i+1},q_ix_i-2,2^{n-i}x_i,2z,
   z^2-2^{n-i-1}x_i u_{2\alpha}^{-1})}.
\]
For $\mathbf a\in D_n$, set
\begin{equation}\label{eq:two-forest-generator}
 g_{A_i}(\mathbf a)=w_i(\mathbf a)z^{\epsilon(\mathbf a)}
 \begin{cases}
  q_i^{t_i(\mathbf a)},&t_i(\mathbf a)\geq0,\\
  x_i^{-t_i(\mathbf a)},&t_i(\mathbf a)<0.
 \end{cases}
\end{equation}
$|g_{A_i}(\mathbf a)|=\rho_2(\mathbf a)$, and
Lemma~\ref{lem:two-normal} gives
\[
 \pi_{\rho_2(\mathbf a)}^G(H^{h_i})
 \cong\mathbb Z/2^{e_{A_i}(\mathbf a)}\{g_{A_i}(\mathbf a)\},
\]
where
\begin{equation*}
 2^{e_{A_i}(\mathbf a)}=\begin{cases}
  2^{n-i},&\epsilon(\mathbf a)=0,\ a_i<0,\\
  2^{n-i+1},&\epsilon(\mathbf a)=0,\ a_i\geq0,\\
  2,&\epsilon(\mathbf a)=1.
 \end{cases}
\end{equation*}

For $i=0$, put
\[
 \Lambda_0^{(2)}=
 \mathbb Z[u_{\lambda_0}^{\pm1},\ldots,
 u_{\lambda_{n-2}}^{\pm1},u_{2\alpha}^{\pm1}],
 \qquad x_0=a_{\lambda_0}u_{\lambda_0}^{-1},
\]
with the conventions $\lambda_0=2\alpha$ and $x_0=a_\alpha^2u_{2\alpha}^{-1}$
when $n=1$.  The normal-form argument in the proof of
Corollary~\ref{cor:two-ordinary} gives
\begin{equation*}
 \pi_\star^G(H^{h_0})\cong
 \frac{\Lambda_0^{(2)}[x_0,z]}
 {(2^nx_0,2z,z^2-2^{n-1}x_0u_{2\alpha}^{-1})}.
\end{equation*}
For $\mathbf a\in D_n^{(\epsilon)}$ with $a_0\leq-\epsilon/2$,
\[
 g_{A_0}(\mathbf a)=w_0(\mathbf a)z^\epsilon x_0^{-t_0(\mathbf a)},
 \qquad
 e_{A_0}(\mathbf a)=
 \begin{cases}
  n,&\epsilon=0,\ a_0<0,\\
  \infty,&\epsilon=0,\ a_0=0,\\
  1,&\epsilon=1.
 \end{cases}
\]

For $0\leq i<n$, $a_{\lambda_i}=x_i u_{\lambda_i}$ with
$u_{\lambda_i}$ invertible. Euler localization gives
\begin{equation*}
 \pi_\star^G(H^{t_i})\cong
 \frac{(\Lambda_i^{(2)}/2^{n-i})[x_i^{\pm1},z]}
 {(2z,z^2-2^{n-i-1}x_i u_{2\alpha}^{-1})}.
\end{equation*}
For $\mathbf a\in D_n^{(\epsilon)}$,
\[
 g_{T_i}(\mathbf a)=w_i(\mathbf a)z^\epsilon x_i^{-t_i(\mathbf a)},
 \qquad
 e_{T_i}(\mathbf a)=
 \begin{cases}
  n-i,&\epsilon=0,\\
  1,&\epsilon=1.
 \end{cases}
\]

Apply \cite[Proposition~3.18]{HHRNonexistence} with $K=G$.
After inverting $a_\alpha$ and $a_{\lambda_j}$ for
$0\leq j\leq n-2$, we obtain
\begin{equation}\label{eq:two-phin-ring}
 \pi_\star^G(H^{\Phi_n})\cong
 \mathbb F_2[
  a_\alpha^{\pm1},a_{\lambda_0}^{\pm1},\ldots,
  a_{\lambda_{n-2}}^{\pm1},u_{2\alpha}].
\end{equation}
For $a_n\geq0$, the degree-$\rho_2(\mathbf a)$ group is generated by
\[
 g_{A_n}(\mathbf a)=u_{2\alpha}^{a_n}a_\alpha^{-2a_{n-1}}
 \prod_{j=0}^{n-2}a_{\lambda_j}^{a_{j+1}-a_j},
\]
which has order $2$. The group is zero for $a_n<0$.

\smallskip
\noindent\emph{(3) Valuation bounds and multiplication.}
For $1\leq i<n$, assign weights to formal monomial expressions,
including powers of $2$, by adding the weights of their factors:
\[
 \operatorname{wt}(2)=1,
 \quad \operatorname{wt}(a_{\lambda_j})=0\ (j<i),
 \quad \operatorname{wt}(u_{\lambda_j})=n-j\ (j\geq i),
\]
\[
 \operatorname{wt}(u_{2\alpha})=1,
 \quad \operatorname{wt}(q_i)=n-i+1,
 \quad \operatorname{wt}(x_i)=-(n-i),
 \quad \operatorname{wt}(z)=-1.
\]
The two sides of each relation
$q_ix_i=2$ and
$z^2=2^{n-i-1}x_i u_{2\alpha}^{-1}$
have the same weight. The weight of
$g_{A_i}(\mathbf a)$ is
\begin{equation*}
 b_{A_i}(\mathbf a)=S_i(\mathbf a)
 +t_i(\mathbf a)^+
 +\epsilon(\mathbf a)\frac{n-i-1}{2}.
\end{equation*}
Reducing products by $q_ix_i=2$ and
$z^2=2^{n-i-1}x_i u_{2\alpha}^{-1}$ gives
\[
 g_{A_i}(\mathbf a)g_{A_i}(\mathbf b)
 =2^{b_{A_i}(\mathbf a)+b_{A_i}(\mathbf b)-b_{A_i}(\mathbf a+\mathbf b)}
  g_{A_i}(\mathbf a+\mathbf b).
\]
Before reducing the integer coefficient modulo
$2^{e_{A_i}(\mathbf a+\mathbf b)}$, the exponent counts
applications of $q_ix_i=2$, with an additional $n-i-1$ for each
application of the $z^2$ relation. Thus the exponent is nonnegative,
and the formula includes zero products. By
\eqref{eq:appendix-fractional-normalization},
\[
 g_{A_i}(\mathbf a)\longmapsto 2^{b_{A_i}(\mathbf a)}
 \pmod {2^{c_{A_i}(\mathbf a)}}
\]
is multiplicative.
For $\mathbf a\in D_n^{(1)}$,
$S_i(\mathbf a)\equiv(n-i-1)/2\pmod{\mathbb Z}$, so all lower bounds
are integral.

For $H^{h_0}$, use the same weights without $q_0$.
After localization, $\operatorname{wt}(x_i^{-1})=n-i$ and
$q_i=2x_i^{-1}$ for $i\geq1$. The weights of $g_{A_0}(\mathbf a)$
and $g_{T_i}(\mathbf a)$ are respectively
\[
 S_0(\mathbf a)+\epsilon\frac{n-1}{2},\qquad
 S_i(\mathbf a)+\epsilon\frac{n-i-1}{2}.
\]
The first formula applies when $A_0(\mathbf a)\ne0$.
For $H^{\Phi_n}$ the bound $a_n$ is additive and $e_{A_n}=1$.
In every case the degree-zero generator maps to $1$.
\end{proof}

\subsection{Restriction, transfer, and Weyl action in valuation language}

Let $G=C_{p^n}$, $K=C_{p^{n-1}}$, and $n\geq2$.
Use the $RO_{\mathrm{sel}}(G)$-grading
$\pi_V^K(E)=\pi_{\Res_K^G V}^K(\res_K^GE)$ fixed above.
For an $H$-algebra $E$ and $x\in\pi_\star^G(E)$,
\begin{equation}\label{eq:appendix-layer-tr-res}
 \tr_K^G\res_K^G(x)
 =x\tr_K^G(1)=px.
\end{equation}
Below, $E$ is one of $H^{h_i}$, $H^{t_i}$ $(0\leq i<n)$, or
$H^{\Phi_n}$, with the generators constructed in
Propositions~\ref{prop:odd-fractional-lattices} and
\ref{prop:two-fractional-lattices}.

For $V=\rho_p(\mathbf a)$ with $\mathbf a\in\mathbb Z^{n+1}$, if
$\res_K^G E\ne0$, $\pi_V^G(E)\ne0$, $\pi_V^K(E)=0$,
then $E=H^{h_{n-1}}$, $a_{n-1}<0$, and $\pi_V^G(E)\cong\mathbb Z/p$.
In this case, for $y\in\pi_V^G(E)$ and
$\widetilde x\in\pi_W^G(E)$, put $x=\res_K^G(\widetilde x)$.
Equation~\eqref{eq:appendix-layer-tr-res} gives
\[
 \res_K^G(y)=0,\qquad py=0,\qquad
 \tr_K^G(x)y=p\widetilde x y=0
 =\tr_K^G\bigl(x\res_K^G(y)\bigr).
\]
Since $\res_K^G H^{t_{n-1}}\simeq\res_K^G H^{\Phi_n}\simeq0$,
Frobenius reciprocity for $H^{t_{n-1}}$ and $H^{\Phi_n}$ reduces to $0=0$.

\paragraph{$p$ odd.}
For $\mathbf a=(a_0,\ldots,a_n)\in\mathbb Z^{n+1}$,
\[
 \Res_K^G\rho_p(\mathbf a)
 =2a_{n-1}+\sum_{j=0}^{n-2}(a_j-a_{j+1})\lambda_j|_K.
\]
Here $\lambda_{n-1}|_K\cong\mathbb R^2$, and the selected $K$-planes are
$\lambda_j|_K$ for $0\leq j\leq n-2$, with their restricted orientations.
The formulas of \eqref{eq:local-height-table} for $C_{p^{n-1}}$ at $(a_0,\ldots,a_{n-1})$
give the $K$-equivariant homotopy groups of $H^{h_i}$ and $H^{t_i}$. For $i=n-1$,
\[
 \pi_{\rho_p(\mathbf a)}^K(H^{h_{n-1}})
 \cong
 \begin{cases}
 p^{a_{n-1}}\mathbb Z/p^{a_{n-1}+1}\mathbb Z,&a_{n-1}\geq0,\\
 0,&a_{n-1}<0,
 \end{cases}
 \qquad
 \pi_{\rho_p(\mathbf a)}^K(H^{t_{n-1}})=0.
\]
For $V=\rho_p(\mathbf a)$ with $\pi_V^G(E),\pi_V^K(E)\ne0$, let
$g_G,g_K$ be the chosen generators and write
$\pi_V^G(E)\cong p^{b_G}\mathbb Z/p^{b_G+e_G}\mathbb Z$, $\pi_V^K(E)\cong p^{b_K}\mathbb Z/p^{b_K+e_K}\mathbb Z$.
Equation~\eqref{eq:local-height-table} gives
\begin{equation}\label{eq:appendix-subgroup-exponents}
 b_K=b_G-a_n,\qquad
 e_K=e_G-1\quad(e_G<\infty),\qquad
 e_K=e_G=\infty\quad(e_G=\infty).
\end{equation}
If $e_G=1$, the group $\pi_V^K(E)$ is zero.

We now give the subgroup maps.

\begin{proposition}[Restriction and transfer maps]\label{prop:odd-local-subgroup-maps}
For $p$ odd, $E=H^{h_i}$ or $H^{t_i}$ $(0\leq i<n)$,
$V=\rho_p(\mathbf a)$, and $\pi_V^G(E),\pi_V^K(E)\ne0$,
\begin{equation*}
 \res_K^G([x])=[p^{-a_n}x],\qquad
 \tr_K^G([y])=[p^{a_n+1}y].
\end{equation*}
Maps involving a zero group are zero. The Weyl action is trivial.
\end{proposition}
\begin{proof}
Naturality of Euler and orientation classes gives
\[
 \res_K^G(a_{\lambda_j})=a_{\lambda_j|_K},\qquad
 \res_K^G(u_{\lambda_j})=u_{\lambda_j|_K},\qquad
 a_{\lambda_{n-1}|_K}=0,\quad u_{\lambda_{n-1}|_K}=1.
\]
Substitution in \eqref{eq:odd-forest-generator} gives
$\res_K^G(g_G)=g_K$ whenever $\pi_V^K(E)\ne0$.
For $i=n-1$ and $a_{n-1}<0$, $\pi_V^K(E)=0$.
Euler localization gives the same formula for $H^{t_i}$.
Under \eqref{eq:appendix-fractional-normalization},
\[
 g_G\longleftrightarrow [p^{b_G}],\qquad
 g_K\longleftrightarrow [p^{b_K}]=[p^{-a_n}p^{b_G}].
\]
Equation~\eqref{eq:appendix-layer-tr-res} gives
$\tr_K^G(g_K)=pg_G$. The transfer formula is well-defined because
$p^{a_n+1}p^{b_K+e_K}=p^{b_G+e_G}$
when $e_G<\infty$; for $e_G=\infty$ there is no annihilator relation.
Every chosen generator at $K$ is a restriction and is therefore $G/K$-fixed.

For $\mathbf a,\mathbf d\in\mathbb Z^{n+1}$ and classes
$[x]\in\pi_{\rho_p(\mathbf a)}^G(E)$,
$[y]\in\pi_{\rho_p(\mathbf d)}^G(E)$, multiplicativity of restriction follows from
$p^{-(a_n+d_n)}xy=(p^{-a_n}x)(p^{-d_n}y)$.
For $[y]\in\pi_{\rho_p(\mathbf a)}^K(E)$ and
$[x]\in\pi_{\rho_p(\mathbf d)}^G(E)$, Frobenius reciprocity follows from
$p^{a_n+d_n+1}\bigl(y\,p^{-d_n}x\bigr) =(p^{a_n+1}y)x$.
The inequalities for $b_E$ and $c_E$ established in the proof of
Theorem~\ref{thm:local-padic-data} make the reductions in the target
quotients well-defined.
\end{proof}

\paragraph{$p=2$.}
Let $\alpha_K$ be the sign representation of $K$ and give
$\lambda_j|_K$, $0\leq j\leq n-3$, the restricted orientations.
Then $\lambda_{n-2}|_K=2\alpha_K$ and
\[
 \Res_K^G\rho_2(\mathbf a)
 =2a_{n-1}+2(a_{n-2}-a_{n-1})\alpha_K
       +\sum_{j=0}^{n-3}(a_j-a_{j+1})\lambda_j|_K.
\]
The sum is empty for $n=2$. A group $\pi_{\rho_2(\mathbf a)}^K(E)$ can be
nonzero only for
\[
 (a_0,\ldots,a_{n-1})\in\mathbb Z^n,
 \qquad a_n\in\tfrac12\mathbb Z,
\]
subject to the conditions in \eqref{eq:local-height-table} with $n-1$ in place of $n$.
The groups at $K$ are given by the valuation quotients for
$C_{2^{n-1}}$ at $(a_0,\ldots,a_{n-1})$. In particular,
\[
 \pi_{\rho_2(\mathbf a)}^K(H^{h_{n-1}})
 \cong
 \begin{cases}
 2^{a_{n-1}}\mathbb Z/2^{a_{n-1}+1}\mathbb Z,&a_{n-1}\geq0,\\
 0,&a_{n-1}<0,
 \end{cases}
 \qquad
 \pi_{\rho_2(\mathbf a)}^K(H^{t_{n-1}})=0.
\]

We have the following formulas for $p=2$.

\begin{proposition}[Restriction and transfer maps and Weyl action for $p=2$]\label{prop:two-local-subgroup-maps}
Let $p=2$ and $V=\rho_2(\mathbf a)$. If $\mathbf a\in\mathbb Z^{n+1}$ and
$\pi_V^G(E),\pi_V^K(E)\ne0$, \eqref{eq:appendix-subgroup-exponents} and the formulas in
Proposition~\ref{prop:odd-local-subgroup-maps} hold with $p=2$.
If $\pi_V^K(E)\ne0$ and $a_n\in\mathbb Z+\tfrac12$, then
$\pi_V^G(E)=0$, so restriction and transfer are zero.
For $\mathbf a\in D_n^{(1)}$, $\pi_V^K(E)=0$.
For $\pi_V^K(E)\ne0$,
\begin{equation*}
 c_{\gamma_{2^n}}([y])=(-1)^{2a_n}[y].
\end{equation*}
Restriction is multiplicative and transfer satisfies Frobenius reciprocity.
\end{proposition}
\begin{proof}
The fixed-point parity conditions give
\[
\begin{aligned}
 \mathbf a\in D_n^{(1)}
 &\Longrightarrow \dim V^K=2a_{n-1}\in2\mathbb Z+1
 \Longrightarrow \pi_V^K(E)=0,\\
 a_n\in\mathbb Z+\tfrac12
 &\Longrightarrow \dim V^G=2a_n\in2\mathbb Z+1
 \Longrightarrow \pi_V^G(E)=0.
\end{aligned}
\]
Naturality gives
\[
 \res_K^G(a_\alpha)=0,\quad
 \res_K^G(u_{2\alpha})=1,\quad
 \res_K^G(a_{\lambda_{n-2}})=a_{\alpha_K}^2,\quad
 \res_K^G(u_{\lambda_{n-2}})=u_{2\alpha_K}.
\]
Thus \eqref{eq:two-forest-generator} restricts to the chosen lower
generator when both groups are nonzero. Moreover,
\[
 \res_K^G(z)=0,\qquad
 \res_K^G(x_{n-1})=\res_K^G(a_\alpha^2u_{2\alpha}^{-1})=0.
\]
Equations~\eqref{eq:appendix-layer-tr-res} and
\eqref{eq:appendix-subgroup-exponents} give the well-defined transfer formula.

To verify products of two classes with degree coordinates in
$D_n^{(1)}$, put
$U=\rho_2(\mathbf a)$, $W=\rho_2(\mathbf d)$ with
$\mathbf a,\mathbf d\in D_n^{(1)}$.
For $E=H^{h_i}$, $1\leq i<n$, the lower bounds satisfy
\[
 \begin{aligned}
 \delta&=b_G(U)+b_G(W)-b_G(U+W)\\
 &=n-i-1+t_i(\mathbf a)^++t_i(\mathbf d)^+
   -(t_i(\mathbf a)+t_i(\mathbf d)-1)^+\\
 &\geq
 \begin{cases}
 n-i-1,&t_i(\mathbf a)+t_i(\mathbf d)-1<0,\\
 n-i,&t_i(\mathbf a)+t_i(\mathbf d)-1\geq0.
 \end{cases}
 \end{aligned}
\]
The corresponding lengths of $\pi_{U+W}^K(E)$ are $n-i-1$ and $n-i$.
For $E=H^{h_0}$ or $H^{t_i}$, the difference $\delta$ is $n-1$ or
$n-i-1$, respectively, equal to the length of $\pi_{U+W}^K(E)$.
When $\pi_{U+W}^K(E)\ne0$, take representatives
$q\in2^{b_G(U)}\mathbb Z$ and $q'\in2^{b_G(W)}\mathbb Z$,
\[
 2^{-(a_n+d_n)}qq'
 \in2^{b_K(U+W)+\delta}\mathbb Z
 \subseteq2^{b_K(U+W)+e_K(U+W)}\mathbb Z.
\]
Thus products of two classes with degree coordinates in
$D_n^{(1)}$ restrict to zero.
Products for which exactly one input has degree coordinates in
$D_n^{(1)}$ have zero target homotopy group at $K$.
When $a_n\in\mathbb Z+\tfrac12$, multiplication by the restriction of
a class at $G$ with integral degree coordinates leaves the last coordinate half-integral, so transfer
is zero on both sides of Frobenius reciprocity.
For $x\in\pi_V^G(E)$ with $\mathbf a\in D_n^{(1)}$,
let $W$ be any degree with $\pi_W^G(E),\pi_W^K(E)\ne0$.
For $y\in\pi_W^K(E)$, choose $\widetilde y\in\pi_W^G(E)$ with
$\res_K^G(\widetilde y)=y$. Then
\[
 \tr_K^G(y)x=2\widetilde yx=0,
\]
since $\res_K^G(x)=0$ and $\pi_{V+W}^G(E)$ is annihilated by $2$.
For degree coordinates in $D_n^{(0)}$ with nonzero $G$- and
$K$-groups, the identities for fractional representatives in the
proof of Proposition~\ref{prop:odd-local-subgroup-maps} apply.

Finally, the chosen oriented isomorphism $\alpha|_K\cong1$
gives a unit in $\pi_{1-\alpha}^K(E)$.
The nontrivial element of $G/K$ reverses the orientation of this line.
The chosen generator of $\pi_{\rho_2(\mathbf a)}^K(E)$
is the product of the $2a_n$th power of this unit and the
restriction of the chosen generator in degree
$\rho_2(a_0,\ldots,a_{n-1},0)$.
Therefore its Weyl action is $(-1)^{2a_n}$.
Since $2(a_n+d_n)=2a_n+2d_n$, this action preserves products.
\end{proof}

We also have the following proposition.

\begin{proposition}[$H^{\Phi_n}$ and $H^{\Phi_{n-1}}$ at $K$]\label{prop:terminal-local-subgroup-values}
For $n\geq2$ and $V\in RO(G)$,
\[
 \pi_V^K(H^{\Phi_n})=0.
\]
For $\mathbf a\in\mathbb Z^{n+1}$ when $p$ is odd, and
$\mathbf a\in\mathbb Z^n\times\tfrac12\mathbb Z$ when $p=2$,
\[
 \pi_{\rho_p(\mathbf a)}^K(H^{\Phi_{n-1}})
 \cong
 \begin{cases}
 p^{a_{n-1}}\mathbb Z/p^{a_{n-1}+1}\mathbb Z,&a_{n-1}\geq0,\\
 0,&a_{n-1}<0.
 \end{cases}
\]
The groups $\pi_V^G(H^{\Phi_{n-1}})$ follow from the pullback
$H^{\Phi_{n-1}}\simeq H^{h_{n-1}}\times^h_{H^{t_{n-1}}}H^{\Phi_n}$.
\end{proposition}
\begin{proof}
Since $K\in\mathcal F_{n-1}$,
$\res_K^G\widetilde E\mathcal F_{n-1}\simeq *$.
For $H^{\Phi_{n-1}}$, restricting the family gives the proper-subgroup
family of $K$. Apply \eqref{eq:odd-phin-ring} or
\eqref{eq:two-phin-ring} to $K$, with $n-1$ in place of $n$;
the formulas in valuation language follow from the preceding subsection.
\end{proof}

We now consider the first inclusion.

\begin{proposition}\label{prop:appendix-prime-base}
Let $G=C_p$ and $K=e$. Then
\[
 \res_e^G(H^{h_0})\simeq H\mathbb Z,\qquad
 \res_e^G(H^{t_0})\simeq0,\qquad \res_e^G(H^{\Phi_1})\simeq0.
\]
For $V\in RO(G)$ with $\dim V=0$ and trivial orientation character,
write $V=\rho_p(\mathbf a)+D$ with the selected representative
$\rho_p(\mathbf a)$ and $v_D$ as in
Lemma~\ref{lem:normalized-character-units}. Put
$g_V=v_Dg_{A_0}(\mathbf a)\in\pi_V^G(H^{h_0})\cong\mathbb Z$,
using the image of $v_D$ under $H\to H^{h_0}$. Then
\[
 \res_e^G(g_V)=1,\qquad
 \tr_e^G(1)=p g_V.
\]
\end{proposition}
\begin{proof}
Restriction to $e$ makes $EG$ contractible and $\widetilde EG$
contractible as a based space, giving the three equivalences.
For $\dim V=0$, the ordinary HFPSS gives
\[
 \res_e^G\colon
 \pi_V^G(H^{h_0})
 \xrightarrow{\cong} H^0(G;\mathbb Z_{\varepsilon_V})
 \lhook\joinrel\longrightarrow \mathbb Z_{\varepsilon_V},
\]
where $\varepsilon_V$ is the orientation character and
$\mathbb Z_{\varepsilon_V}$ is $\mathbb Z$ with $G$-action through
$\varepsilon_V$. The displayed inclusion is the identity for trivial
$\varepsilon_V$ and has zero source otherwise.
Frobenius reciprocity gives
$\tr_e^G\res_e^G(g_V)=p g_V$.
The remaining assertions follow from $\pi_d(H\mathbb Z)=0$ for $d\ne0$
and the action on the oriented representation sphere.
\end{proof}

\begin{remark}
For $p=2$, if $\dim V=0$ and the orientation character is nontrivial,
then $\pi_V^G(H^{h_0})=0$, so restriction and transfer are zero.
If $\dim V\ne0$, the group $\pi_V^e(H^{h_0})$ is zero. All restriction and
transfer maps for $H^{t_0}$ and $H^{\Phi_1}$ are zero.
On the nonzero group $\pi_V^e(H^{h_0})$, an element of $G$ acts by the
orientation character of $V$.
\end{remark}

\subsection{Comparison maps}

We need the following lemma.

\begin{lemma}[The Euler class in valuation language]
\label{lem:local-euler-verification}
For $0\leq i<n-1$, the $H^{t_i}$-module map
$H^{t_i}\to S^{\lambda_{i+1}}\wedge H^{t_i}$ induced by
$a_{\lambda_{i+1}}$ corresponds to the class of $1$ in
\[
 p^{-1}\mathbb Z/p^{n-i-1}\mathbb Z
 \cong \pi_{-\lambda_{i+1}}^G(H^{t_i}).
\]
\end{lemma}
\begin{proof}
Put $x_i=a_{\lambda_i}/u_{\lambda_i}$. The adjacent gold relation gives
\[
 a_{\lambda_{i+1}}=p\,u_{\lambda_{i+1}}x_i.
\]
For $p=2$, $i=n-2$, this is the top gold relation with
$\lambda_{n-1}=2\alpha$ and $u_{\lambda_{n-1}}=u_{2\alpha}$.
In degree $-\lambda_{i+1}$, the lower bound is $S_i=-1$, so
\[
 u_{\lambda_{i+1}}x_i\longleftrightarrow[p^{-1}]
 \in p^{-1}\mathbb Z/p^{n-i-1}\mathbb Z,
 \qquad p[p^{-1}]=[1].\qedhere
\]
\end{proof}

We can now give the formulas for the comparison maps.

\begin{lemma}[Comparison maps]
\label{lem:local-height-normalization}
In valuation language, each map
$\ell_i:A_i(\mathbf a)\to T_i(\mathbf a)$ or
$\rho_i:A_{i+1}(\mathbf a)\to T_i(\mathbf a)$ sends $[q]$ to $[q]$
in $T_i(\mathbf a)$.
\end{lemma}
\begin{proof}
We compute in $H^{t_i}$ with $x_i=a_{\lambda_i}/u_{\lambda_i}$.
For odd $p$, we take $\epsilon=0$ and $t_j=a_j$. At $p=2$, we retain
$\epsilon$, $t_j=a_j+\epsilon/2$, and $z=a_\alpha u_{2\alpha}^{-1}$.
The factor $z^\epsilon$ is omitted for odd $p$.
We use $w_i(\mathbf a)$ for the Laurent monomial in
\eqref{eq:odd-forest-w} or \eqref{eq:two-forest-w}.
Under \eqref{eq:appendix-fractional-normalization},
\[
 g_{A_i}(\mathbf a)\longleftrightarrow[p^{b_{A_i}(\mathbf a)}],
 \qquad
 g_{T_i}(\mathbf a)\longleftrightarrow[p^{b_{T_i}(\mathbf a)}].
\]
The chosen generator of $T_i(\mathbf a)$ is
\[
 g_{T_i}(\mathbf a)=w_i(\mathbf a)z^\epsilon x_i^{-t_i}.
\]
For $\ell_i$, $1\leq i<n$, the relation $q_ix_i=p$ gives
\[
 \begin{aligned}
 \ell_i(g_{A_i}(\mathbf a))&=p^{t_i^+}g_{T_i}(\mathbf a),\\
 b_{A_i}(\mathbf a)-b_{T_i}(\mathbf a)&=t_i^+.
 \end{aligned}
\]
Thus $\ell_i([q])=[q]$.
If $\epsilon=1$ and $t_i>0$, then $\ell_i(g_{A_i}(\mathbf a))=0$
since $e_{T_i}=1$.

For $\rho_i$, $0\leq i<n-1$, the monomials satisfy
\[
 \frac{w_{i+1}(\mathbf a)}{w_i(\mathbf a)}
 =x_i^{a_{i+1}-a_i},\qquad
 q_{i+1}=x_i^{-1},\qquad x_{i+1}=p x_i.
\]
The last relation also holds for $p=2$, $i=n-2$.
If $t_{i+1}\geq0$, the exponent after substitution is
$a_{i+1}-a_i-t_{i+1}=-t_i$. If $t_{i+1}<0$, the substitution
contributes $p^{-t_{i+1}}$. Therefore
\[
 \begin{aligned}
 \rho_i(g_{A_{i+1}}(\mathbf a))
 &=p^{(-t_{i+1})^+}g_{T_i}(\mathbf a),\\
 b_{A_{i+1}}(\mathbf a)-b_{T_i}(\mathbf a)
 &=(-t_{i+1})^+.
 \end{aligned}
\]
Lemma~\ref{lem:local-euler-height}, proved using
Lemma~\ref{lem:local-euler-verification}, gives the Euler nilpotence
needed for the factorization
\[
 H^{\Phi_{i+1}}\longrightarrow H^{h_{i+1}}
 \xrightarrow{\rho_i}H^{t_i}.
\]
By \eqref{eq:complete-mapping-equivalence}, $\rho_i$ sends the
images of Euler and orientation classes to their images in $H^{t_i}$.
Thus the displayed identities compute $\rho_i$.

For $i=0$ and $A_0(\mathbf a)\ne0$, the generators $g_{A_0}(\mathbf a)$
and $g_{T_0}(\mathbf a)$ are both $w_0z^\epsilon x_0^{-t_0}$, and
$b_{A_0}(\mathbf a)=b_{T_0}(\mathbf a)$. Thus $\ell_0([q])=[q]$.

For $\rho_{n-1}$, the generators $g_{A_n}(\mathbf a)$ and $g_{T_{n-1}}(\mathbf a)$
are represented in $\pi_\star^G(H^{\Phi_n})$ and
$\pi_\star^G(H^{t_{n-1}})$, respectively, as follows. For odd $p$,
\[
 \left(\prod_{j<n-1}a_{\lambda_j}^{a_{j+1}-a_j}\right)
 u_{\lambda_{n-1}}^{a_n}a_{\lambda_{n-1}}^{-a_{n-1}},
\]
and for $p=2$,
\[
 u_{2\alpha}^{a_n}a_\alpha^{-2a_{n-1}}
 \prod_{j=0}^{n-2}a_{\lambda_j}^{a_{j+1}-a_j}.
\]
These formulas hold for $\mathbf a\in D_n$ with $a_n\geq0$;
for $a_n<0$, $A_n(\mathbf a)=0$.
The bounds $b_{A_n}(\mathbf a)$ and $b_{T_{n-1}}(\mathbf a)$
both equal $a_n$.
Hence $\rho_{n-1}([q])=[q]$.
\end{proof}

\section*{Acknowledgements}

I would like to thank Professor Guchuan Li for helpful discussions on this problem.
\enlargethispage{2\baselineskip}
\bibliographystyle{alpha}
\bibliography{reference}

@incollection{Boardman,
  author  = {J. Michael Boardman},
  title   = {Conditionally convergent spectral sequences},
  booktitle = {Homotopy Invariant Algebraic Structures},
  series  = {Contemporary Mathematics},
  volume  = {239},
  publisher = {American Mathematical Society},
  address = {Providence, RI},
  year    = {1999},
  pages   = {49--84},
  doi     = {10.1090/conm/239/03597}
}

@article{SikoraYan,
  author  = {Igor Sikora and Guoqi Yan},
  title   = {On the structure of the {$RO(G)$}-graded homotopy of {$H\underline M$} for cyclic {$p$}-groups},
  journal = {Algebraic \& Geometric Topology},
  volume  = {25},
  year    = {2025},
  pages   = {2933--2980},
  doi     = {10.2140/agt.2025.25.2933},
  url     = {https://arxiv.org/abs/2309.17242}
}

@article{Yan,
  author  = {Guoqi Yan},
  title   = {The {$RO(C_{2^n})$}-graded homotopy of {$H\underline{\mathbb{Z}}$} through generalized {Tate} squares},
  journal = {La Matematica},
  volume  = {4},
  year    = {2025},
  pages   = {364--418},
  doi     = {10.1007/s44007-025-00153-1},
  url     = {https://doi.org/10.1007/s44007-025-00153-1}
}

@misc{Zeng,
  author        = {Mingcong Zeng},
  title         = {Equivariant {Eilenberg--Mac Lane} spectra in cyclic {$p$}--groups},
  year          = {2018},
  eprint        = {1710.01769},
  archivePrefix = {arXiv},
  primaryClass  = {math.AT},
  note = {\href{https://arxiv.org/abs/1710.01769}{arXiv:1710.01769}},
  url           = {https://arxiv.org/abs/1710.01769}
}

@misc{Angeltveit,
  author = {Vigleik Angeltveit},
  title = {Recovering a cohomological {Mackey} functor from its restriction to {Sylow} subgroups},
  year = {2022},
  eprint = {2201.01444},
  archivePrefix = {arXiv},
  primaryClass = {math.AT},
  note = {\href{https://arxiv.org/abs/2201.01444}{arXiv:2201.01444}},
  url = {https://arxiv.org/abs/2201.01444}
}

@misc{AMGR,
  author        = {David Ayala and Aaron Mazel-Gee and Nick Rozenblyum},
  title         = {Derived {Mackey} functors and {$C_{p^n}$}-equivariant cohomology},
  year          = {2023},
  eprint        = {2105.02456},
  archivePrefix = {arXiv},
  primaryClass  = {math.AT},
  note          = {\href{https://arxiv.org/abs/2105.02456}{arXiv:2105.02456}},
  url           = {https://arxiv.org/abs/2105.02456}
}

@misc{DuggerHazelAlg,
  author        = {Daniel Dugger and Christy Hazel},
  title         = {{Mackey} homological algebra over cyclic groups},
  year          = {2026},
  eprint        = {2602.18670},
  archivePrefix = {arXiv},
  primaryClass  = {math.AT},
  note = {\href{https://arxiv.org/abs/2602.18670}{arXiv:2602.18670}},
  url           = {https://arxiv.org/abs/2602.18670}
}

@misc{DuggerHazelPoint,
  author        = {Daniel Dugger and Christy Hazel},
  title         = {The {Bredon} equivariant cohomology of a point for cyclic groups},
  year          = {2026},
  eprint        = {2602.18680},
  archivePrefix = {arXiv},
  primaryClass  = {math.AT},
  note = {\href{https://arxiv.org/abs/2602.18680}{arXiv:2602.18680}},
  url           = {https://arxiv.org/abs/2602.18680}
}

@article{QuigleyShah,
  author        = {J. D. Quigley and Jay Shah},
  title         = {On the parameterized {Tate} construction},
  journal       = {Journal of Topology},
  volume        = {18},
  number        = {1},
  year          = {2025},
  pages         = {e70018},
  doi           = {10.1112/topo.70018},
  eprint        = {2110.07707},
  archivePrefix = {arXiv},
  primaryClass  = {math.AT},
  url           = {https://doi.org/10.1112/topo.70018}
}

@article{BasuDey,
  author        = {Samik Basu and Pinka Dey},
  title         = {Equivariant cohomology for cyclic groups},
  journal       = {International Mathematics Research Notices},
  volume        = {2025},
  number        = {12},
  year          = {2025},
  pages         = {rnaf150},
  doi           = {10.1093/imrn/rnaf150},
  eprint        = {2403.00362},
  archivePrefix = {arXiv},
  primaryClass  = {math.AT},
  url           = {https://doi.org/10.1093/imrn/rnaf150}
}

@article{AngeltveitBohmann,
  author  = {Vigleik Angeltveit and Anna Marie Bohmann},
  title   = {Graded {Tambara} functors},
  journal = {Journal of Pure and Applied Algebra},
  volume  = {222},
  number  = {12},
  year    = {2018},
  pages   = {4126--4150},
  doi     = {10.1016/j.jpaa.2018.02.023},
  eprint  = {1504.00668},
  archivePrefix = {arXiv},
  primaryClass = {math.AT}
}

@article{BlumbergHill,
  author  = {Andrew J. Blumberg and Michael A. Hill},
  title   = {Operadic multiplications in equivariant spectra, norms, and transfers},
  journal = {Advances in Mathematics},
  volume  = {285},
  year    = {2015},
  pages   = {658--708},
  doi     = {10.1016/j.aim.2015.07.013},
  eprint  = {1309.1750},
  archivePrefix = {arXiv},
  primaryClass = {math.AT}
}

@article{Balderrama,
  author  = {William Balderrama},
  title   = {Total power operations in spectral sequences},
  journal = {Transactions of the American Mathematical Society},
  volume  = {377},
  number  = {7},
  year    = {2024},
  pages   = {4779--4823},
  doi     = {10.1090/tran/9073},
  eprint  = {2205.07409},
  note    = {Theorem numbering refers to \href{https://arxiv.org/abs/2205.07409}{arXiv:2205.07409}},
  archivePrefix = {arXiv},
  primaryClass = {math.AT}
}

@article{Carrick,
  author        = {Christian Carrick},
  title         = {Smashing localizations in equivariant stable homotopy},
  journal       = {Journal of Homotopy and Related Structures},
  volume        = {17},
  year          = {2022},
  pages         = {355--392},
  doi           = {10.1007/s40062-022-00310-1},
  eprint        = {1909.08771},
  archivePrefix = {arXiv},
  primaryClass  = {math.AT}
}

@misc{UllmanTambara,
  author        = {John Ullman},
  title         = {Tambara Functors and Commutative Ring Spectra},
  year          = {2013},
  eprint        = {1304.4912},
  archivePrefix = {arXiv},
  primaryClass  = {math.AT},
  note          = {arXiv:1304.4912}
}

@book{HHRbook,
  author    = {Michael A. Hill and Michael J. Hopkins and Douglas C. Ravenel},
  title     = {Equivariant Stable Homotopy Theory and the Kervaire Invariant Problem},
  publisher = {Cambridge University Press},
  year      = {2021},
  doi       = {10.1017/9781108917278},
  isbn      = {9781108831444}
}

@article{HillAQ,
  author  = {Michael A. Hill},
  title   = {On the {Andr{\'e}--Quillen} homology of {Tambara} functors},
  journal = {Journal of Algebra},
  volume  = {489},
  year    = {2017},
  pages   = {115--137},
  doi     = {10.1016/j.jalgebra.2017.06.029},
  eprint  = {1701.06219},
  archivePrefix = {arXiv},
  primaryClass = {math.AT}
}

@article{HillMazur,
  author  = {Michael A. Hill and Kristen Mazur},
  title   = {An equivariant tensor product on {Mackey} functors},
  journal = {Journal of Pure and Applied Algebra},
  volume  = {223},
  number  = {12},
  year    = {2019},
  pages   = {5310--5345},
  doi     = {10.1016/j.jpaa.2019.04.001},
  eprint  = {1508.04062},
  archivePrefix = {arXiv},
  note = {Theorem numbering refers to \href{https://arxiv.org/abs/1508.04062}{arXiv:1508.04062}}
}

@article{MNN,
  author = {Akhil Mathew and Niko Naumann and Justin Noel},
  title = {Nilpotence and descent in equivariant stable homotopy theory},
  journal = {Advances in Mathematics},
  volume = {305},
  year = {2017},
  pages = {994--1084},
  eprint = {1507.06869},
  archivePrefix = {arXiv},
  primaryClass = {math.AT},
  url = {https://arxiv.org/abs/1507.06869}
}

@article{HHRNonexistence,
  author = {Michael A. Hill and Michael J. Hopkins and Douglas C. Ravenel},
  title = {On the nonexistence of elements of {Kervaire} invariant one},
  journal = {Annals of Mathematics},
  volume = {184},
  number = {1},
  year = {2016},
  pages = {1--262},
  doi = {10.4007/annals.2016.184.1.1}
}

@article{HSWXHeightFour,
  author = {Michael A. Hill and XiaoLin Danny Shi and Guozhen Wang and Zhouli Xu},
  title = {The Slice Spectral Sequence of a {$C_4$}-Equivariant Height-4 {Lubin--Tate} Theory},
  journal = {Memoirs of the American Mathematical Society},
  volume = {288},
  number = {1429},
  year = {2023}
}

@misc{ChanVogeli,
  author = {David Chan and Chase Vogeli},
  title = {The {Galois}-equivariant {$K$}-theory of finite fields},
  year = {2024},
  eprint = {2406.19481},
  archivePrefix = {arXiv},
  primaryClass = {math.KT},
  note = {\href{https://arxiv.org/abs/2406.19481}{arXiv:2406.19481}}
}

@article{HellerOrmsby,
  author = {Jeremiah Heller and Kyle Ormsby},
  title = {{Galois} equivariance and stable motivic homotopy theory},
  journal = {Transactions of the American Mathematical Society},
  volume = {368},
  number = {11},
  year = {2016},
  pages = {8047--8077},
  doi = {10.1090/tran6647}
}

@article{HahnShiRealOrientations,
  author = {Jeremy Hahn and XiaoLin Danny Shi},
  title = {Real orientations of {Lubin--Tate} spectra},
  journal = {Inventiones Mathematicae},
  volume = {221},
  year = {2020},
  pages = {731--776},
  doi = {10.1007/s00222-020-00960-z}
}

@article{HellerVoineaguOstvaer,
  author  = {Heller, Jeremiah and Voineagu, Mircea
             and {\O}stv{\ae}r, Paul Arne},
  title   = {Topological comparison theorems for
             {Bredon} motivic cohomology},
  journal = {Transactions of the American Mathematical Society},
  volume  = {371},
  number  = {4},
  pages   = {2875--2921},
  year    = {2019},
  doi     = {10.1090/tran/7553}
}

@misc{Georgakopoulos,
  author        = {Nick Georgakopoulos},
  title         = {The {$RO(C_4)$} integral homology of a point},
  year          = {2019},
  eprint        = {1912.06758},
  archivePrefix = {arXiv},
  primaryClass  = {math.AT},
  note          = {\href{https://arxiv.org/abs/1912.06758}{arXiv:1912.06758}},
  url           = {https://arxiv.org/abs/1912.06758}
}

@misc{BasuGhosh,
  author        = {Samik Basu and Surojit Ghosh},
  title         = {Non-trivial extensions in equivariant cohomology with constant coefficients},
  year          = {2021},
  eprint        = {2108.12763},
  archivePrefix = {arXiv},
  primaryClass  = {math.AT},
  note          = {\href{https://arxiv.org/abs/2108.12763}{arXiv:2108.12763}},
  url           = {https://arxiv.org/abs/2108.12763}
}

@incollection{Lewis1988,
  author    = {L. Gaunce Lewis, Jr.},
  title     = {The {$RO(G)$}-graded equivariant ordinary cohomology of complex projective spaces with linear {$\mathbb Z/p$} actions},
  booktitle = {Algebraic Topology and Transformation Groups (G\"ottingen, 1987)},
  series    = {Lecture Notes in Mathematics},
  volume    = {1361},
  pages     = {53--122},
  publisher = {Springer},
  year      = {1988}
}

@book{FerlandLewis,
  author    = {Kevin K. Ferland and L. Gaunce Lewis, Jr.},
  title     = {The {$RO(G)$}-Graded Equivariant Ordinary Homology of {$G$}-Cell Complexes with Even-Dimensional Cells for {$G=\mathbb Z/p$}},
  series    = {Memoirs of the American Mathematical Society},
 note = {No. 794},
  publisher = {American Mathematical Society},
  year      = {2004},
  doi       = {10.1090/memo/0794}
}

@article{BasuGhoshSquarefree,
  author        = {Samik Basu and Surojit Ghosh},
  title         = {Equivariant cohomology for cyclic groups of square-free order},
  journal       = {Research in the Mathematical Sciences},
  volume        = {11},
  number        = {29},
  year          = {2024},
  doi           = {10.1007/s40687-024-00443-0},
  eprint        = {2006.09669},
  archivePrefix = {arXiv},
  primaryClass  = {math.AT}
}

@misc{AngeltveitPicard,
  author = {Vigleik Angeltveit},
  title = {The Picard group of the category of
    $C_n$-equivariant stable homotopy theory},
  year = {2021},
  eprint = {2109.04625},
  archivePrefix = {arXiv},
  primaryClass = {math.AT}
}

@book{GreenleesMayTate,
  author    = {J. P. C. Greenlees and J. P. May},
  title     = {Generalized {Tate} Cohomology},
  series    = {Memoirs of the American Mathematical Society},
  volume    = {113},
  note      = {No. 543},
  publisher = {American Mathematical Society},
  address   = {Providence, RI},
  year      = {1995}
}

@misc{StricklandTambara,
  author        = {Neil P. Strickland},
  title         = {{Tambara} functors},
  year          = {2012},
  eprint        = {1205.2516},
  archivePrefix = {arXiv},
  primaryClass  = {math.CT},
  note          = {arXiv:1205.2516},
  url           = {https://arxiv.org/abs/1205.2516}
}

@misc{BHKKS25,
  author = {Viktor Burghardt and Po Hu and Daniel Kriz
            and Igor Kriz and Petr Somberg},
  title  = {Chain-level models of equivariant topological
            {Hochschild} homology of semiperfect and smooth
            algebras over finite fields},
  year   = {2025},
  note   = {To appear in Homology, Homotopy and Applications},
  url    = {https://krizigor.github.io/chain25083.pdf}
}

@article{Bohme19,
  author        = {Benjamin B{\"o}hme},
  title         = {Multiplicativity of the idempotent splittings
                   of the {Burnside} ring and the {$G$}-sphere spectrum},
  journal       = {Advances in Mathematics},
  volume        = {347},
  pages         = {904--939},
  year          = {2019},
  doi           = {10.1016/j.aim.2019.03.008},
  eprint        = {1802.01938},
  archivePrefix = {arXiv},
  primaryClass  = {math.AT}
}

@article{ABGHLM18,
  author        = {Vigleik Angeltveit and Andrew J. Blumberg
                   and Teena Gerhardt and Michael A. Hill
                   and Tyler Lawson and Michael A. Mandell},
  title         = {Topological cyclic homology via the norm},
  journal       = {Documenta Mathematica},
  volume        = {23},
  pages         = {2101--2163},
  year          = {2018},
  doi           = {10.4171/DM/671},
  eprint        = {1401.5001},
  archivePrefix = {arXiv},
  primaryClass  = {math.KT}
}

@misc{BlumbergMandell24,
  author        = {Andrew J. Blumberg and Michael A. Mandell},
  title         = {The homotopy theory of cyclotomic spectra,
                   10 years later},
  year          = {2024},
  eprint        = {2412.15928},
  archivePrefix = {arXiv},
  primaryClass  = {math.AT},
  note          = {Theorem numbering refers to
                   \href{https://arxiv.org/abs/2412.15928}
                   {arXiv:2412.15928}},
  url           = {https://arxiv.org/abs/2412.15928}
}
\end{document}